\documentclass[a4paper,leqno,10t]{article}

\usepackage[off]{auto-pst-pdf}
\usepackage{amsmath,amsfonts,amssymb}
\usepackage{fancyhdr} 
\usepackage{epstopdf}
\usepackage{latexsym}
\usepackage{mathrsfs,amssymb} 
\usepackage[intlimits]{empheq} 
\usepackage[arrow, matrix, curve]{xy}
\usepackage{amsmath}
\usepackage{amssymb}
\usepackage{authblk}
\usepackage{empheq}
\usepackage{mathrsfs}
\usepackage[noindentafter]{titlesec}
\usepackage{blindtext}
\usepackage[a4paper]{geometry}
\usepackage[latin1]{inputenc}
\usepackage[T1]{fontenc}
\usepackage{amsthm}
\usepackage{dsfont}
\usepackage{pgf}
\usepackage{graphicx}
\usepackage[nottoc]{tocbibind}
\usepackage[linktocpage=true,colorlinks=false]{hyperref}
\usepackage{xcolor}
\usepackage{pfnote}
\usepackage{fnpos} 
\usepackage{color}
 \usepackage{fancyhdr,blindtext}
\fancypagestyle{plain}{%
\fancyhf{}}

\renewcommand{\sectionmark}[1]%
{\markboth{#1}{}}

\makeindex

\theoremstyle{definition} 

\newtheorem{ex}{\bfseries \upshape Example}[section]
\newtheorem{dfn}[ex]{\bfseries \upshape Definition}
\newtheorem{rem}[ex]{\bfseries \upshape Remark}

\theoremstyle{plain}
\newtheorem{prop}[ex]{\bfseries \upshape Proposition}
\newtheorem{lem}[ex]{\bfseries \upshape Lemma}
\newtheorem{theo}{\bfseries \upshape Theorem}
\newtheorem{cor}{\bfseries \upshape Corollary}
\newenvironment{pr}{\begin{proof}[{\bf Proof}]}{\end{proof}}

\begin{document}
\title{{ \bf Weakly mixing diffeomorphisms preserving a measurable Riemannian metric are dense in $\mathcal{A}_{\alpha}\left(M\right)$ for arbitrary Liouvillean number $\alpha$}}
\author[1]{{\sc Roland Gunesch}}
\author[2]{{\sc Philipp Kunde}}
\affil[1]{University of Education Vorarlberg, Feldkirch, Austria,}
\affil[2]{Department of Mathematics, University of Hamburg, Hamburg, Germany}
\date{}

\maketitle

\begin{abstract} 
We show that on any smooth compact connected manifold of dimension $m\geq 2$ admitting a smooth non-trivial circle action $\mathcal{S} = \left\{S_t\right\}_{t \in \mathbb{R}}$, $S_{t+1}=S_t$, the set of weakly mixing $C^{\infty}$-diffeomorphisms which preserve both a smooth volume $\nu$ and a measurable Riemannian metric is dense in $\mathcal{A}_{\alpha} \left(M\right)= \overline{ \left\{h \circ S_{\alpha} \circ h^{-1} : h \in \text{Diff}^{\infty}\left(M, \nu\right) \right\}}^{C^{\infty}}$ for every Liouvillean number $\alpha$. The proof is based on a quantitative version of the Anosov-Katok-method with explicitly constructed conjugation maps and partitions.
\end{abstract}





\section{Introduction}
To begin, recall that a dynamical system $\left(X,T,\nu\right)$ is ergodic if and only if every measurable complex-valued function $h$ on $\left(X, \nu\right)$ which is invariant (i.e. such that $h\left(Tx\right) = h\left(x\right)$ for every $x \in X$) must necessarily be constant. We define $\left(X,T,\nu\right)$ to be weakly mixing if it satisfies the stronger condition that there is no non-constant measurable complex valued function $h$ on $\left(X, \nu\right)$ such that $h\left(Tx\right) = \lambda \cdot h\left(x\right)$ for some $\lambda \in \mathbb{C}$. Equivalently there is  an increasing sequence $\left(m_n\right)_{n \in \mathbb{N}}$ of natural numbers such that $\lim_{n \rightarrow \infty} \left| \nu\left(B \cap T^{-m_n}\left(A\right)\right) - \nu\left(A\right) \cdot \nu\left(B\right) \right| = 0$ for every pair of measurable sets $A,B \subseteq X$ (see \cite{Skl} or \cite[Theorem 5.1]{AK}). We call a circle action $\left\{S_t\right\}_{t \in \mathbb{R}}$ on a manifold $M$ non-trivial if there exists $t \in \mathbb{R}$ and $x\in M$ with $S_t(x)\neq x$; in other words, not all orbits are fixed points (even though some may be).   \\
Until 1970 it was an open question if there exists an ergodic area-preserving smooth diffeomorphism on the disc $\mathbb{D}^2$. This problem was solved by the so-called ``approximation by conjugation''-method developed by D. Anosov and A. Katok in \cite{AK}. In fact, on every smooth compact connected manifold $M$ of dimension $m\geq 2$ admitting a non-trivial circle action $\mathcal{S} = \left\{S_t\right\}_{t \in \mathbb{S}^1}$ preserving a smooth volume $\nu$ this method enables the construction of smooth diffeomorphisms with specific ergodic properties (e.\,g.~weakly mixing ones in \cite[section 5]{AK}) or non-standard smooth realizations of measure-preserving systems (e.\,g.~\cite[section 6]{AK} and \cite{FSW}). These diffeomorphisms are constructed as limits of conjugates $f_n = H_n \circ S_{\alpha_{n+1}} \circ H^{-1}_n$, where $\alpha_{n+1} = \alpha_n + \frac{1}{k_n \cdot l_n \cdot q^2_n} \in \mathbb{Q}$, where $H_n = H_{n-1} \circ h_n$ and where $h_n$ are measure-preserving diffeomorphisms satisfying $S_{\frac{1}{q_n}} \circ h_n = h_n \circ S_{\frac{1}{q_n}}$. In each step the conjugation map $h_n$ and the parameter $k_n$ are chosen such that the diffeomorphism $f_n$ imitates the desired property with a certain precision. In a final step of the construction, the parameter $l_n$ is chosen large enough to guarantee closeness of $f_{n}$ to $f_{n-1}$ in the $C^{\infty}$-topology, and so the convergence of the sequence $\left(f_n\right)_{n \in \mathbb{N}}$ to a limit diffeomorphism is provided. It is even possible to keep this limit diffeomorphism within any given $C^{\infty}$-neighbourhood of the initial element $S_{\alpha_1}$ or, by applying a fixed diffeomorphism $g$ first, of $g\circ S_{\alpha_1} \circ g^{-1}$. So the construction can be carried out in a neighbourhood of any diffeomorphism conjugate to an element of the action. Thus, $\mathcal{A}\left(M\right) = \overline{\left\{h \circ S_t \circ h^{-1} \ : \ t \in \mathbb{S}^1, h \in \text{Diff}^{\infty}\left(M, \nu \right)\right\}}^{C^{\infty}}$ is a natural space for the produced diffeomorphisms. Moreover, we will consider the restricted space $\mathcal{A}_{\alpha}\left(M\right) = \overline{\left\{h \circ S_{\alpha} \circ h^{-1} \ : h \in \text{Diff}^{\infty}\left(M, \nu \right)\right\}}^{C^{\infty}}$ for $\alpha \in \mathbb{S}^1$. \\
In the following let $M$ be a smooth compact connected manifold of dimension $m \geq 2$ admitting a non-trivial circle action $\mathcal{S} = \left\{S_t \right\}_{t \in \mathbb{R}}$, $S_{t+1} = S_t$. Note that any such action possesses a smooth invariant volume: Every smooth manifold carries a Riemannian metric and hence a smooth Riemannian volume form $\hat{\nu}$. Any smooth volume form is given by $f \cdot \hat{\nu}$, where $f$ is a positive scalar function. If $\bar{f}$ is the fiberwise average of $f$, then $\bar{f} \cdot \hat{\nu}$ is a smooth volume form which is invariant under $\mathcal{S}$. In case of a manifold with boundary by a smooth diffeomorphism we mean infinitely differentiable in the interior and such that all the derivatives can be extended to the boundary continuously. \\
In their influential paper \cite{AK} Anosov and Katok proved amongst others that in $\mathcal{A} \left(M\right)$ the set of weakly mixing diffeomorphisms is generic (i.\,e. it is a dense $G_{\delta}$-set) in the $C^{\infty}\left(M\right)$-topology. For this they used the ``approximation by conjugation''-method. In \cite{GK} the conjugation maps are constructed more explicitly such that they can be equipped with the additional structure of being locally very close to an isometry, thus showing that there exists a weakly mixing smooth diffeomorphism preserving a smooth measure and a measurable Riemannian metric on any manifold with non-trivial circle action. Actually, it follows from the respective proofs that both results are true in $\mathcal{A}_{\alpha}\left(M\right)$ for a $G_{\delta}$-set of $\alpha \in \mathbb{R}$. However, both proofs do not give a full description of the set of $\alpha \in \mathbb{R}$ for which the particular result holds in $\mathcal{A}_{\alpha}\left(M\right)$. Such an investigation is started in \cite{FS}: B. Fayad and M. Saprykina showed in case of dimension $2$ that if $\alpha \in \mathbb{S}^1$ is Liouville, the set of weakly mixing diffeomorphisms is generic in the $C^{\infty}\left(M\right)$-topology in $\mathcal{A}_{\alpha}\left(M\right)$. Here an irrational number $\alpha$ is called Liouville if and only if for every $C \in \mathbb{R}_{>0}$ and for every $n \in \mathbb{N}$ there are infinitely many pairs of coprime integers $p,q$ such that $\left| \alpha - \frac{p}{q}\right| < \frac{C}{q^n}$. \\
In this article we prove the following theorem generalizing the results of \cite{GK} as well as \cite{FS}:
\begin{theo} \label{theo:main}
Let $M$ be a smooth compact and connected manifold of dimension $m\geq2$ with a non-trivial circle action $\mathcal{S} = \left\{S_t\right\}_{t \in \mathbb{R}}$, $S_{t+1}=S_t$. For any $\mathcal{S}$-invariant smooth volume $\nu$ the following is true: If $\alpha \in \mathbb{R}$ is Liouville, then the set of volume-preserving diffeomorphisms, that are weakly mixing and preserve a measurable Riemannian metric, is dense in the $C^{\infty}$-topology in $\mathcal{A}_{\alpha}\left(M\right)$.
\end{theo}
See \cite[section 3]{GK} for a comprehensive consideration of IM-diffeomorphisms (i.\,e. diffeomorphisms preserving an absolutely continuous probability measure and a measurable Riemannian metric) and IM-group actions. In particular, the existence of a measurable invariant metric for a diffeomorphism is equivalent to the existence of an invariant measure for the projectivized derivative extension which is absolutely continuous in the fibers. It is a natural question to ask about the ergodic properties of the derivative extension with respect to such a measure. While in our construction the projectivized derivative extension is as non-ergodic as possible (in fact, the derivative cocycle is cohomologous to the identity), it is work in progress to realize ergodic behaviour. This would provide the only known examples of measure-preserving diffeomorphisms whose differential is ergodic with respect to a smooth measure in the projectivization of the tangent bundle. Recently, it has been proven that for every $\rho >0$ and $m \geq 2$ there exists a weakly mixing real-analytic diffeomorphism $f \in \text{Diff}^{\omega}_{\rho} \left(\mathbb{T}^m, \mu \right)$ preserving a measurable Riemannian metric (\cite{Kana}).\\
We want to point out that Theorem \ref{theo:main} is in some sense the best we can obtain:  
\begin{itemize}
	\item By \cite[corollary 1.4]{FS}, whose proof uses Herman's last geometric result (\cite{FKr}), we have the following dichotomy in case of $M = \mathbb{S}^1 \times \left[0,1\right]$: A number $\alpha \in \mathbb{R} \setminus \mathbb{Q}$ is Diophantine if and only if there is no ergodic diffeomorphism of $M$ whose rotation number (on at least one of the boundaries) is equal to $\alpha$. Since weakly mixing diffeomorphisms are ergodic, there cannot be a weakly mixing $f \in \mathcal{A}_{\alpha}\left(\mathbb{S}^1 \times \left[0,1\right]\right)$ for $\alpha \in \mathbb{R} \setminus \mathbb{Q}$ Diophantine.
	\item By a result of A. Furman (appendix to \cite{GK}) a weakly mixing diffeomorphism cannot preserve a Riemannian metric with $L^{2}$-distortion (i.e. both the norm and its inverse are $L^2$-functions). Moreover, it is conjectured that a weakly mixing diffeomorphism cannot preserve a Riemannian metric with $L^{1}$-distortion (see \cite[Conjecture 3.7.]{GK}).
\end{itemize}

Using the standard techniques to prove genericity of the weak mixing-property and Theorem \ref{theo:main} we conclude in subsection \ref{subsection:first}:
\begin{cor} \label{cor:A}
Let $M$ be a smooth compact and connected manifold of dimension $m\geq2$ with a non-trivial circle action $\mathcal{S} = \left\{S_t\right\}_{t \in \mathbb{R}}$, $S_{t+1}=S_t$, preserving a smooth volume $\nu$. If $\alpha \in \mathbb{R}$ is Liouville, the set of volume-preserving weakly mixing diffeomorphisms is a dense $G_{\delta}$-set in the $C^{\infty}$-topology in $\mathcal{A}_{\alpha}\left(M\right)$.
\end{cor}
Hence, we obtain the result of \cite{FS} in arbitrary dimension at least $2$.

\section{Preliminaries}

\subsection{Definitions and notations}
In this chapter we want to introduce advantageous definitions and notations. Initially we discuss topologies on the space of smooth diffeomorphisms on the manifold $M = \mathbb{S}^1 \times \left[0,1\right]^{m-1}$. Note that for diffeomorphisms $f = \left(f_1,...,f_m\right): \mathbb{S}^1 \times \left[0,1\right]^{m-1} \rightarrow \mathbb{S}^1 \times \left[0,1\right]^{m-1}$ the coordinate function $f_1$ understood as a map $\mathbb{R}\times \left[0,1\right]^{m-1} \rightarrow \mathbb{R}$ has to satisfy the condition $f_1\left(\theta +n, r_1,...,r_{m-1}\right) = f_1 \left(\theta, r_1,...,r_{m-1}\right) + l$ for $n \in \mathbb{Z}$, where either $l =n$ or $l= -n$. Moreover, for $i \in \left\{2,...,m\right\}$ the coordinate function $f_i$ has to be $\mathbb{Z}$-periodic in the first component, i.e. $f_i \left( \theta + n, r_1,...,r_{m-1}\right) = f_i \left(\theta,r_1,...,r_{m-1}\right)$ for every $n \in \mathbb{Z}$. \\
To define explicit metrics on Diff$^k\left(\mathbb{S}^1 \times \left[0,1\right]^{m-1}\right)$ and in the following, the subsequent notations will be useful:
\begin{dfn}
\begin{enumerate}
	\item For a sufficiently differentiable function $f: \mathbb{R}^m \rightarrow \mathbb{R}$ and a multi-index $\vec{a} = \left(a_1,...,a_m\right) \in \mathbb{N}^m_0$
\begin{equation*}
D_{\vec{a}}f := \frac{\partial^{\left|\vec{a}\right|}}{\partial x_1^{a_1}...\partial x_m^{a_m}} f,
\end{equation*}
where $\left|\vec{a}\right| = \sum^{m}_{i=1} a_i$ is the order of $\vec{a}$.
  \item For a continuous function $F: \left(0,1\right)^m \rightarrow \mathbb{R}$
\begin{equation*}
\left\|F\right\|_0 := \sup_{z \in \left(0,1\right)^m} \left|F\left(z\right)\right|.
\end{equation*}
\end{enumerate}
\end{dfn}

Diffeomorphisms on $\mathbb{S}^1 \times \left[0,1\right]^{m-1}$ can be regarded as maps from $\left[0,1\right]^m$ to $\mathbb{R}^m$. In this spirit the expressions $\left\|f_i\right\|_0$ as well as  $\left\|D_{\vec{a}}f_i\right\|_0$ for any multi-index $\vec{a}$ with $\left| \vec{a}\right|\leq k$ have to be understood for $f=\left(f_1,...,f_m\right) \in$ Diff$^k\left(\mathbb{S}^1 \times \left[0,1\right]^{m-1}\right)$. Since such a diffeomorphism is a continuous map on the compact manifold and every partial derivative can be extended continuously to the boundary, all these expressions are finite. Thus the subsequent definition makes sense:
\begin{dfn}
\begin{enumerate}
	\item For $f,g \in$ Diff$^k\left(\mathbb{S}^1 \times \left[0,1\right]^{m-1}\right)$ with coordinate functions $f_i$ resp. $g_i$ we define
\begin{equation*}
\tilde{d}_0\left(f,g\right) = \max_{i=1,..,m} \left\{ \inf_{p \in \mathbb{Z}} \left\| \left(f - g\right)_i + p\right\|_0\right\}
\end{equation*}
as well as
\begin{equation*}
\tilde{d}_k\left(f,g\right) = \max \left\{ \tilde{d}_0\left(f,g\right), \left\|D_{\vec{a}}\left(f-g\right)_i\right\|_0 \ : \ i=1,...,m \ , \ 1\leq \left|\vec{a}\right| \leq k \right\}.
\end{equation*}
  \item Using the definitions from 1. we define for $f,g \in$ Diff$^k\left(\mathbb{S}^1 \times \left[0,1\right]^{m-1}\right)$:
  \begin{equation*}
  d_k\left(f,g\right) = \max \left\{ \tilde{d}_k\left(f,g\right) \ , \ \tilde{d}_k\left(f^{-1},g^{-1}\right)\right\}.
  \end{equation*}
\end{enumerate}
\end{dfn}

Obviously $d_k$ describes a metric on Diff$^k\left(\mathbb{S}^1 \times \left[0,1\right]^{m-1}\right)$ measuring the distance between the diffeomorphisms as well as their inverses. As in the case of a general compact manifold the following definition connects to it:

\begin{dfn}
\begin{enumerate}
	\item A sequence of Diff$^{\infty}\left(\mathbb{S}^1 \times \left[0,1\right]^{m-1}\right)$-diffeomorphisms is called convergent in Diff$^{\infty}\left(\mathbb{S}^1 \times \left[0,1\right]^{m-1}\right)$ if it converges in Diff$^k\left(\mathbb{S}^1 \times \left[0,1\right]^{m-1}\right)$ for every $k \in \mathbb{N}$.
	\item On Diff$^{\infty}\left(\mathbb{S}^1 \times \left[0,1\right]^{m-1}\right)$ we declare the following metric
	\begin{equation*}
	d_{\infty}\left(f,g\right) = \sum^{\infty}_{k=1} \frac{d_k\left(f,g\right)}{2^k \cdot \left(1 + d_k\left(f,g\right)\right)}.
	\end{equation*}
\end{enumerate}
\end{dfn}

It is a general fact that Diff$^{\infty}\left(\mathbb{S}^1 \times \left[0,1\right]^{m-1}\right)$ is a complete metric space with respect to this metric $d_{\infty}$. \\
Again considering diffeomorphisms on $\mathbb{S}^1 \times \left[0,1\right]^{m-1}$ as maps from $\left[0,1\right]^m$ to $\mathbb{R}^m$ we add the adjacent notation:
\begin{dfn}
Let $f \in$ Diff$^k\left(\mathbb{S}^1 \times \left[0,1\right]^{m-1}\right)$ with coordinate functions $f_i$ be given. Then
\begin{equation*}
\left\| Df \right\|_0 := \max_{i,j \in \left\{1,...,m\right\}} \left\| D_j f_i \right\|_0
\end{equation*}
and
\begin{equation*}
||| f ||| _k := \max \left\{ \left\|D_{\vec{a}} f_i \right\|_0 , \left\|D_{\vec{a}} \left(f^{-1}_{i}\right)\right\|_0 \ : \ i = 1,...,m, \ \vec{a} \text{ multi-index with } 0\leq \left| \vec{a}\right| \leq k \right\}.
\end{equation*}
\end{dfn}

\begin{rem}
By the above-mentioned observations for every multi-index $\vec{a}$ with $\left|\vec{a}\right|\geq 1$ and every $i \in \left\{1,...,m\right\}$ the derivative $D_{\vec{a}} h_i$ is $\mathbb{Z}$-periodic in the first variable. Since in case of a diffeomorphism $g= \left(g_1,...,g_m\right)$ on $\mathbb{S}^1 \times \left[0,1\right]^{m-1}$ regarded as a map $\left[0,1\right]^m \rightarrow \mathbb{R}^m$ the coordinate functions $g_j$ for $j \in \left\{2,...,m\right\}$ satisfy $g_j\left(\left[0,1\right]^m\right) \subseteq \left[0,1\right]$, it holds:
\begin{equation*}
\sup_{z \in \left(0,1\right)^m} \left| \left(D_{\vec{a}} h_i \right) \left(g\left(z\right)\right)\right| \leq |||h|||_{\left|\vec{a}\right|}.
\end{equation*}
\end{rem}

Furthermore, we introduce the notion of a partial partition of a compact manifold $M$, which is a pairwise disjoint countable collection of measurable subsets of $M$.
\begin{dfn}
\begin{itemize}
	\item A sequence of partial partitions $\nu_n$ converges to the decomposition into points if and only if for a given measurable set $A$ and for every $n \in \mathbb{N}$ there exists a measurable set $A_n$, which is a union of elements of $\nu_n$, such that $\lim_{n \rightarrow \infty} \mu \left( A \Delta A_n \right) = 0$. We often denote this by $\nu_n \rightarrow \varepsilon$.
	\item A partial partition $\nu$ is the image under a diffeomorphism $F:M \rightarrow M$ of a partial partition $\eta$  if and only if $\nu = \left\{ F \left(I\right) \: : \: I \in \eta \right\}$. We write this as $\nu = F\left(\eta\right)$.
\end{itemize}
\end{dfn}

\subsection{First steps of the proof} \label{subsection:first}
First of all we show how constructions on $\mathbb{S}^1 \times \left[0,1\right]^{m-1}$ can be transfered to a general compact connected smooth manifold M with a non-trivial circle action $\mathcal{S} = \left\{S_{t}\right\}_{t \in \mathbb{R}}$, $S_{t+1} = S_t$. By \cite[Proposition 2.1.]{AK}, we can assume that $1$ is the smallest positive number satisfying $S_t = \text{id}$. Hence, we can assume $\mathcal{S}$ to be effective. We denote the set of fixed points of $\mathcal{S}$ by $F$ and for $q \in \mathbb{N}$ $F_q$ is the set of fixed points of the map $S_{\frac{1}{q}}$. \\
On the other hand, we consider $\mathbb{S}^1 \times \left[0,1\right]^{m-1}$ with Lebesgue measure $\mu$. Furthermore, let $\mathcal{R} = \left\{R_{\alpha} \right\}_{\alpha \in \mathbb{S}^1}$ be the standard action of $\mathbb{S}^1$ on $\mathbb{S}^1 \times \left[0,1\right]^{m-1}$, where the map $R_{\alpha}$ is given by $R_{\alpha}\left(\theta, r_1,...,r_{m-1}\right) = \left(\theta + \alpha, r_1,...,r_{m-1}\right)$. Hereby, we can formulate the following result (see \cite[Proposition 1]{FSW}):
\begin{prop}
Let $M$ be a $m$-dimensional smooth, compact and connected manifold admitting an effective circle action $\mathcal{S} = \left\{S_{t}\right\}_{t \in \mathbb{R}}$, $S_{t+1} = S_t$, preserving a smooth volume $\nu$. Let $B \coloneqq \partial M \cup F \cup \left( \bigcup_{q\geq 1} F_q \right)$. There exists a continuous surjective map $G: \mathbb{S}^1 \times \left[0,1\right]^{m-1} \rightarrow M$ with the following properties:
\begin{enumerate}
	\item The restriction of $G$ to $\mathbb{S}^{1} \times \left(0,1\right)^{m-1}$ is a $C^{\infty}$-diffeomorphic embedding.
	\item $\nu\left(G\left(\partial\left(\mathbb{S}^{1} \times \left[0,1\right]^{m-1}\right)\right)\right) = 0$
	\item $G\left(\partial\left(\mathbb{S}^{1} \times \left[0,1\right]^{m-1}\right)\right) \supseteq B$
	\item $G_{*}\left(\mu\right) = \nu$
	\item $\mathcal{S} \circ G = G \circ \mathcal{R}$
\end{enumerate}
\end{prop}

By the same reasoning as in \cite[section 2.2.]{FSW}, this proposition allows us to carry a construction from $\left(\mathbb{S}^{1} \times \left[0,1\right]^{m-1}, \mathcal{R}, \mu\right)$ to the general case $\left(M, \mathcal{S}, \nu\right)$: \\
Suppose $f: \mathbb{S}^{1} \times \left[0,1\right]^{m-1} \rightarrow \mathbb{S}^{1} \times \left[0,1\right]^{m-1}$ is a weakly mixing diffeomorphism sufficiently close to $R_{\alpha}$ in the $C^{\infty}$-topology with $f$-invariant measurable Riemannian metric $\omega$ obtained by $f = \lim_{n \rightarrow \infty} f_n$ with $f_n = H_n \circ R_{\alpha_{n+1}} \circ H^{-1}_{n}$, where $f_n = R_{\alpha_{n+1}}$ in a neighbourhood of the boundary (in Proposition \ref{prop:satz2} we will see that these conditions can be satisfied in the constructions of this article). Then we define a sequence of diffeomorphisms:
\begin{equation*}
\tilde{f_n}: M \rightarrow M \ \ \ \ \ \ \tilde{f}_n\left(x\right) = \begin{cases}
G \circ f_n \circ G^{-1}\left(x\right)& \textit{if $x \in G\left(\mathbb{S}^1 \times \left(0,1\right)^{m-1}\right)$} \\
S_{\alpha_{n+1}}\left(x\right)& \textit{if $x \in G\left(\partial\left(\mathbb{S}^1 \times \left(0,1\right)^{m-1}\right)\right)$}
\end{cases}
\end{equation*}
Constituted in \cite[section 5.1.]{FK}, this sequence is convergent in the $C^{\infty}$-topology to the diffeomorphism 
\begin{equation*}
\tilde{f}: M \rightarrow M \ \ \ \ \ \ \tilde{f}\left(x\right) = \begin{cases}
G \circ f \circ G^{-1}\left(x\right)& \textit{if $x \in G\left(\mathbb{S}^1 \times \left(0,1\right)^{m-1}\right)$} \\
S_{\alpha}\left(x\right)& \textit{if $x \in G\left(\partial\left(\mathbb{S}^1 \times \left(0,1\right)^{m-1}\right)\right)$}
\end{cases}
\end{equation*}
provided the closeness from $f$ to $R_{\alpha}$ in the $C^{\infty}$-topology. \\
We observe that $f$ and $\tilde{f}$ are measure-theoretically isomorphic. Then $\tilde{f}$ is weakly mixing because the weak mixing-property is invariant under isomorphisms. \\
Moreover, we want to show how we can construct a $\tilde{f}$-invariant measurable Riemannian metric $\tilde{\omega}$ out of the $f$-invariant metric $\omega$. Since $\tilde{\omega}$ only needs to be a measurable metric and $\nu\left(G\left(\partial\left(\mathbb{S}^{1} \times \left[0,1\right]^{m-1}\right)\right)\right) = 0$, we only have to construct it on $G\left(\mathbb{S}^{1} \times \left(0,1\right)^{m-1}\right)$. Using the diffeomorphic embedding $G$ we consider $\tilde{\omega} |_{G\left(\mathbb{S}^{1} \times \left(0,1\right)^{m-1}\right)} \coloneqq \left(G^{-1}\right)^{\ast} \omega |_{G\left(\mathbb{S}^{1} \times \left(0,1\right)^{m-1}\right)}$ and show that it is $\tilde{f}$-invariant: On $G\left(\mathbb{S}^{1} \times \left(0,1\right)^{m-1}\right)$ we have $\tilde{f}=G \circ f \circ G^{-1}$ and thus we can compute:
\begin{equation*}
\tilde{f}^{\ast} \tilde{\omega} = \left(G \circ f \circ G^{-1}\right)^{\ast} \left(\left(G^{-1}\right)^{\ast} \omega \right) = \left(G^{-1}\right)^{\ast} \circ f^{\ast} \circ G^{\ast} \circ \left(G^{-1}\right)^{\ast} \omega = \left(G^{-1}\right)^{\ast} \circ f^{\ast} \omega = \left(G^{-1}\right)^{\ast} \omega = \tilde{\omega}
\end{equation*}
Altogether the construction done in the case of $\left(\mathbb{S}^{1} \times \left[0,1\right]^{m-1}, \mathcal{R}, \mu\right)$ is transfered to $\left(M, \mathcal{S}, \nu\right)$. Hence, it suffices to consider constructions on $M = \mathbb{S}^{1} \times \left[0,1\right]^{m-1}$ with circle action $\mathcal{R}$ subsequently. In this case we will prove the following result:
\begin{prop} \label{prop:satz2}
For every Liouvillean number $\alpha$ there is a sequence $\left(\alpha_n\right)_{n \in \mathbb{N}}$ of rational numbers $\alpha_n = \frac{p_n}{q_n}$ satisfying $\lim_{n\rightarrow \infty} \left| \alpha-\alpha_n\right|=0$ monotonically, and there are sequences $\left(g_n\right)_{n \in \mathbb{N}}$, $\left(\phi_n\right)_{n \in \mathbb{N}}$ of measure-preserving diffeomorphisms satisfying $g_n\circ R_{\frac{1}{q_n}}=R_{\frac{1}{q_n}}\circ g_n$ as well as $\phi_n \circ R_{\frac{1}{q_n}} =  R_{\frac{1}{q_n}} \circ \phi_n$ such that the diffeomorphisms $f_n = H_n \circ R_{\alpha_{n+1}} \circ H^{-1}_{n}$ with $H_n := h_1 \circ h_2 \circ ... \circ h_n$, where $h_n := g_n \circ \phi_n$, coincide with $R_{\alpha_{n+1}}$ in a neighbourhood of the boundary, converge in the Diff$^{\infty}\left(M\right)$-topology, and the diffeomorphism $f = \lim_{n \rightarrow \infty} f_n$ is weakly mixing, has an invariant measurable Riemannian metric, and satisfies $f \in \mathcal{A}_{\alpha}\left(M\right)$. \\
Furthermore, for every $\varepsilon > 0$ the parameters in the construction can be chosen in such a way that $d_{\infty}\left(f, R_{\alpha}\right) < \varepsilon$.
\end{prop}
By this Proposition weakly mixing diffeomorphisms preserving a measurable Riemannian metric are dense in $\mathcal{A}_{\alpha}\left(M\right)$: \\
Because of $\mathcal{A}_{\alpha} \left(M\right) = \overline{\left\{ h \circ R_{\alpha} \circ h^{-1} : h \in \text{Diff}^{\infty}\left(M, \mu\right)\right\}}^{C^{\infty}}$ it is enough to show that for every diffeomorphism $h\in\text{Diff}^{\infty}\left(M,\mu\right)$ and every $\epsilon > 0$ there is a weakly mixing diffeomorphism $\tilde{f}$ preserving a measurable Riemannian metric such that $d_{\infty}\left(\tilde{f}, h \circ R_{\alpha} \circ h^{-1} \right) < \epsilon$. For this purpose, let $h \in \text{Diff}^{\infty}\left(M, \mu\right)$ and $\epsilon > 0$ be arbitrary. By \cite[p. 3]{Omori} and \cite[Theorem 43.1.]{KM}, Diff$^{\infty}\left(M\right)$ is a Lie group. In particular, the conjugating map $g \mapsto h \circ g \circ h^{-1}$ is continuous with respect to the metric $d_{\infty}$. Continuity in the point $R_{\alpha}$ yields the existence of $\delta >0$, such that $d_{\infty}\left(g, R_{\alpha}\right)<\delta$ implies $d_{\infty}\left(h \circ g \circ h^{-1}, h \circ R_{\alpha} \circ h^{-1}\right)< \epsilon$. By Proposition \ref{prop:satz2} we can find a weakly mixing diffeomorphism $f$ with $f$-invariant measurable Riemannian metric $\omega$ and $d_{\infty}(f,R_{\alpha}) < \delta$. Hence $\tilde{f} := h \circ f \circ h^{-1}$ satisfies $d_{\infty}\left(\tilde{f}, h \circ R_{\alpha} \circ h^{-1} \right) < \epsilon$. Note that $\tilde{f}$ is weakly mixing and preserves the measurable Riemannian metric $\tilde{\omega} \coloneqq \left(h^{-1}\right)^{\ast}\omega$. \\
Hence, Theorem \ref{theo:main} is deduced from Proposition \ref{prop:satz2}.

\begin{rem}
Moreover we can show that the set of weakly mixing diffeomorphisms is generic in $\mathcal{A}_{\alpha}\left(M\right)$ (i.e. it is a dense $G_{\delta}$-set) using Proposition \ref{prop:satz2} and the same technique as in \cite{Ha}, section \textit{Category}. \\
Using Proposition \ref{prop:satz2} we can show that the set of weakly mixing diffeomorphisms is generic in $\mathcal{A}_{\alpha}\left(M\right)$ (i.e. it is a dense $G_{\delta}$-set). Thereby, we consider a countable dense set $\left\{\varphi_n\right\}_{n \in \mathbb{N}}$ in $L^2\left(M, \mu\right)$, which is a separable space, and define the sets:
\begin{equation*}
O\left(i,j,k,n\right) = \left\{ T \in \mathcal{A}_{\alpha}\left(M\right) \ : \ \left| \left(U^{n}_{T} \varphi_i , \varphi_j\right) - \left(\varphi_i,1\right) \cdot \left(1, \varphi_j\right) \right| < \frac{1}{k}\right\}
\end{equation*}
Since $\left(U_T \varphi, \psi\right)$ depends continuously on $T$, each $O\left(i,j,k,n\right)$ is open. Hence,
\begin{equation*}
K \coloneqq \bigcap_{i \in \mathbb{N}} \bigcap_{j \in \mathbb{N}} \bigcap_{k \in \mathbb{N}} \bigcup_{n \in \mathbb{N}} O\left(i,j,k,n\right)
\end{equation*}
is a $G_{\delta}$-set. \\
By another equivalent characterisation a measure-preserving transformation $T$ is weakly mixing if and only if for every $\varphi,\psi \in L^2\left(M, \mu\right)$ there is a sequence $\left(m_n\right)_{n \in \mathbb{N}}$ of density one such that $\lim_{n \rightarrow \infty} \left(U^{m_n}_{T} \varphi , \psi \right) = \left( \varphi,1\right)\cdot\left(1,\psi\right)$. Thus, every weakly mixing diffeomorphism is contained in $K$. On the other hand, we show that a transformation, that is not weakly mixing, does not belong to $K$: If $T$ is not weakly mixing, $U_T$ has a non-trivial eigenfunction. W.l.o.g. we can assume the existence of $f \in L^2\left(M, \mu\right)$ and $c \in \mathbb{C}$ of absolute value 1 satisfying $U_T f = c \cdot f$, $\left\|f\right\|_{L^2} = 1$ and $\left(1,f\right)=0$. Since $\left\{\varphi_n\right\}_{n \in \mathbb{N}}$ is dense in $L^2\left(M, \mu\right)$, there is an index $i$ such that $\left\|f-\varphi_i\right\|_{L^2} < 0.1$. Obviously $\left\|\varphi_i\right\|_{L^2} \leq \left\|f\right\|_{L^2} + \left\|f-\varphi_i\right\|_{L^2} < 1.1$ and $\left| \left(U^{n}_{T} f , f\right) - \left(f,1\right) \cdot \left(1, f\right) \right| = \left|\left(c^n \cdot f,f\right)\right|= \left|c^n\right| \cdot \left\|f\right\|^2_{L^2} = 1$. Consequently we can estimate:
\begin{align*}
1 & = \left| \left(U^{n}_{T} f , f\right) - \left(f,1\right) \cdot \left(1, f\right) \right| \\
& \leq \left| \left(U^{n}_{T} f , f\right) - \left(U^{n}_{T}f, \varphi_i\right) \right|+\left| \left(U^{n}_{T} f, \varphi_i\right) - \left(U^{n}_{T} \varphi_i,\varphi_i\right)\right|+\left| \left(U^{n}_{T} \varphi_i , \varphi_i\right) - \left(\varphi_i,1\right) \cdot \left(1, \varphi_i\right) \right| \\
& \ \ \ \ + \left| \left(\varphi_i,1\right) \cdot \left(1 , \varphi_i\right) - \left(\varphi_i,1\right) \cdot \left(1, f\right) \right| + \left| \left(\varphi_i , 1\right) \cdot \left(1,f\right)- \left(f,1\right) \cdot \left(1, f\right) \right| \\
& \leq \left|c\right|^n \cdot \left\|f\right\|_{L^2} \cdot \left\|f-\varphi_i\right\|_{L^2} + \left\|f-\varphi_i\right\|_{L^2} \cdot \left\|\varphi_i\right\|_{L^2} + \left| \left(U^{n}_{T} \varphi_i, \varphi_i\right) - \left(\varphi_i,1\right) \cdot \left(1, \varphi_i\right) \right| \\
& \ \ \ \ + \left\|\varphi_i\right\|_{L^2} \cdot \left\|f-\varphi_i\right\|_{L^2} \\
& \leq 0.1 + 0.11 + \left| \left(U^{n}_{T} \varphi_i, \varphi_i\right) - \left(\varphi_i,1\right) \cdot \left(1, \varphi_i\right) \right| + 0.11 \\
& < \left| \left(U^{n}_{T} \varphi_i, \varphi_i\right) - \left(\varphi_i,1\right) \cdot \left(1, \varphi_i\right) \right| + 0.5
\end{align*}
Thus $\left| \left(U^{n}_{T} \varphi_i, \varphi_i\right) - \left(\varphi_i,1\right) \cdot \left(1, \varphi_i\right) \right|$ has to be larger than $\frac{1}{2}$. Hence $T$ does not belong to $O\left(i,i,2,n\right)$ for any value of $n$ and accordingly does not belong to $K$.
So $K$ coincides with the set of weakly mixing diffeomorphisms in $\mathcal{A}_{\alpha}\left(M\right)$. By the observations above we know that this set is dense. In conclusion the set of weakly mixing diffeomorphisms is a dense $G_{\delta}$-set in $\mathcal{A}_{\alpha}\left(M\right)$. Thus Corollary \ref{cor:A} is proven.
\end{rem}

\subsection{Outline of the proof}
The constructions are based on the ``approximation by conjugation''-method developed by D.V. Anosov and A. Katok in \cite{AK}. As indicated in the introduction, one constructs successively a sequence of measure preserving diffeomorphisms $f_n= H_n \circ R_{\alpha_{n+1}} \circ H^{-1}_{n}$, where the conjugation maps $H_n = h_1 \circ ... \circ h_n$ and the rational numbers $\alpha_n = \frac{p_n}{q_n}$ are chosen in such a way that the functions $f_n$ converge to a diffeomorphism $f$ with the desired properties. \\
First of all we will define two sequences of partial partitions, which converge to the decomposition into points in each case. The first type of partial partition, called $\eta_n$, will satisfy the requirements in the proof of the weak mixing-property. On the partition elements of the even more refined second type, called $\zeta_n$, the conjugation map $h_n$ will act as an isometry, and this will enable us to construct an invariant measurable Riemannian metric. Afterwards we will construct these conjugating diffeomorphisms $h_n = g_n \circ \phi_n$, which are composed of two step-by-step defined smooth measure-preserving diffeomorphisms. In this construction the map $g_n$ should introduce shear in the $\theta$-direction as in \cite{FS}. So $\tilde{g}_{\left[n q^{\sigma}_{n}\right]} \left(\theta,r_1,..., r_{m-1}\right) = \left(\theta + \left[n \cdot q^{\sigma}_{n}\right] \cdot r_1, r_1,..., r_{m-1}\right)$ might seem an obvious candidate. Unfortunately, that map is not an isometry. Therefore, the map $g_n$ will be constructed in such a way that $g_n$ is an isometry on the image under $\phi_n$ of any partition element $\check{I}_n \in \zeta_n$, and $g_n \left(\hat{I}_n\right) = \tilde{g}_{\left[n q^{\sigma}_{n}\right]}\left(\hat{I}_n\right)$ as well as $g_n \left( \Phi_n \left(\hat{I}_n\right) \right) = \tilde{g}_{\left[n q^{\sigma}_{n}\right]} \left( \Phi_n \left(\hat{I}_n\right) \right)$ for every $\hat{I}_n \in \eta_n$, where $\Phi_n = \phi_n \circ R^{m_n}_{\alpha_{n+1}} \circ \phi^{-1}_{n}$ with a specific sequence $\left(m_n\right)_{n\in \mathbb{N}}$  of natural numbers (see section \ref{section:distri}) is important in the proof of the weak mixing property. Likewise the conjugation map $\phi_n$ will be built such that it acts on the elements of $\zeta_n$ as an isometry and on the elements of $\eta_n$ in such a way that it satisfies the requirements of the desired criterion for weak mixing. This criterion is established in section \ref{section:crit}. It is similar to the criterion in \cite{FS} but modified in many places because of the new conjugation map $g_n$ and the new type of partitions. The construction presented here combines the advantages of shearing maps and local isometries, and it combines local maps in such a way that the derivatives of the resultant conjugation maps can be suitably bounded. Unfortunately, this requires a fairly elaborate and slightly technical construction. \\
In section \ref{section:conv} we will show convergence of the sequence $\left(f_n\right)_{n \in \mathbb{N}}$ in $\mathcal{A}_{\alpha}\left(M\right)$ for a given Liouville number $\alpha$ by the same approach as in \cite{FS}. To do so, we have to estimate the norms $||| H_n |||_k$ very carefully. Furthermore, we will see at the end of section \ref{section:conv} that the criterion for weak mixing applies to the obtained diffeomorphism $f = \lim_{n\rightarrow \infty}f_n$. Finally, we will construct the desired $f$-invariant measurable Riemannian metric in section \ref{section:metric}.

\section{Explicit constructions} \label{section:constr}

\subsection{Sequences of partial partitions}
In this subsection we define the two announced sequences of partial partitions $\left(\eta_n\right)_{n \in \mathbb{N}}$ and $\left(\zeta_n\right)_{n \in \mathbb{N}}$ of $M = \mathbb{S}^{1} \times \left[0,1\right]^{m-1}$.

\subsubsection{Partial partition $\eta_n$} \label{subsubsection:eta}
\begin{rem}
For convenience we will use the notation $\prod^{m}_{i=2} \left[ a_i, b_i\right]$ for $\left[a_2, b_2\right] \times ... \times \left[a_m, b_m\right]$.
\end{rem}

Initially, $\eta_n$ will be constructed on the fundamental sector $\left[0, \frac{1}{q_n}\right] \times \left[0,1\right]^{m-1}$. For this purpose we divide the fundamental sector into $n$ sections: 
\begin{itemize}
	\item In case of $k \in \mathbb{N}$ and $2\leq k \leq n-1$ on $\left[\frac{k-1}{n \cdot q_n}, \frac{k}{n \cdot q_n} \right] \times \left[0,1\right]^{m-1}$ the partial partition $\eta_n$ consists of all multidimensional intervals of the following form:
	
\begin{equation*}
\begin{split}
& \Bigg{[} \frac{k-1}{n \cdot q_n} + \frac{j^{(1)}_{1}}{n \cdot q^{2}_n}+...+ \frac{j^{\left((m-1) \cdot \frac{\left(k+1\right) \cdot k}{2}\right)}_{1}}{n \cdot q^{1+(m-1) \cdot \frac{\left(k+1\right) \cdot k}{2}}_n} + \frac{1}{10 \cdot n^5 \cdot q^{1+(m-1) \cdot \frac{\left(k+1\right) \cdot k}{2}}_n}, \\
& \quad \frac{k-1}{n \cdot q_n}+ \frac{j^{(1)}_{1}}{n \cdot q^{2}_n}+...+\frac{j^{\left((m-1) \cdot \frac{\left(k+1\right) \cdot k}{2}\right)}_{1} + 1}{n \cdot q^{1+(m-1) \cdot \frac{\left(k+1\right) \cdot k}{2}}_n} - \frac{1}{10 \cdot n^5 \cdot q^{1+(m-1) \cdot \frac{\left(k+1\right) \cdot k}{2}}_n} \Bigg{]} \\
\times & \prod^{m}_{i=2} \left[ \frac{j^{(1)}_i}{q_n} + ...+\frac{j^{(k+1)}_i}{q^{k+1}_n} +  \frac{1}{26 \cdot n^4 \cdot q^{k+1}_n}, \frac{j^{(1)}_i}{q_n} + ...+\frac{j^{(k+1)}_i +1}{q^{k+1}_n} -  \frac{1}{26 \cdot n^4 \cdot q^{k+1}_n}\right],
\end{split}
\end{equation*}
where $j^{(l)}_1 \in \mathbb{Z}$ und $\left\lceil \frac{q_n}{10n^4} \right\rceil \leq j^{(l)}_1 \leq q_n - \left\lceil \frac{q_n}{10n^4} \right\rceil - 1$ for $l=1,...,(m-1) \cdot \frac{\left(k+1\right) \cdot k}{2}$ as well as $j^{(l)}_{i} \in \mathbb{Z}$ and $\left\lceil \frac{q_n}{10n^4} \right\rceil \leq j^{(l)}_{i} \leq q_n - \left\lceil \frac{q_n}{10n^4} \right\rceil - 1$ for $i=2,...,m$ and $l=1,...,k+1$.
   \item On $\left[0, \frac{1}{n \cdot q_n}\right] \times \left[0,1\right]^{m-1}$ as well as $\left[\frac{n-1}{n \cdot q_n}, \frac{1}{q_n} \right] \times \left[0,1\right]^{m-1}$ there are no elements of the partial partition $\eta_n$.
\end{itemize}

By applying the map $R_{l/q_n}$ with $l \in \mathbb{Z}$, this partial partition of $\left[0, \frac{1}{q_n}\right] \times \left[0,1\right]^{m-1}$ is extended to a partial partition of $\mathbb{S}^{1} \times \left[0,1\right]^{m-1}$.

\begin{rem}
By construction this sequence of partial partitions converges to the decomposition into points.
\end{rem}

\subsubsection{Partial partition $\zeta_n$} \label{subsubsection:zeta}
As in the previous case we will construct the partial partition $\zeta_n$ on the fundamental sector $\left[0, \frac{1}{q_n}\right] \times \left[0,1\right]^{m-1}$ initially and therefore divide this sector into $n$ sections: In case of $k \in \mathbb{N}$ and $1\leq k \leq n$ on $\left[\frac{k-1}{n \cdot q_n}, \frac{k}{n \cdot q_n} \right] \times \left[0,1\right]^{m-1}$ the partial partition $\zeta_n$ consists of all multidimensional intervals of the following form:

\begin{equation*}
\begin{split}
& \Bigg{[}\frac{k-1}{n \cdot q_n} + \frac{j^{(1)}_{1}}{n \cdot q^2_n}+...+ \frac{j^{\left((m-1) \cdot \frac{k \cdot \left(k+1\right)}{2}\right)}_{1}}{n \cdot q^{1+(m-1) \cdot \frac{k \cdot \left(k+1\right)}{2}}_n}+\frac{1}{n^5 \cdot q^{1+(m-1) \cdot \frac{k \cdot \left(k+1\right)}{2}}_n},\\
& \quad \frac{k-1}{n \cdot q_n} + \frac{j^{(1)}_{1}}{n \cdot q^2_n}+...+ \frac{j^{\left((m-1) \cdot \frac{k \cdot \left(k+1\right)}{2}\right)}_{1}+1}{n \cdot q^{1+(m-1) \cdot \frac{k \cdot \left(k+1\right)}{2}}_n}-\frac{1}{n^5 \cdot q^{1+(m-1) \cdot \frac{k \cdot \left(k+1\right)}{2}}_n} \Bigg{]} \\
\times & \Bigg{[} \frac{j^{(1)}_2}{q_n}+ ... +\frac{j^{\left((m-1) \cdot \frac{k \cdot \left(k+1\right)}{2}+1\right)}_2}{q^{1+(m-1) \cdot \frac{k \cdot \left(k+1\right)}{2}}_n}+ \frac{j^{\left((m-1) \cdot \frac{k \cdot \left(k+1\right)}{2}+2\right)}_2 }{8 n^5 \cdot q^{1+(m-1) \cdot \frac{k \cdot \left(k+1\right)}{2}}_n \cdot \left[n q^{\sigma}_{n}\right]} +  \frac{1}{8  n^9 \cdot q^{1+(m-1) \cdot \frac{k \cdot \left(k+1\right)}{2}}_n \cdot \left[ n q^{\sigma}_{n}\right]}, \\
& \quad \frac{j^{(1)}_2}{q_n} + ... +\frac{j^{\left((m-1) \cdot \frac{k \cdot \left(k+1\right)}{2}+1\right)}_2}{q^{1+(m-1) \cdot \frac{k \cdot \left(k+1\right)}{2}}_n}+ \frac{j^{\left((m-1) \cdot \frac{k \cdot \left(k+1\right)}{2}+2\right)}_2 + 1}{8 n^5 \cdot q^{1+(m-1) \cdot \frac{k \cdot \left(k+1\right)}{2}}_n \cdot \left[n q^{\sigma}_{n}\right]} -  \frac{1}{8  n^9 \cdot q^{1+(m-1) \cdot \frac{k \cdot \left(k+1\right)}{2}}_n \cdot \left[ n q^{\sigma}_{n}\right]}\Bigg{]} \\
\times & \prod^{m}_{i=3} \left[ \frac{j^{(1)}_i}{q_n} + ... + \frac{j^{(k)}_i}{q^{k}_n} + \frac{1}{n^4 \cdot q^{k}_n}, \frac{j^{(1)}_i}{q_n} +...+ \frac{j^{(k)}_i +1}{q^{k}_n}- \frac{1}{n^4 \cdot q^{k}_n}\right], 
\end{split}
\end{equation*}
where $j^{(l)}_i \in \mathbb{Z}$ and $\left\lceil \frac{q_n}{n^4} \right\rceil \leq j^{(l)}_i \leq q_n - \left\lceil \frac{q_n}{n^4} \right\rceil - 1$ for $i=3,...,m$ and $l=1,..,k$ as well as $j^{(l)}_{1} \in \mathbb{Z}$, $\left\lceil \frac{q_n}{n^4} \right\rceil \leq j^{(l)}_{1} \leq q_n - \left\lceil \frac{q_n}{n^4} \right\rceil - 1$ for $l=1,...,(m-1) \cdot \frac{k \cdot \left(k+1\right)}{2}$ as well as $j^{(l)}_2 \in \mathbb{Z}$ and $\left\lceil \frac{q_n}{n^4} \right\rceil \leq j^{(l)}_2 \leq q_n - \left\lceil \frac{q_n}{n^4} \right\rceil - 1$ for $l=1,...,(m-1) \cdot \frac{k \cdot \left(k+1\right)}{2}+1$ as well as $j^{\left((m-1) \cdot \frac{k \cdot \left(k+1\right)}{2}+2\right)}_2\in \mathbb{Z}$ and $8 \cdot n \cdot \left[n \cdot q^{\sigma}_{n}\right] \leq j^{\left((m-1) \cdot \frac{k\cdot \left(k+1\right)}{2}+2\right)}_2 \leq 8 \cdot n^5  \cdot \left[n \cdot q^{\sigma}_{n}\right] - 8 \cdot n \cdot \left[n \cdot q^{\sigma}_{n}\right] - 1$.

\begin{rem} \label{rem:überd}
For every $n\geq 3$ the partial partition $\zeta_n$ consists of disjoint sets, covers a set of measure at least $1-\frac{4 \cdot m}{n^2}$, and the sequence $\left(\zeta_n\right)_{n \in \mathbb{N}}$ converges to the decomposition into points. 
\end{rem}

\subsection{The conjugation map $g_n$} \label{subsection:g}
Let $\sigma \in (0,1)$. As mentioned in the sketch of the proof we aim for a smooth measure-preserving diffeomorphism $g_n$ which satisfies $g_n \left(\hat{I}_n\right) = \tilde{g}_{\left[n q^{\sigma}_{n}\right]} \left(\hat{I}_n\right)$ as well as $g_n \left( \Phi_n \left(\hat{I}_n\right) \right) = \tilde{g}_{\left[n q^{\sigma}_{n}\right]} \left( \Phi_n \left(\hat{I}_n\right) \right)$ for every $\hat{I}_n \in \eta_n$ and is an isometry on the image under $\phi_n$ of any partition element $\check{I}_n \in \zeta_n$. \\
Let $a,b \in \mathbb{Z}$ and $\varepsilon \in \left(0, \frac{1}{16}\right]$ such that $\frac{1}{\varepsilon} \in \mathbb{Z}$. Moreover, we consider $\delta >0$ such that $\frac{1}{\delta} \in \mathbb{Z}$ and $\frac{a\cdot b \cdot \delta}{\varepsilon} \in \mathbb{Z}$. We denote $\left[0,1\right]^2$ by $\Delta$ and $\left[\varepsilon, 1- \varepsilon \right]^2$ by $\Delta\left(\varepsilon\right)$.
\begin{lem} \label{lem:g}
For every $\varepsilon \in \left(0, \frac{1}{16}\right]$ there exists a smooth measure-preserving diffeomorphism $g_{\varepsilon}: \left[0,1\right]^2 \rightarrow \left\{ \left(x + \varepsilon \cdot y, y \right) \ : \ x,y \in \left[0,1\right]\right\}$ that is the identity on $\Delta\left(4 \varepsilon \right)$ and coincides with the map $\left(x,y\right) \mapsto \left(x + \varepsilon \cdot y, y\right)$ on $\Delta \setminus \Delta\left(\varepsilon\right)$.
\end{lem}

\begin{pr}
First of all let $\psi_{\varepsilon}:\mathbb{R}^2\to\mathbb{R}^2$ be a smooth diffeomorphism satisfying
\begin{equation*}
\psi_{\varepsilon}\left(x,y\right) = \begin{cases}
\left(x,y\right)& \textit{on $\mathbb{R}^2 \setminus \Delta\left(2\varepsilon\right)$} \\
\left(\frac{1}{2}+\frac{1}{5} \cdot \left(x-\frac{1}{2}\right), \frac{1}{2}+\frac{1}{5} \cdot \left(y-\frac{1}{2}\right)\right)& \textit{on $\Delta\left(4 \varepsilon\right)$}
\end{cases}
\end{equation*}
Furthermore, let $\tau_{\varepsilon}$ be a smooth diffeomorphism with the following properties
\begin{equation*}
\tau_{\varepsilon}\left(x,y\right) = \begin{cases}
\left(x + \varepsilon \cdot y,y \right)& \textit{on $\left\{\left(x-\frac{1}{2}\right)^2+\left(y-\frac{1}{2}\right)^2 \geq \left(\frac{5}{16}\right)^2\right\}$} \\
\left(x,y\right)& \textit{on $\left\{\left(x-\frac{1}{2}\right)^2+\left(y-\frac{1}{2}\right)^2 \leq \frac{1}{50}\right\}$}
\end{cases}
\end{equation*}
We define $\bar{g}_{\varepsilon} \coloneqq \psi^{-1}_{\varepsilon} \circ \tau_{\varepsilon} \circ \psi_{\varepsilon}$. Then the diffeomorphism $\bar{g}_{\varepsilon}$ coincides with the identity on $\Delta\left(4 \varepsilon\right)$ and with the map $\left(x,y\right) \mapsto \left(x + \varepsilon \cdot y, y\right)$ on $\mathbb{R}^2 \setminus \Delta\left(\varepsilon\right)$. From this we conclude that $\det \left(D\bar{g}_{\varepsilon}\right) > 0$. Moreover, $\bar{g}_{\varepsilon}$ is measure-preserving on $U_{\varepsilon} \coloneqq \left(\mathbb{R}^2 \setminus \Delta\left(\varepsilon\right)\right) \cup \Delta\left(4 \varepsilon \right)$. \\
With the aid of ``Moser's trick'' we want to construct a diffeomorphism $g_{\varepsilon}$ which is measure-preserving on the whole $\mathbb{R}^2$ and agrees with $\bar{g}_{\varepsilon}$ on $U_{\varepsilon}$. To do so, we consider the canonical volume form $\Omega_0$ on $\mathbb{R}^2$: $\Omega_0 = dx \wedge dy$; in other words, $\Omega_0 = d\omega_0$ using the 1-form $\omega_0 = \frac{1}{2} \cdot \left( x \cdot dy - y \cdot dx\right)$. Additionally we introduce the volume form $\Omega_1 \coloneqq \bar{g}^{*}_{\varepsilon} \Omega_0$. \\
At first we note that $\bar{g}_{\varepsilon}$ preserves the 1-form $\omega_0$ on $U_{\varepsilon}$: Clearly this holds on $\Delta\left(4 \varepsilon\right)$, where $\bar{g}_{\varepsilon}$ is the identity. On $\mathbb{R}^2 \setminus \Delta \left(\varepsilon \right)$ we have $D\bar{g}_{\varepsilon} \left(x,y\right) = \left(x + \varepsilon y, y\right)$, and thus we get
\begin{equation*}
\bar{g}^{*}_{\varepsilon} \omega_0 = \omega_0 \left(x+ \varepsilon \cdot y, y\right) = \frac{1}{2} \cdot \left(\left(x+ \varepsilon \cdot y\right) dy - y \cdot d\left(x+\varepsilon \cdot y\right)\right) = \frac{1}{2} \cdot \left(x \cdot dy - y \cdot dx\right) = \omega_0\left(x,y\right).
\end{equation*}
Furthermore, we introduce $\Omega' \coloneqq \Omega_1 - \Omega_0$. Since the exterior derivative commutes with the pull-back, it holds that $\Omega' = d\left(\bar{g}^{*}_{\varepsilon} \omega_0 - \omega_0\right)$. In addition we consider the volume form $\Omega_t \coloneqq \Omega_0 + t \cdot \Omega'$ and note that $\Omega_t$ is non-degenerate for $t\in[0,1]$. Thus, we get a uniquely defined vector field $X_t$ such that $\Omega_t\left(X_t, \cdot \right) = \left(\omega_0 - \bar{g}^{*}_{\varepsilon} \omega_0\right) \left(\cdot\right)$. Since $\Delta$ is a compact manifold, the non-autonomous differential equation $\frac{d}{dt} u(t) = X_t \left(u(t)\right)$ with initial values in $\Delta$ has a solution defined on $\mathbb{R}$. Hence, we get a one-parameter family of diffeomorphisms $\left\{\nu_t\right\}_{t \in \left[0,1\right]}$ on $\Delta$ satisfying $\dot{\nu}_t = X_t\left(\nu_t\right)$, $\nu_0 = \text{id}$. \\
Referring to \cite[Lemma 2.2]{Be}, it holds that
\begin{equation*}
\frac{d}{dt} \nu^{\ast}_t \Omega_t = d\left(\nu^{\ast}_t\left( i\left(X_t\right) \Omega_t \right)\right) + \nu^{\ast}_t\left(\frac{d}{dt} \Omega_t + i\left(X_t\right)d\Omega_t\right).
\end{equation*}
Because of $d\left(\nu^{\ast}_t\left( i\left(X_t\right) \Omega_t \right)\right) = \nu^{\ast}_t\left(d\left( i\left(X_t\right) \Omega_t \right)\right)$ and $d \Omega_t = d\left( d \omega_0 + t \cdot \left(d \left(\bar{g}^{\ast}_{\varepsilon} \omega_0\right) - d \omega_0\right)\right) = 0$ we compute:
\begin{align*}
\frac{d}{dt} \nu^{\ast}_t \Omega_t & = \nu^{\ast}_t\left(d\left( i\left(X_t\right) \Omega_t \right)\right) + \nu^{\ast}_t\left(\frac{d}{dt} \Omega_t\right) = \nu^{\ast}_t d \left( \Omega_t\left(X_t, \cdot \right) \right) + \nu^{\ast}_t \Omega' \\
& = \nu^{\ast}_t d\left(\omega_0 - \bar{g}^{\ast}_{\varepsilon}\omega_0\right) + \nu^{\ast}_t \Omega' = \nu^{\ast}_t \left( \Omega_0 - \Omega_1 \right) + \nu^{\ast}_t \left( \Omega_1 - \Omega_0 \right) = 0.
\end{align*}
Consequently $\nu^{\ast}_1 \Omega_1 = \nu^{\ast}_0 \Omega_0 = \Omega_0$ (using $\nu_0 = \text{id}$ in the last step). As we have seen, it holds that $\bar{g}^{*}_{\varepsilon} \omega_0 = \omega_0$ on $U_{\varepsilon}$. Therefore, on $U_{\varepsilon}$ it holds that $\Omega_t \left( X_t, \cdot \right) = 0$. Since $\Omega_t$ is non-degenerate, we conclude that $X_t = 0$ on $U_{\varepsilon}$ and hence $\nu_1 = \nu_0 = \text{id}$ on $U_{\varepsilon} \cap \Delta$. Now we can extend $\nu_1$ smoothly to $\mathbb{R}^2$ as the identity. \\
Denote $g_{\varepsilon} \coloneqq \bar{g}_{\varepsilon} \circ \nu_1$. Because of $\nu_1 = \text{id}$ on $U_{\varepsilon}$, the map $g_{\varepsilon}$ coincides with $\bar{g}_{\varepsilon}$ on $U_{\varepsilon}$ as announced. Furthermore we have
\begin{equation*}
g^{\ast}_{\varepsilon} \Omega_0 = \left(\bar{g}_{\varepsilon} \circ \nu_1 \right)^{\ast} \Omega_0 = \nu^{\ast}_1 \left( \bar{g}^{\ast}_{\varepsilon} \Omega_0\right)= \nu^{\ast}_1 \Omega_1 = \Omega_0.
\end{equation*}
Using the transformation formula we compute for an arbitrary measurable set $A \subseteq \mathbb{R}^2$:
\begin{equation*}
\mu\left( g_{\varepsilon}\left(A\right)\right) = \int_{g_{\varepsilon}\left(A\right)} \Omega_0 = \int_{A} \left|\det\left(Dg_{\varepsilon}\right)\right| \cdot \Omega_0.
\end{equation*}
We know $\det\left(D\nu_1\right) > 0$ (because $\nu_0 = \text{id}$ and all the maps $\nu_t$ are diffeomorphisms) as well as $\det\left(D \bar{g}_{\varepsilon}\right) > 0$, and thus $\left|\det\left(Dg_{\varepsilon}\right)\right| = \det\left(Dg_{\varepsilon}\right)$. Since $g^{\ast}_{\varepsilon} \Omega_0 = \left(\det\left(Dg_{\varepsilon}\right) \right) \cdot \Omega_0$ (compare with \cite[proposition 5.1.3.]{KH}) we finally conclude:
\begin{equation*}
\mu\left( g_{\varepsilon}\left(A\right)\right) = \int_{A} \left(\det\left(Dg_{\varepsilon}\right)\right) \cdot \Omega_0 = \int_{A} g^{\ast}_{\varepsilon} \Omega_0 = \int_{A} \Omega_0 = \mu\left(A\right).
\end{equation*}
Consequently $g_{\varepsilon}$ is a measure-preserving diffeomorphism on $\mathbb{R}^2$ satisfying the desired properties.
\end{pr}

Let $\tilde{g}_b: \mathbb{S}^1 \times \left[0,1\right]^{m-1} \rightarrow \mathbb{S}^1 \times \left[0,1\right]^{m-1}$ be the smooth measure-preserving diffeomorphism $\tilde{g}_b \left(\theta, r_1,...,r_{m-1}\right) = \left(\theta + b \cdot r_1, r_1,...,r_{m-1}\right)$ and denote $\left[0, \frac{1}{a}\right] \times \left[0, \frac{\varepsilon}{b \cdot a}\right] \times \left[\delta,1- \delta\right]^{m-2}$ by $\Delta_{a,b, \varepsilon, \delta}$. Using the map $D_{a,b,\varepsilon}: \mathbb{R}^m \rightarrow \mathbb{R}^m, \left(\theta, r_1,..., r_{m-1}\right) \mapsto \left(a \cdot \theta, \frac{b \cdot a}{\varepsilon} \cdot r_1, r_2,..., r_{m-1} \right)$ and $g_{\varepsilon}$ from Lemma \ref{lem:g} we define the measure-preserving diffeomorphism $g_{a,b, \varepsilon, \delta}: \Delta_{a,b, \varepsilon, \delta} \rightarrow \tilde{g}_b\left(\Delta_{a,b, \varepsilon, \delta}\right)$ by setting $g_{a,b, \varepsilon, \delta} = D^{-1}_{a,b, \varepsilon} \circ \left(g_{\varepsilon}, \text{id}_{\mathbb{R}^{m-2}}\right) \circ D_{a,b, \varepsilon}$. Using the fact that $\frac{a\cdot b \cdot \delta}{\varepsilon} \in \mathbb{Z}$ we extend it to a smooth diffeomorphism $g_{a,b, \varepsilon, \delta}: \left[0, \frac{1}{a}\right] \times \left[\delta,1-\delta\right]^{m-1} \rightarrow \tilde{g}_b \left( \left[0, \frac{1}{a}\right] \times \left[\delta,1-\delta\right]^{m-1} \right)$ by the description:
\begin{equation*}
g_{a,b, \varepsilon, \delta} \left(\theta , r_1 + l \cdot \frac{\varepsilon}{b \cdot a},r_2,..., r_{m-1} \right) = \left(l \cdot \frac{\varepsilon}{a}, l \cdot \frac{\varepsilon}{b \cdot a}, \vec{0}\right) + g_{a,b, \varepsilon, \delta} \left( \theta, r_1,...,r_{m-1}\right)
\end{equation*}
for $r_1 \in \left[0, \frac{\varepsilon}{b \cdot a}\right]$ and some $l \in \mathbb{Z}$ satisfying $\frac{b \cdot a \cdot \delta}{\varepsilon} \leq l \leq \frac{b \cdot a}{\varepsilon} - \frac{b \cdot a \cdot \delta}{\varepsilon} - 1$. Since this map coincides with the map $\tilde{g}_b$ in a neighbourhood of the boundary we can extend it to a map $g_{a,b,\varepsilon, \delta}: \left[0,\frac{1}{a}\right] \times \left[0,1\right]^{m-1} \rightarrow \tilde{g}_b\left(\left[0,\frac{1}{a}\right] \times \left[0,1\right]^{m-1}\right)$ by setting it equal to $\tilde{g}_b$ on $\left[0,\frac{1}{a}\right] \times \left(\left[0,1\right]^{m-1} \setminus \left[\delta,1-\delta\right]^{m-1}\right)$. \\

We initially construct the smooth measure-preserving diffeomorphism $g_n$ on the fundamental sector. For this, we divide the sector into $n$ sections: On $\left[\frac{k-1}{n \cdot q_n}, \frac{k}{n \cdot q_n}\right] \times \left[0,1\right]^{m-1}$ in case of $k \in \mathbb{Z}$ and $1 \leq k \leq n$: 
\begin{equation*}
g_n = g_{n \cdot q^{1+(m-1) \cdot \frac{\left(k+1\right) \cdot k}{2}}_n, \left[n \cdot q^{\sigma}_{n}\right], \frac{1}{8n^4}, \frac{1}{32n^4}}.
\end{equation*}
Since $g_n$ coincides with the map $\tilde{g}_{\left[n \cdot q^{\sigma}_{n}\right]}$ in a neighbourhood of the boundary of the different sections on the $\theta$-axis, this yields a smooth map, and we can extend it to a smooth measure-preserving diffeomorphism on $\mathbb{S}^1 \times \left[0,1\right]^{m-1}$ using the description $g_n \circ R_{\frac{l}{q_n}} = R_{\frac{l}{q_n}} \circ g_n$ for $l \in \mathbb{Z}$. Furthermore, we note that the subsequent constructions are done in such a way that $260n^4$ divides $q_n$ (see Lemma \ref{lem:conv}) and so the assumption $\frac{a\cdot b \cdot \delta}{\varepsilon}=\frac{a \cdot b}{4} \in \mathbb{Z}$ is satisfied. Indeed, this map $g_n$ satisfies the following desired property:

\begin{lem} \label{lem:outer}
For every element $\hat{I}_n \in \eta_n$ we have $g_n\left(\hat{I}_n\right) = \tilde{g}_{\left[n q^{\sigma}_{n}\right]} \left(\hat{I}_n\right)$.
\end{lem}

\begin{pr}
We consider a partition element $\hat{I}_{n,k} \in \eta_n$ on $\left[\frac{k-1}{n \cdot q_n}, \frac{k}{n \cdot q_n}\right] \times \left[0,1\right]^{m-1}$ in case of $k \in \mathbb{Z}$ and $2 \leq k \leq n-1$ and want to examine the effect of $g_n = g_{n \cdot q^{1+(m-1) \cdot \frac{\left(k+1\right) \cdot k}{2}}_n, \left[n \cdot q^{\sigma}_{n}\right], \frac{1}{8n^4}, \frac{1}{32n^4}}$ on it. \\
In the $r_1$-coordinate we use the fact that there is $u \in \mathbb{Z}$ such that 
\begin{equation*}
\frac{1}{26n^4 q^{k+1}_n} =u \cdot \frac{\varepsilon}{b \cdot a} = u \cdot \frac{1}{8 n^4 \cdot \left[nq^{\sigma}_n\right] \cdot nq^{1+(m-1) \cdot \frac{\left(k+1\right) \cdot k}{2}}_n}, 
\end{equation*}
where we use the fact that $260n^4$ divides $q_n$ (Lemma \ref{lem:conv}). Also, with respect to the $\theta$-coordinate the boundary of this element lies in the domain where $g_{a,b, \varepsilon, \delta} = \tilde{g}_{\left[nq^{\sigma}_n\right]}$ because $\frac{1}{10 \cdot n^4} < \varepsilon =\frac{1}{8 \cdot n^4}$.
\end{pr}

\subsection{The conjugation map $\phi_n$} \label{subsection:phi}

\begin{lem} \label{lem:varphi}
For every $\varepsilon \in \left(0, \frac{1}{4}\right)$ and every $i,j \in \left\{1,...,m\right\}$ there exists a smooth measure-preserving diffeomorphism $\varphi_{\varepsilon, i, j}$ on $\mathbb{R}^m$ which is the rotation in the $x_i - x_{j}$-plane by $\pi/2$ about the point $\left(\frac{1}{2},..., \frac{1}{2}\right) \in \mathbb{R}^m$ on $\left[2 \varepsilon, 1-2 \varepsilon\right]^{m}$ and coincides with the identity outside of $\left[\varepsilon, 1-\varepsilon\right]^m$.
\end{lem}

\begin{pr}
The proof is similar to the proof of Lemma \ref{lem:g}. (See also \cite[section 4.6]{GK} for a geometrical argument of the proof.)
\end{pr}

Furthermore, for $\lambda \in \mathbb{N}$ we define the maps $C_{\lambda} \left(x_1,x_2,...,x_m\right) = \left( \lambda \cdot x_1, x_2,..., x_m \right)$ and $D_{\lambda} \left(x_1,...,x_m \right) = \left( \lambda \cdot x_1, \lambda \cdot x_2,..., \lambda \cdot x_m \right)$. Let $\mu \in \mathbb{N}$, $\frac{1}{\delta} \in \mathbb{N}$ and assume $\frac{1}{\delta}$ divides $\mu$. We construct a diffeomorphism $\psi_{\mu, \delta, i, j, \varepsilon_2}$ in the following way:
\begin{itemize}
	\item Consider $\left[0, 1-2 \cdot \delta\right]^m$: Since $\frac{1}{\delta}$ divides $\mu$, we can divide $\left[0, 1-2 \cdot \delta\right]^m$ into cubes of side length $\frac{1}{\mu}$.
	\item Under the map $D_{\mu}$ any of these cubes of the form $\prod^{m}_{i=1} \left[\frac{l_i}{\mu}, \frac{l_i +1}{\mu}\right]$ with $l_i \in \mathbb{N}$ is mapped onto $\prod^{m}_{i=1} \left[ l_i, l_i +1\right]$.
	\item On $\left[0,1\right]^{m}$ we will use the diffeomorphism $\varphi^{-1}_{\varepsilon_2, i, j}$ constructed in Lemma \ref{lem:varphi} . Since this is the identity outside of $\Delta\left(\varepsilon_2\right)$, we can extend it to a diffeomorphism $\bar{\varphi}^{-1}_{\varepsilon_2, i, j}$ on $\mathbb{R}^m$ using the instruction $\bar{\varphi}^{-1}_{\varepsilon_2, i, j} \left(x_1 +k_1, x_2 + k_2,...,x_m + k_m\right) = \left(k_1,...,k_m\right) + \varphi^{-1}_{\varepsilon_2, i, j}\left(x_1, x_2,...,x_m\right)$, where $k_i \in \mathbb{Z}$ and $x_i \in \left[0,1\right]$.
	\item Now we define the smooth measure-preserving diffeomorphism 
	\begin{equation*}
	\tilde{\psi}_{\mu, \delta, i, j, \varepsilon_2} = D^{-1}_{\mu} \circ \bar{\varphi}^{-1}_{\varepsilon_2, i, j} \circ D_{\mu} \ \ \ :\ \ \ \left[0,1-2\delta\right]^{m} \rightarrow \left[0,1-2\delta\right]^{m}
	\end{equation*}
	\item With this we define 
\begin{align*}
& \psi_{\mu, \delta, i, j, \varepsilon_2} \left(x_1,...,x_m \right) = \\
& \begin{cases}
\left(\left[\tilde{\psi}_{\mu, \delta, i, j, \varepsilon_2}\left(x_1-\delta,...,x_m-\delta\right)\right]_1 + \delta,...,\left[\tilde{\psi}_{\mu, \delta, i, j, \varepsilon_2}\left(x_1-\delta,...,x_m-\delta\right)\right]_m + \delta \right)& \text{on $\left[\delta,1-\delta\right]^m$} \\
\left(x_1,...,x_m\right)& \text{otherwise}
\end{cases}
\end{align*}
	This is a smooth map because $\tilde{\psi}_{\mu, \delta, i, j, \varepsilon_2}$ is the identity in a neighbourhood of the boundary by the construction.
\end{itemize}

\begin{rem} \label{rem:W}
For every set $W = \prod^{m}_{i=1} \left[\frac{l_i}{\mu} + r_i, \frac{l_i +1}{\mu} - r_i\right]$ where $l_i \in \mathbb{Z}$ and $r_i \in \mathbb{R}$ satisfies $\left|r_i \cdot \mu \right| \leq \varepsilon_2$ we have $\psi_{\mu, \delta, i, j, \varepsilon_2} \left(W\right) = W$.
\end{rem}

Using these maps we build the following smooth measure-preserving diffeomorphism:
\begin{equation*}
\tilde{\phi}_{\lambda, \varepsilon, i,j, \mu, \delta, \varepsilon_2} : \left[0, \frac{1}{\lambda} \right] \times \left[0,1\right]^{m-1} \rightarrow \left[0, \frac{1}{\lambda} \right] \times \left[0,1\right]^{m-1}, \ \ \ \tilde{\phi}_{\lambda, \varepsilon, i,j, \mu, \delta, \varepsilon_2} = C^{-1}_{\lambda} \circ \psi_{\mu, \delta, i, j, \varepsilon_2} \circ \varphi_{\varepsilon, i, j} \circ C_{\lambda}
\end{equation*}
Afterwards, $\tilde{\phi}_{\lambda, \varepsilon, i,j, \mu, \delta, \varepsilon_2}$ is extended to a diffeomorphism on $\mathbb{S}^1 \times \left[0,1\right]^{m-1}$ by the description $\tilde{\phi}_{\lambda, \varepsilon, i,j, \mu, \delta, \varepsilon_2}\left(x_1 + \frac{1}{\lambda}, x_2,..., x_m\right)= \left(\frac{1}{\lambda}, 0, ..., 0\right) + \tilde{\phi}_{\lambda, \varepsilon, i,j, \mu, \delta, \varepsilon_2}\left(x_1, x_2,..., x_m\right)$.

For convenience we will use the notation $\phi^{(j)}_{\lambda, \mu} = \tilde{\phi}_{\lambda, \frac{1}{60n^4}, 1,j, \mu, \frac{1}{10n^4}, \frac{1}{22n^4}}$. With this we define the diffeomorphism $\phi_n$ on the fundamental sector: On $\left[\frac{k-1}{n \cdot q_n}, \frac{k}{n \cdot q_n}\right] \times \left[0,1\right]^{m-1}$ in case of $k \in \mathbb{N}$ and $1 \leq k \leq n$: 	
\begin{equation*}
\phi_n = \tilde{\phi}^{(m)}_{n \cdot q^{1+(m-1) \cdot \frac{k \cdot \left(k-1\right)}{2} + \left(m-2\right) \cdot k}_n, q^k_n} \circ \tilde{\phi}^{(m-1)}_{n \cdot q^{1+(m-1) \cdot \frac{k \cdot \left(k-1\right)}{2} + \left(m-3\right) \cdot k}_n, q^{k}_{n}} \circ ... \circ \tilde{\phi}^{(2)}_{n \cdot q^{1+(m-1) \cdot \frac{k \cdot \left(k-1\right)}{2}}_n, q^{k}_{n}}
\end{equation*}
This is a smooth map because $\phi_n$ coincides with the identity in a neighbourhood of the different sections. \\
Now we extend $\phi_n$ to a diffeomorphism on $\mathbb{S}^1 \times \left[0,1\right]^{m-1}$ using the description $\phi_n \circ R_{\frac{1}{q_n}} = R_{\frac{1}{q_n}} \circ \phi_n$.
\begin{figure}[hbtp] 
\begin{center}
\includegraphics{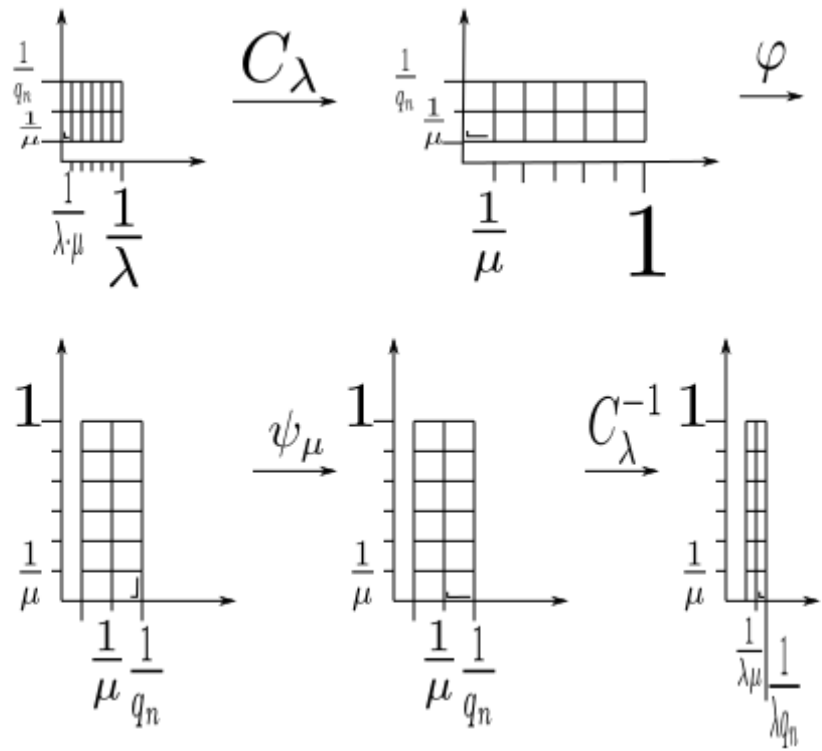}
\caption{Qualitative shape of the action of $\phi_n$ on a partition element $\hat{I} \in \eta_n$ and tangent vectors in case of dimension $m=2$.}
\end{center}
\end{figure}

\section{$\left(\gamma, \delta, \epsilon\right)$-distribution} \label{section:distri}
We introduce the central notion of the criterion for weak mixing deduced in the next section:
\begin{dfn}
Let $\Phi: M \rightarrow M$ be a diffeomorphism. We say $\Phi$ $\left(\gamma, \delta, \epsilon\right)$-distributes an element $\hat{I}$ of a partial partition, if the following properties are satisfied:
\begin{itemize}
	\item $\pi_{\vec{r}}\left(\Phi\left(\hat{I}\right)\right)$ is a $\left(m-1\right)$-dimensional interval $J$, i.e. $J = I_1 \times ... \times I_{m-1}$ with intervals $I_k \subseteq \left[0,1\right]$, and $1-\delta \leq \lambda\left(I_k\right) \leq 1$ for $k=1,...,m-1$. Here $\pi_{\vec{r}}$ denotes the projection on the $\left(r_1,...,r_{m-1}\right)$-coordinates (i.e., the last $m-1$ coordinates; the first one is the $\theta$-coordinate).
	\item $\Phi\left(\hat{I}\right)$ is contained in a set of the form $\left[c,c+\gamma\right] \times J$ for some $c \in \mathbb{S}^{1}$.
	\item For every $\left(m-1\right)$-dimensional interval $\tilde{J} \subseteq J$ it holds:
	\begin{equation*}
	\left| \frac{\mu\left(\hat{I}\cap \Phi^{-1}\left(\mathbb{S}^{1}\times \tilde{J}\right)\right)}{\mu\left(\hat{I}\right)} - \frac{\mu^{(m-1)}\left(\tilde{J}\right)}{\mu^{(m-1)}\left(J\right)} \right| \leq \epsilon \cdot \frac{\mu^{(m-1)}\left(\tilde{J}\right)}{\mu^{(m-1)}\left(J\right)},
	\end{equation*}
	where $\mu^{(m-1)}$ is the Lebesgue measure on $\left[0,1\right]^{m-1}$. 
\end{itemize}
\end{dfn}

\begin{rem}
Analogous to \cite{FS} we will call the third property ``almost uniform distribution'' of $\hat{I}$ in the $r_1,..,r_{m-1}$-coordinates. In the following we will often write it in the form of
\begin{equation*}
\left| \mu\left(\hat{I}\cap \Phi^{-1}\left(\mathbb{S}^{1}\times \tilde{J}\right)\right) \cdot \mu^{(m-1)}\left(J\right) - \mu\left(\hat{I}\right) \cdot \mu^{(m-1)}\left(\tilde{J}\right) \right| \leq \epsilon \cdot \mu \left(\hat{I} \right) \cdot \mu^{(m-1)}\left(\tilde{J}\right).
\end{equation*}
\end{rem}

In the next step we define the sequence of natural numbers $\left(m_n\right)_{n \in \mathbb{N}}$:
\begin{align*}
m_n & = \min \left\{ m \leq q_{n+1} \ \  : \ \ m \in \mathbb{N},\ \  \inf_{k \in \mathbb{Z}} \left| m \cdot \frac{p_{n+1}}{q_{n+1}} - \frac{1}{n \cdot q_n} + \frac{k}{q_n}\right| \leq \frac{260 \cdot (n+1)^4}{q_{n+1}}\right\} \\
& = \min \left\{ m \leq q_{n+1} \ \ : \ \ m \in \mathbb{N},\ \ \inf_{k \in \mathbb{Z}} \left| m \cdot \frac{q_n \cdot p_{n+1}}{q_{n+1}} - \frac{1}{n} + k \right| \leq \frac{260 \cdot (n+1)^4 \cdot q_n}{q_{n+1}}\right\}
\end{align*}

\begin{lem}
The set $\left\{ m \leq q_{n+1} \ \ : \ \ m \in \mathbb{N},\ \ \inf_{k \in \mathbb{Z}} \left| m \cdot \frac{q_n \cdot p_{n+1}}{q_{n+1}} - \frac{1}{n} + k \right| \leq \frac{260(n+1)^4 \cdot q_n}{q_{n+1}}\right\}$ is nonempty for every $n \in \mathbb{N}$, i.e., $m_n$ exists.
\end{lem}

\begin{pr}
In Lemma \ref{lem:conv} we will construct the sequence $\alpha_n = \frac{p_n}{q_n}$ in such a way that $q_n = 260n^4 \cdot \tilde{q}_n$ and $p_n = 260n^4 \cdot \tilde{p}_n$ with $\tilde{p}_n, \tilde{q}_n$ relatively prime. Therefore, the set $\left\{ j \cdot \frac{q_n \cdot p_{n+1}}{q_{n+1}} \ : \ j=1,...,q_{n+1} \right\}$ contains $\frac{q_{n+1}}{260(n+1)^4 \cdot \text{gcd}\left(q_n, \tilde{q}_{n+1}\right)}$ different equally distributed points on $\mathbb{S}^1$. Hence there are at least $\frac{q_{n+1}}{260(n+1)^4 \cdot q_n}$ different such points and so for every $x \in \mathbb{S}^1$ there is a $j \in \left\{1,...,q_{n+1} \right\}$ such that $\inf_{k \in \mathbb{Z}} \left| x - j \cdot \frac{q_n \cdot p_{n+1}}{q_{n+1}} + k \right| \leq \frac{260(n+1)^4 \cdot q_n}{q_{n+1}}$. In particular, this is true for $x=\frac{1}{n}$.
\end{pr}

\begin{rem} \label{rem:an}
We define
\begin{equation*}
a_n = \left(m_n \cdot \frac{p_{n+1}}{q_{n+1}} - \frac{1}{n \cdot q_n}\right) \text{ mod } \frac{1}{q_n}
\end{equation*}
By the above construction of $m_n$ it holds that $\left|a_n \right| \leq \frac{260 \cdot (n+1)^4}{q_{n+1}}$. In Lemma \ref{lem:conv} we will see that it is possible to choose $q_{n+1} \geq 64 \cdot 260 \cdot \left(n+1\right)^4 \cdot n^{11} \cdot q^{(m-1) \cdot n^2+3}_n$. Thus, we get:
\begin{equation*}
\left|a_n \right| \leq \frac{1}{64 \cdot n^{11} \cdot q^{(m-1) \cdot n^2+3}_n}.
\end{equation*}
\end{rem}

Our constructions are done in such a way that the following property is satisfied:

\begin{lem} \label{lem:distri}
The map $\Phi_n \coloneqq \phi_n \circ R^{m_n}_{\alpha_{n+1}} \circ \phi^{-1}_{n}$ with the conjugating maps $\phi_n$ defined in section \ref{subsection:phi} $\left(\frac{1}{n \cdot q^{m}_{n}}, \frac{1}{n^4}, \frac{1}{n}\right)$-distributes the elements of the partition $\eta_n$.
\end{lem}

\begin{pr}
We consider a partition element $\hat{I}_{n,k}$ on $\left[\frac{k-1}{n \cdot q_n}, \frac{k}{n \cdot q_n}\right] \times \left[0,1\right]^{m-1}$. When applying the map $\phi^{-1}_n$ we observe that this element is positioned in such a way that all the occuring maps $\varphi^{-1}_{\varepsilon,1,j}$ and $\varphi_{\varepsilon_2,1,j}$ act as the respective rotations. Then we compute $\phi^{-1}_n\left(\hat{I}_{n,k}\right)$:
\begin{equation*}
\begin{split}
& \Bigg{[} \frac{k-1}{n \cdot q_n} + \frac{j^{(1)}_{1}}{n \cdot q^{2}_n}+...+ \frac{j^{\left((m-1) \cdot \frac{\left(k-1\right) \cdot k}{2}\right)}_{1}}{n \cdot q^{(m-1) \cdot \frac{\left(k-1\right) \cdot k}{2}+1}_n} + \frac{j^{(1)}_2}{n \cdot q^{(m-1) \cdot \frac{\left(k-1\right) \cdot k}{2}+2}_n}+ ... + \frac{j^{(k)}_2}{n \cdot q^{(m-1) \cdot \frac{\left(k-1\right) \cdot k}{2}+k+1}_n}  \\
& \quad + \frac{j^{(1)}_3}{n \cdot q^{(m-1) \cdot \frac{\left(k-1\right) \cdot k}{2}+k+2}_n}+...+\frac{j^{(k)}_m}{n \cdot q^{(m-1) \cdot \frac{\left(k+1\right) \cdot k}{2}+1}_n}+\frac{1}{10 \cdot n^5 \cdot q^{(m-1) \cdot \frac{\left(k+1\right) \cdot k}{2}+1}_n}, \\
& \quad \frac{k-1}{n \cdot q_n}+ \frac{j^{(1)}_{1}}{n \cdot q^{2}_n}+...+\frac{j^{(k)}_m+1}{n \cdot q^{(m-1) \cdot \frac{\left(k+1\right) \cdot k}{2}+1}_n}-\frac{1}{10 \cdot n^5 \cdot q^{(m-1) \cdot \frac{\left(k+1\right) \cdot k}{2}+1}_n} \Bigg{]} \\
\times & \prod^{m}_{i=2} \Bigg{[} 1-\frac{j^{\left((m-1) \cdot \frac{\left(k-1\right) \cdot k}{2}+(i-2) \cdot k +1\right)}_{1}}{q_n}-...-\frac{j^{\left((m-1) \cdot \frac{\left(k-1\right) \cdot k}{2}+(i-1) \cdot k \right)}_{1}+1}{q^k_n}+\frac{j^{(k+1)}_i}{q^{k+1}_n} +  \frac{1}{26 \cdot n^4 \cdot q^{k+1}_n}, \\
& \quad 1-\frac{j^{\left((m-1) \cdot \frac{\left(k-1\right) \cdot k}{2}+(i-2) \cdot k +1\right)}_{1}}{q_n}-...-\frac{j^{\left((m-1) \cdot \frac{\left(k-1\right) \cdot k}{2}+(i-1) \cdot k \right)}_{1}+1}{q^k_n}+\frac{j^{(k+1)}_i+1}{q^{k+1}_n} -  \frac{1}{26 \cdot n^4 \cdot q^{k+1}_n} \Bigg{]}.
\end{split}
\end{equation*}
By our choice of the number $m_n$ the subsequent application of $R^{m_n}_{\alpha_{n+1}}$ yields modulo $\frac{1}{q_n}$:

\begin{equation*}
\begin{split}
& \Bigg{[} \frac{k}{n \cdot q_n} + \frac{j^{(1)}_{1}}{n \cdot q^{2}_n}+...+ \frac{j^{\left((m-1) \cdot \frac{\left(k-1\right) \cdot k}{2}\right)}_{1}}{n \cdot q^{(m-1) \cdot \frac{\left(k-1\right) \cdot k}{2}+1}_n} + \frac{j^{(1)}_2}{n \cdot q^{(m-1) \cdot \frac{\left(k-1\right) \cdot k}{2}+2}_n}+ ... + \frac{j^{(k)}_2}{n \cdot q^{(m-1) \cdot \frac{\left(k-1\right) \cdot k}{2}+k+1}_n}  \\
& \quad + \frac{j^{(1)}_3}{n \cdot q^{(m-1) \cdot \frac{\left(k-1\right) \cdot k}{2}+k+2}_n}+...+\frac{j^{(k)}_m}{n \cdot q^{(m-1) \cdot \frac{\left(k+1\right) \cdot k}{2}+1}_n}+\frac{1}{10 \cdot n^5 \cdot q^{(m-1) \cdot \frac{\left(k+1\right) \cdot k}{2}+1}_n} + a_n, \\
& \quad \frac{k}{n \cdot q_n}+ \frac{j^{(1)}_{1}}{n \cdot q^{2}_n}+...+\frac{j^{(k)}_m+1}{n \cdot q^{(m-1) \cdot \frac{\left(k+1\right) \cdot k}{2}+1}_n}-\frac{1}{10 \cdot n^5 \cdot q^{(m-1) \cdot \frac{\left(k+1\right) \cdot k}{2}+1}_n} + a_n \Bigg{]} \\
\times & \prod^{m}_{i=2} \Bigg{[} 1-\frac{j^{\left((m-1) \cdot \frac{\left(k-1\right) \cdot k}{2}+(i-2) \cdot k +1\right)}_{1}}{q_n}-...-\frac{j^{\left((m-1) \cdot \frac{\left(k-1\right) \cdot k}{2}+(i-1) \cdot k \right)}_{1}+1}{q^k_n} +\frac{j^{(k+1)}_i}{q^{k+1}_n} + \frac{1}{26 \cdot n^4 \cdot q^{k+1}_n}, \\
& \quad 1-\frac{j^{\left((m-1) \cdot \frac{\left(k-1\right) \cdot k}{2}+(i-2) \cdot k +1\right)}_{1}}{q_n}-...-\frac{j^{\left((m-1) \cdot \frac{\left(k-1\right) \cdot k}{2}+(i-1) \cdot k \right)}_{1}+1}{q^k_n} +\frac{j^{(k+1)}_i+1}{q^{k+1}_n}-  \frac{1}{26 \cdot n^4 \cdot q^{k+1}_n} \Bigg{]},
\end{split}
\end{equation*}
at which $a_n$ is the ``error term'' introduced in Remark \ref{rem:an}. Under $\varphi_{\frac{1}{60n^4}, 1,2} \circ C_{n \cdot q^{(m-1) \cdot \frac{\left(k+1\right) \cdot k}{2}+1}_n}$ this is mapped to

\begin{equation*}
\begin{split}
& \left(\frac{k}{n \cdot q_n} + \frac{j^{(1)}_{1}}{n \cdot q^{2}_n}+...+\frac{j^{(k)}_m}{n \cdot q^{(m-1) \cdot \frac{\left(k+1\right) \cdot k}{2}+1}_n}, \vec{0}\right) + \\
& \Bigg{[}\frac{j^{\left((m-1) \cdot \frac{\left(k-1\right) \cdot k}{2}+1\right)}_1}{q_n}+ ... +\frac{j^{\left((m-1) \cdot \frac{\left(k-1\right) \cdot k}{2}+k\right)}_1+1}{q^{k}_n}- \frac{j^{(k+1)}_2+1}{q^{k+1}_n}+ \frac{1}{26 \cdot n^4 \cdot q^{k+1}_n}, \\
& \quad \frac{j^{\left((m-1) \cdot \frac{\left(k-1\right) \cdot k}{2}+1\right)}_1}{q_n}+ ... +\frac{j^{\left((m-1) \cdot \frac{\left(k-1\right) \cdot k}{2}+k\right)}_1+1}{q^{k}_n}- \frac{j^{(k+1)}_2}{q^{k+1}_n}- \frac{1}{26 \cdot n^4 \cdot q^{k+1}_n}, \Bigg{]} \\
\times & \left[\frac{1}{10 \cdot n^4}+ n \cdot q^{(m-1) \cdot \frac{\left(k+1\right) \cdot k}{2}+1}_n \cdot a_n, 1- \frac{1}{10 \cdot n^4}+ n \cdot q^{(m-1) \cdot \frac{\left(k+1\right) \cdot k}{2}+1}_n \cdot a_n\right] \\
\times & \prod^{m}_{i=3} \Bigg{[} 1-\frac{j^{\left((m-1) \cdot \frac{\left(k-1\right) \cdot k}{2}+(i-2) \cdot k +1\right)}_{1}}{q_n}-...-\frac{j^{\left((m-1) \cdot \frac{\left(k-1\right) \cdot k}{2}+(i-1) \cdot k \right)}_{1}+1}{q^k_n} +\frac{j^{(k+1)}_i}{q^{k+1}_n} + \frac{1}{26 \cdot n^4 \cdot q^{k+1}_n}, \\
& \quad 1-\frac{j^{\left((m-1) \cdot \frac{\left(k-1\right) \cdot k}{2}+(i-2) \cdot k +1\right)}_{1}}{q_n}-...-\frac{j^{\left((m-1) \cdot \frac{\left(k-1\right) \cdot k}{2}+(i-1) \cdot k \right)}_{1}+1}{q^k_n} +\frac{j^{(k+1)}_i+1}{q^{k+1}_n}-  \frac{1}{26 \cdot n^4 \cdot q^{k+1}_n} \Bigg{]}
\end{split}
\end{equation*}
using the bound on $a_n$. With the aid of Remark \ref{rem:W}, the bound on $a_n$ from Remark \ref{rem:an} and the fact that $10n^4$ divides $q^{k+1}_n$ by Lemma \ref{lem:conv} we can compute the image of $\hat{I}_{n,k}$ under $\tilde{\phi}^{(2)}_{n \cdot q^{(m-1) \cdot \frac{\left(k+1\right) \cdot k}{2}+1}_n,q^{k+1}_n} \circ R^{m_n}_{\alpha_{n+1}} \circ \phi^{-1}_n$:

\begin{equation*}
\begin{split}
& \Bigg{[} \frac{k}{n \cdot q_n} + \frac{j^{(1)}_{1}}{n \cdot q^{2}_n}+...+\frac{j^{(k)}_m}{n \cdot q^{(m-1) \cdot \frac{\left(k+1\right) \cdot k}{2}+1}_n}+\frac{j^{\left((m-1) \cdot \frac{\left(k-1\right) \cdot k}{2}+1\right)}_1}{n \cdot q^{(m-1) \cdot \frac{\left(k+1\right) \cdot k}{2}+2}_n}+... \\
& \quad +\frac{j^{\left((m-1) \cdot \frac{\left(k-1\right) \cdot k}{2}+k\right)}_1+1}{n \cdot q^{(m-1) \cdot \frac{\left(k+1\right) \cdot k}{2}+k+1}_n}-\frac{j^{(k+1)}_2+1}{n \cdot q^{(m-1) \cdot \frac{\left(k+1\right) \cdot k}{2}+k+2}_n}+ \frac{1}{26 \cdot n^5 \cdot q^{(m-1) \cdot \frac{\left(k+1\right) \cdot k}{2}+k+2}_n} , \\
& \quad \frac{k}{n \cdot q_n}+ \frac{j^{(1)}_{1}}{n \cdot q^{2}_n}+...-\frac{j^{(k+1)}_2}{n \cdot q^{(m-1) \cdot \frac{\left(k+1\right) \cdot k}{2}+k+2}_n}- \frac{1}{26 \cdot n^5 \cdot q^{(m-1) \cdot \frac{\left(k+1\right) \cdot k}{2}+k+2}_n} \Bigg{]} \\
\times & \left[\frac{1}{10 \cdot n^4}+ n \cdot q^{(m-1) \cdot \frac{\left(k+1\right) \cdot k}{2}+1}_n \cdot a_n, 1- \frac{1}{10 \cdot n^4}+ n \cdot q^{(m-1) \cdot \frac{\left(k+1\right) \cdot k}{2}+1}_n \cdot a_n\right] \\
\times & \prod^{m}_{i=3} \Bigg{[} 1-\frac{j^{\left((m-1) \cdot \frac{\left(k-1\right) \cdot k}{2}+(i-2) \cdot k +1\right)}_{1}}{q_n}-...-\frac{j^{\left((m-1) \cdot \frac{\left(k-1\right) \cdot k}{2}+(i-1) \cdot k \right)}_{1}+1}{q^k_n} +\frac{j^{(k+1)}_i}{q^{k+1}_n} + \frac{1}{26 \cdot n^4 \cdot q^{k+1}_n}, \\
& \quad 1-\frac{j^{\left((m-1) \cdot \frac{\left(k-1\right) \cdot k}{2}+(i-2) \cdot k +1\right)}_{1}}{q_n}-...-\frac{j^{\left((m-1) \cdot \frac{\left(k-1\right) \cdot k}{2}+(i-1) \cdot k \right)}_{1}+1}{q^k_n} +\frac{j^{(k+1)}_i+1}{q^{k+1}_n}-  \frac{1}{26 \cdot n^4 \cdot q^{k+1}_n} \Bigg{]}.
\end{split}
\end{equation*}

Continuing in the same way we obtain that $\Phi_n\left(\hat{I}_{n,k}\right)$ is equal to

\begin{equation*}
\begin{split}
& \Bigg{[} \frac{k}{n \cdot q_n} + \frac{j^{(1)}_{1}}{n \cdot q^{2}_n}+...+ \frac{j^{\left((m-1) \cdot \frac{\left(k-1\right) \cdot k}{2}\right)}_{1}}{n \cdot q^{(m-1) \cdot \frac{\left(k-1\right) \cdot k}{2}+1}_n} + \frac{j^{(1)}_2}{n \cdot q^{(m-1) \cdot \frac{\left(k-1\right) \cdot k}{2}+2}_n}+ ...+ \frac{j^{(k)}_2}{n \cdot q^{(m-1) \cdot \frac{\left(k-1\right) \cdot k}{2}+k+1}_n}  \\
& \quad + \frac{j^{(1)}_3}{n \cdot q^{(m-1) \cdot \frac{\left(k-1\right) \cdot k}{2}+k+2}_n}+...+\frac{j^{(k)}_m}{n \cdot q^{(m-1) \cdot \frac{\left(k+1\right) \cdot k}{2}+1}_n} \\
& \quad +\frac{j^{\left((m-1) \cdot \frac{\left(k-1\right) \cdot k}{2}+1\right)}_1}{n \cdot q^{(m-1) \cdot \frac{\left(k+1\right) \cdot k}{2}+2}_n}+...+\frac{j^{\left((m-1) \cdot \frac{\left(k-1\right) \cdot k}{2}+k\right)}_1+1}{n \cdot q^{(m-1) \cdot \frac{\left(k+1\right) \cdot k}{2}+k+1}_n}-\frac{j^{(k+1)}_2+1}{n \cdot q^{(m-1) \cdot \frac{\left(k+1\right) \cdot k}{2}+k+2}_n} \\
& \quad  +\frac{j^{\left((m-1) \cdot \frac{\left(k-1\right) \cdot k}{2}+k+1\right)}_1}{n \cdot q^{(m-1) \cdot \frac{\left(k+1\right) \cdot k}{2}+k+3}_n}+...+\frac{j^{\left((m-1) \cdot \frac{\left(k-1\right) \cdot k}{2}+2k\right)}_1+1}{n \cdot q^{(m-1) \cdot \frac{\left(k+1\right) \cdot k}{2}+2k+2}_n}-\frac{j^{(k+1)}_3+1}{n \cdot q^{(m-1) \cdot \frac{\left(k+1\right) \cdot k}{2}+2k+3}_n}+... \\
& \quad +\frac{j^{\left((m-1) \cdot \frac{\left(k+1\right) \cdot k}{2}\right)}_1+1}{n \cdot q^{(m-1) \cdot \frac{\left(k+1\right) \cdot \left(k+2\right)}{2}}_n}-\frac{j^{(k+1)}_m+1}{n \cdot q^{(m-1) \cdot \frac{\left(k+1\right) \cdot \left(k+2\right)}{2}+1}_n}+\frac{1}{26 \cdot n^5 \cdot q^{(m-1) \cdot \frac{\left(k+1\right) \cdot \left(k+2\right)}{2}+1}_n} , \\
& \quad \frac{k}{n \cdot q_n}+ \frac{j^{(1)}_{1}}{n \cdot q^{2}_n}+...- \frac{j^{(k+1)}_m}{n \cdot q^{(m-1) \cdot \frac{\left(k+1\right) \cdot \left(k+2\right)}{2}+1}_n}-\frac{1}{26 \cdot n^5 \cdot q^{(m-1) \cdot \frac{\left(k+1\right) \cdot \left(k+2\right)}{2}+1}_n} \Bigg{]} \\
\times & \left[\frac{1}{10 \cdot n^4}+ n \cdot q^{(m-1) \cdot \frac{\left(k+1\right) \cdot k}{2}+1}_n \cdot a_n, 1- \frac{1}{10 \cdot n^4}+ n \cdot q^{(m-1) \cdot \frac{\left(k+1\right) \cdot k}{2}+1}_n \cdot a_n\right] \times  \prod^{m}_{i=3} \left[\frac{1}{26n^4},1-\frac{1}{26n^4}\right].
\end{split}
\end{equation*}
Thus, such a set $\Phi_n\left(\hat{I}_{n}\right)$ with $\hat{I}_n \in \eta_n$ has a $\theta$-witdth of at most $\frac{1}{n \cdot q^{3m+1}_n}$. \\
Moreover, we see that we can choose $\epsilon=0$ in the definition of a $\left(\gamma, \delta, \epsilon\right)$-distribution: With the notation $A_{\theta} \coloneqq \pi_{\theta}\left(\Phi_n\left(\hat{I}_n\right)\right)$ we have $\Phi_n\left(\hat{I}_n\right)= A_{\theta} \times J$ and so for every $(m-1)$-dimensional interval $\tilde{J}\subseteq J$: 
\begin{equation*}
\frac{\mu\left(\hat{I}_n \cap \Phi^{-1}_n\left(\mathbb{S}^1 \times \tilde{J}\right)\right)}{\mu\left(\hat{I}_n\right)} =  \frac{\mu\left(\Phi_n\left(\hat{I}_n\right) \cap \mathbb{S}^1 \times \tilde{J}\right)}{\mu\left(\Phi_n\left(\hat{I}_n\right)\right)} = \frac{\tilde{\lambda}\left(A_{\theta}\right) \cdot \mu^{(m-1)}\left(\tilde{J}\right)}{\tilde{\lambda}\left(A_{\theta}\right) \cdot \mu^{(m-1)}\left(J\right)} = \frac{\mu^{(m-1)}\left(\tilde{J}\right)}{\mu^{(m-1)}\left(J\right)}
\end{equation*}
because $\Phi_n$ is measure-preserving.
\end{pr}

Furthermore, we show the next property concerning the conjugating map $g_n$ constructed in section \ref{subsection:g}:

\begin{lem} \label{lem:g-phi}
For every $\hat{I}_n \in \eta_n$ we have: $g_n \left(\Phi_n\left(\hat{I}_n\right)\right) = \tilde{g}_{\left[n q^{\sigma}_{n}\right]} \left(\Phi_n\left(\hat{I}_n\right)\right)$.
\end{lem}

\begin{pr}
In the proof of the precedent Lemma \ref{lem:distri} we computed $\Phi_n\left(\hat{I}_{n,k}\right)$ for a partition element $\hat{I}_{n,k}$. Now we have to examine the effect of $g_n = g_{n \cdot q^{1+(m-1) \cdot \frac{\left(k+1\right) \cdot \left(k+2\right)}{2}}_n, \left[n \cdot q^{\sigma}_{n}\right], \frac{1}{8n^4}, \frac{1}{32n^4}}$ on it. \\
Since $260n^4$ divides $q_n$ by Lemma \ref{lem:conv}, there is $u \in \mathbb{Z}$ such that 
\begin{equation*}
\frac{1}{10n^4}=u \cdot \frac{\varepsilon}{b \cdot a} = u \cdot \frac{1}{8 n^4 \cdot \left[nq^{\sigma}_n\right] \cdot n q^{1+(m-1) \cdot \frac{\left(k+1\right) \cdot \left(k+2\right)}{2}}_n}. 
\end{equation*}
By $\frac{1}{26n^4} < \varepsilon = \frac{1}{8n^4}$ and the bound on $a_n$ the boundary of $\Phi_n\left(\hat{I}_{n,k}\right)$ lies in the domain where $g_{n \cdot q^{1+(m-1) \cdot \frac{\left(k+1\right) \cdot \left(k+2\right)}{2}}_n, \left[n \cdot q^{\sigma}_{n}\right], \frac{1}{8n^4}, \frac{1}{32n^4}} = \tilde{g}_{\left[nq^{\sigma}_n\right]}$.
\end{pr}

\section{Criterion for weak mixing} \label{section:crit}

In this section we will prove a criterion for weak mixing on $M=\mathbb{S}^1 \times \left[0,1\right]^{m-1}$ in the setting of the beforehand constructions. For the derivation we need a couple of lemmas. The first one expresses the weak mixing property on the elements of a partial partition $\eta_n$ generally:

\begin{lem} \label{lem:app1}
Let $f \in$ Diff$^{\infty}\left(M, \mu\right)$, $\left(m_n\right)_{n \in \mathbb{N}}$ be a sequence of natural numbers and $\left(\nu_n\right)_{n \in \mathbb{N}}$ be a sequence of partial partitions, where $\nu_n \rightarrow \varepsilon$ and for every $n \in \mathbb{N}$ $\nu_n$ is  the image of a partial partition $\eta_n$ under a measure-preserving diffeomorphism $F_n$, satisfying the following property: For every $m$-dimensional cube $A \subseteq \mathbb{S}^1 \times \left(0,1\right)^{m-1}$ and for every $\epsilon > 0$ there exists $N \in \mathbb{N}$ such that for every $n\geq N$ and for every $\Gamma_n \in \nu_n$ we have
\begin{equation} \label{1}
\left| \mu \left(\Gamma_n \cap f^{-m_n}\left(A\right)\right) - \mu\left(\Gamma_n\right) \cdot \mu\left(A\right) \right| \leq 3 \cdot \epsilon \cdot \mu \left(\Gamma_n\right) \cdot \mu\left(A \right). 
\end{equation}
Then $f$ is weakly mixing.
\end{lem}

\begin{pr}
A diffeomorphism $f$ is weakly mixing if for all measurable sets $A,B \subseteq M$ it holds:
\begin{equation*}
\lim_{n\rightarrow \infty} \left| \mu\left(B\cap f^{-m_n}\left(A\right)\right) - \mu\left(B\right) \cdot \mu\left(A\right) \right| = 0.
\end{equation*}
Since every measurable set in $M = \mathbb{S}^{1} \times \left[0,1\right]^{m-1}$ can be approximated by a countable disjoint union of $m$-dimensional cubes in $\mathbb{S}^{1} \times \left(0,1 \right)^{m-1}$ in arbitrary precision, we only have to prove the statement in case that $A$ is a $m$-dimensional cube in $\mathbb{S}^{1} \times \left(0,1 \right)^{m-1}$. \\
Hence, we consider an arbitrary $m$-dimensional cube $A \subset \mathbb{S}^{1} \times \left(0,1 \right)^{m-1}$. Moreover, let $B \subseteq M$ be a measurable set. Since $\nu_n \rightarrow \varepsilon$ for every $\epsilon \in \left(0,1\right] $ there are $ n \in \mathbb{N}$ and a set $\hat{B} = \bigcup_{i \in \Lambda} \Gamma^{i}_{n}$, where $\Gamma^{i}_{n} \in \nu_n$ and $\Lambda$ is a countable set of indices, such that $\mu\left(B\triangle \hat{B}\right) < \epsilon \cdot \mu\left(B\right) \cdot \mu\left(A\right)$. We obtain for sufficiently large $n$:
\begin{align*}
& \left| \mu\left(B\cap f^{-m_n}\left(A\right)\right) - \mu\left(B\right) \cdot \mu\left(A\right) \right| \\
& \leq \left| \mu\left(B \cap f^{-m_n}\left(A\right)\right) - \mu\left(\hat{B} \cap f^{-m_n}\left(A\right)\right) \right| + \left| \mu\left(\hat{B} \cap f^{-m_n}\left(A\right)\right) - \mu\left(\hat{B}\right) \cdot \mu\left(A\right) \right|\\
 & \quad  + \left|\mu\left(\hat{B}\right) \cdot \mu\left(A\right) - \mu\left(B\right) \cdot \mu\left(A\right) \right|  \\ \displaybreak[0]
& = \left| \mu\left(B \cap f^{-m_n}\left(A\right)\right) - \mu\left(\hat{B} \cap f^{-m_n}\left(A\right)\right) \right| \\
 & \quad + \left| \mu\left(\bigcup_{i \in \Lambda} \left(\Gamma^{i}_{n} \cap f^{-m_n}\left(A\right)\right)\right) - \mu\left(\bigcup_{i \in \Lambda} \Gamma^{i}_{n}\right) \cdot \mu\left(A\right) \right| + \mu\left(A\right) \cdot \left|\mu\left(\hat{B}\right)  - \mu\left(B\right)  \right|  \\ \displaybreak[0]
& \leq  \mu\left(\hat{B} \triangle B\right) + \left| \sum_{i \in \Lambda} \mu\left(\Gamma^{i}_{n} \cap f^{-m_n}\left(A\right)\right) - \mu\left(\Gamma^{i}_{n}\right) \cdot \mu\left(A\right) \right| + \mu\left(A\right) \cdot \mu\left(\hat{B} \triangle B\right) \\ \displaybreak[0]
& \leq  \epsilon \cdot \mu(B) \cdot \mu(A) + \sum_{i \in \Lambda} \left( \left| \mu\left(\Gamma^{i}_{n} \cap f^{-m_n}(A)\right) - \mu\left(\Gamma^{i}_{n}\right) \cdot \mu(A) \right| \right) + \epsilon \cdot \mu(A)^{2} \cdot \mu(B) \\ \displaybreak[0]
& \leq   \sum_{i \in \Lambda} \left( 3 \cdot \epsilon \cdot \mu\left(\Gamma^{i}_{n}\right) \cdot \mu(A) \right) + 2 \cdot \epsilon \cdot \mu(A) \cdot \mu(B) = 3 \cdot \epsilon \cdot \mu(A) \cdot \mu\left(\bigcup_{i \in \Lambda} \hat{I}^{i}_{n}\right) + 2 \cdot \epsilon \cdot \mu(A) \cdot \mu(B) \\ \displaybreak[0]
& =  3 \cdot \epsilon \cdot \mu(A) \cdot \mu\left(\hat{B}\right) + 2 \cdot \epsilon \cdot \mu(A) \cdot \mu(B)\leq  3 \cdot \epsilon \cdot \mu(A) \cdot \left( \mu(B) + \mu\left(\hat{B}\triangle B\right)\right) + 2 \cdot \epsilon \cdot \mu(A) \cdot \mu(B) \\
& \leq  5 \cdot \epsilon \cdot \mu(A) \cdot \mu(B) + 3 \cdot \epsilon^{2} \cdot \mu(A)^{2} \cdot \mu(B).
\end{align*}
This estimate shows $\lim_{n\rightarrow \infty} \left| \mu\left(B\cap f^{-m_n}\left(A\right)\right) - \mu\left(B\right) \cdot \mu\left(A\right) \right| = 0$, because $\epsilon$ can be chosen arbitrarily small.
\end{pr}

In property (\ref{1}) we want to replace $f$ by $f_n$:
\begin{lem} \label{lem:app2}
Let $f = \lim_{n \rightarrow \infty} f_n$ be a diffeomorphism obtained by the constructions in the preceding sections and $\left(m_n\right)_{n \in \mathbb{N}}$ be a sequence of natural numbers fulfilling $d_0\left(f^{m_n},f^{m_n}_{n}\right) < \frac{1}{2^{n}}$. Furthermore, let $\left(\nu_n\right)_{n \in \mathbb{N}}$ be a sequence of partial partitions, where $\nu_n \rightarrow \varepsilon$ and for every $n \in \mathbb{N}$ $\nu_n$ is  the image of a partial partition $\eta_n$ under a measure-preserving diffeomorphism $F_n$, satisfying the following property: For every $m$-dimensional cube $A \subseteq \mathbb{S}^1 \times \left(0,1\right)^{m-1}$ and for every $\epsilon \in \left(0,1\right]$ there exists $N \in \mathbb{N}$ such that for every $n\geq N$ and for every $\Gamma_n \in \nu_n$ we have
\begin{equation} \label{2}
\left| \mu \left(\Gamma_n \cap f^{-m_n}_{n}\left(A\right)\right) - \mu\left(\Gamma_n\right) \cdot \mu\left(A\right) \right| \leq \epsilon \cdot \mu \left(\Gamma_n\right) \cdot \mu\left(A\right).
\end{equation}
Then $f$ is weakly mixing.
\end{lem}

\begin{pr}
We want to show that the requirements of Lemma \ref{lem:app1} are fulfilled. This implies that $f$ is weakly mixing. \\
For it let $A \subseteq \mathbb{S}^{1} \times \left(0,1 \right)^{m-1}$ be an arbitrary $m$-dimensional cube and $\epsilon \in \left(0,1\right]$. \\
We consider two $m$-dimensional cubes $A_1, A_2 \subset \mathbb{S}^{1} \times \left(0,1 \right)^{m-1}$ with $A_1 \subset A \subset A_2$ as well as $\mu\left(A\triangle A_i\right) < \epsilon \cdot \mu\left(A\right)$ and for sufficiently large $n$: dist$\left(\partial A, \partial A_i\right) > \frac{1}{2^{n}}$ for $i=1,2$.

If $n$ is sufficiently large, we obtain for $\Gamma_n \in \nu_n$ and for $i=1,2$ by the assumptions of this Lemma:
\begin{equation*}
\left| \mu \left(\Gamma_n \cap f^{-m_n}_{n}\left(A_i\right)\right) - \mu \left(\Gamma_n\right) \cdot \mu\left(A_i\right) \right| \leq \epsilon \cdot \mu \left(\Gamma_n\right) \cdot \mu\left(A_i\right).
\end{equation*}

Herefrom we conclude $\left(1-\epsilon\right) \cdot \mu \left(\Gamma_n\right) \cdot \mu\left(A_1\right) \leq \mu \left(\Gamma_n \cap f^{-m_n}_{n}\left(A_1\right)\right)$ on the one hand and $\mu \left(\Gamma_n \cap f^{-m_n}_{n}\left(A_2\right)\right) \leq \left(1+\epsilon\right) \cdot \mu \left(\Gamma_n\right) \cdot \mu\left(A_2\right)$ on the other hand. Because of $d_0\left(f^{m_n},f^{m_n}_{n}\right) < \frac{1}{2^{n}}$ the following relations are true:
\begin{align*}
 f^{m_n}_{n}(x) \in A_1 & \Longrightarrow f^{m_n}(x) \in A, \\
 f^{m_n}(x) \in A & \Longrightarrow f^{m_n}_{n}(x) \in A_2.
\end{align*}
 Thus: $\mu \left(\Gamma_n \cap f^{-m_n}_{n}\left(A_1\right)\right) \leq \mu \left(\Gamma_n \cap f^{-m_n}\left(A\right)\right) \leq \mu \left(\Gamma_n \cap f^{-m_n}_{n}\left(A_2\right)\right)$. \\
Altogether, it holds: $\left(1-\epsilon\right) \cdot \mu \left(\Gamma_n\right) \cdot \mu\left(A_1\right) \leq \mu \left(\Gamma_n \cap f^{-m_n}\left(A\right)\right) \leq \left(1+\epsilon\right) \cdot \mu \left(\Gamma_n \right) \cdot \mu\left(A_2\right)$. Therewith, we obtain the following estimate from above:
\begin{align*}
& \mu \left(\Gamma_n \cap f^{-m_n}\left(A\right)\right) - \mu \left(\Gamma_n\right) \cdot \mu\left(A\right) \\
& \leq \left(1+\epsilon \right) \cdot \mu \left(\Gamma_n \right) \cdot \mu \left(A_2\right) - \mu \left(\Gamma_n\right) \cdot \mu\left(A_2\right) + \mu \left(\Gamma_n \right) \cdot \left( \mu\left(A_2 \right) - \mu\left(A\right) \right) \\
& \leq \epsilon \cdot \mu \left(\Gamma_n\right) \cdot \mu\left(A_2\right) + \mu \left(\Gamma_n\right) \cdot \mu\left(A_2 \triangle A\right) \leq \epsilon \cdot \mu \left(\Gamma_n\right) \cdot \left(\mu(A) + \mu\left(A_2 \triangle A\right)\right) + \epsilon \cdot \mu \left(\Gamma_n\right) \cdot \mu\left(A\right) \\
& \leq 2 \cdot \epsilon \cdot \mu \left(\Gamma_n\right) \cdot \mu\left(A\right) + \epsilon^{2} \cdot \mu \left(\Gamma_n\right) \cdot \mu\left(A\right) \leq 3 \cdot \epsilon \cdot \mu \left(\Gamma_n\right) \cdot \mu\left(A\right)
\end{align*}

Furthermore, we deduce the following estimate from below in an analogous way:
\begin{equation*}
\mu \left(\Gamma_n \cap f^{-m_n}\left(A\right)\right) - \mu \left(\Gamma_n\right) \cdot \mu\left(A\right) \geq - 3 \cdot \epsilon \cdot \mu \left(\Gamma_n\right) \cdot \mu\left(A\right)
\end{equation*}

Hence, we get: $\left| \mu \left(\Gamma_n \cap f^{-m_n}\left(A\right)\right) - \mu \left(\Gamma_n\right) \cdot \mu\left(A\right) \right| \leq 3 \cdot \epsilon \cdot \mu \left(\Gamma_n\right) \cdot \mu\left(A\right)$, i.e. the requirements of Lemma \ref{lem:app1} are met.
\end{pr}

Now we concentrate on the setting of our explicit constructions:

\begin{lem} \label{lem:points}
Consider  the sequence of partial partitions $\left(\eta_n\right)_{n \in \mathbb{N}}$ constructed in section \ref{subsubsection:eta} and the diffeomorphisms $g_n$ from chapter \ref{subsection:g}. Furthermore, let $\left(H_n\right)_{n \in \mathbb{N}}$ be a sequence of measure-preserving smooth diffeomorphisms satisfying $\left\|DH_{n-1}\right\| \leq \frac{\ln\left(q_n\right)}{n}$ for every $n \in \mathbb{N}$ and define the partial partitions $\nu_n = \left\{ \Gamma_n = H_{n-1} \circ g_n \left( \hat{I}_n \right) \ : \ \hat{I}_n \in \eta_n \right\}$. \\
Then we get $\nu_n \rightarrow \varepsilon$.
\end{lem}

\begin{pr}
By construction $\eta_n = \left\{ \hat{I}^{i}_{n} : i \in \Lambda_n \right\}$, where $\Lambda_n$ is a countable set of indices. Because of $\eta_n \rightarrow \varepsilon$ it holds $\lim_{n\rightarrow\infty} \mu \left(\bigcup_{i \in \Lambda_n} \hat{I}^{i}_{n} \right) = 1$. Since $H_{n-1} \circ g_n$ is measure-preserving, we conclude:  
\begin{equation*}
\lim_{n\rightarrow\infty} \mu \left(\bigcup_{i \in \Lambda_n} \Gamma^{i}_{n}\right) = \lim_{n\rightarrow\infty} \mu \left(\bigcup_{i \in \Lambda_n} H_{n-1} \circ g_n\left(\hat{I}^{i}_{n} \right)\right) = \lim_{n\rightarrow\infty} \mu \left(H_{n-1} \circ g_n \left(\bigcup_{i \in \Lambda_n} \hat{I}^{i}_{n} \right)\right) = 1.
\end{equation*}
For any $m$-dimensional cube with side length $l_n$ it holds: diam$\left(W_n\right) = \sqrt{m} \cdot l_n$. Because every element of the partition $\eta_n$ is contained in a cube of side length $\frac{1}{q_n}$ it follows for every $i \in \Lambda_n$: diam$\left(\hat{I}^{i}_{n}\right) \leq \sqrt{m} \cdot \frac{1}{q_n}$. Furthermore, we saw in Lemma \ref{lem:outer}: $g_n\left(\hat{I}^i_n\right) = \tilde{g}_{\left[n q^{\sigma}_{n}\right]} \left(\hat{I}^i_n\right)$ for every $i \in \Lambda_n$. Hence, for every $\Gamma^{i}_{n} = H_{n-1} \circ \tilde{g}_{\left[n q^{\sigma}_{n}\right]} \left(I^{i}_{n}\right)$ :
\begin{equation*}
\text{diam}\left(\Gamma^{i}_{n}\right) \leq \left\|DH_{n-1} \right\|_0 \cdot \left\| D\tilde{g}_{\left[n q^{\sigma}_{n}\right]} \right\|_0 \cdot \text{diam}\left(\hat{I}^{i}_{n}\right) \leq \frac{\ln\left(q_n\right)}{n} \cdot \left[n \cdot q^{\sigma}_{n} \right] \cdot \frac{\sqrt{m} }{q_n} \leq \sqrt{m} \cdot q^{\sigma-1}_{n} \cdot \ln\left(q_n\right).
\end{equation*}
Because of $\sigma < 1$ we conclude $\lim_{n\rightarrow\infty} $diam$\left(\Gamma^{i}_{n}\right) = 0$ and consequently $\nu_n \rightarrow \varepsilon$.
\end{pr}

In the following the Lebesgue measures on $\mathbb{S}^1$, $\left[0,1\right]^{m-2}$, $\left[0,1\right]^{m-1}$ are denoted by $\tilde{\lambda}$, $\mu^{(m-2)}$ and $\tilde{\mu}$ respectively. The next technical result is needed in the proof of Lemma \ref{lem:cube}.

\begin{lem} \label{lem:help}
Given an interval on the $r_1$-axis of the form $K = \bigcup_{k \in \mathbb{Z}, k_1 \leq k \leq k_2} \left[\frac{k \cdot \varepsilon}{b \cdot a}, \frac{(k+1) \cdot \varepsilon}{b \cdot a} \right]$, where $k_1, k_2 \in \mathbb{Z}$ with $\frac{b \cdot a}{\varepsilon} \cdot \delta \leq k_1 < k_2 \leq \frac{b \cdot a}{\varepsilon} - \frac{b \cdot a}{\varepsilon} \cdot \delta - 1$, and a $(m-2)$-dimensional interval $Z$ in $\left(r_2,...,r_{m-1}\right)$, let $K_{c, \gamma}$ denote the cuboid $\left[c, c + \gamma \right] \times K \times Z$ for some $\gamma > 0$. We consider the diffeomorphism $g_{a,b, \varepsilon,\delta}$ constructed in subsection \ref{subsection:g}  and an interval $L=\left[l_1, l_2 \right]$ of $\mathbb{S}^1$ satisfying $\tilde{\lambda}\left(L\right)\geq 4 \cdot \frac{1-2 \varepsilon}{a} - \gamma$. \\
If $b \cdot \lambda (K) > 2$, then for the set $Q:= \pi_{\vec{r}}\left(K_{c, \gamma} \cap g^{-1}_{a,b, \varepsilon,\delta}\left(L \times K \times Z\right)\right)$ we have:
\begin{align*}
& \left| \tilde{\mu}\left(Q\right) - \lambda \left(K\right) \cdot \tilde{\lambda} \left(L\right) \cdot \mu^{\left(m-2\right)}\left(Z\right) \right| \\
& \leq \left(\frac{2}{b}\cdot \tilde{\lambda} \left(L \right) + \frac{2 \cdot \gamma}{b} + \gamma \cdot \lambda \left(K\right) + 4 \cdot \frac{1-2\varepsilon}{a} \cdot \lambda (K) + 8 \cdot \frac{1-2\varepsilon}{b \cdot a}\right) \cdot \mu^{\left(m-2\right)}\left(Z\right).
\end{align*}
\end{lem}

\begin{pr}
We consider the diffeomorphism $\tilde{g}_b: M \rightarrow M$, $\left(\theta, r_1,...,r_{m-1}\right) \mapsto \left( \theta + b \cdot r_1, r_1,...,r_{m-1}\right)$ and the set:
\begin{align*}
Q_b & \coloneqq \pi_{\vec{r}}\left( K_{c,\gamma} \cap \tilde{g}^{-1}_b\left(L \times K \times Z\right)\right) \\
& = \left\{ \left(r_1, r_2, ...,r_{m-1}\right) \in K \times Z \ : \  \left( \theta + b \cdot r_1, \vec{r}\right) \in L \times K \times Z , \theta \in \left[c, c+ \gamma \right] \right\} \\ & = \left\{ \left(r_1, r_2, ..., r_{m-1}\right) \in K \times Z \ : \  b \cdot r_1 \in \left[l_1 - c - \gamma , l_2 - c \right] \ \text{ mod } 1 \right\}.
\end{align*}
The interval $b \cdot K$ seen as an interval in $\mathbb{R}$ does not intersect more than $b \cdot \lambda(K) + 2$ and not less than $b \cdot \lambda \left(K\right) - 2$ intervals of the form $\left[ i, i+1 \right]$ with $i \in \mathbb{Z}$. By construction of the map $g_{a,b, \varepsilon,\delta}$ it holds for $\Delta_l \coloneqq \left[\frac{l \cdot \varepsilon}{b \cdot a}, \frac{(l+1) \cdot \varepsilon}{b \cdot a}\right]$ in consideration: $\pi_{\vec{r}} \left(g_{a,b, \varepsilon,\delta}\left(\left[c,c+ \gamma \right] \times \Delta_l \times Z\right)\right) = \Delta_l \times Z$. \\
\textbf{Claim: } A resulting interval on the $r_1$-axis of $K_{c,\gamma} \cap \tilde{g}^{-1}_b\left(L \times K \times Z\right)$ and the corresponding $r_1$-projection of $K_{c,\gamma} \cap g^{-1}_{a,b,\varepsilon}\left(L \times K \times Z\right)$ can differ by a length of at most $4 \cdot \frac{1-2 \varepsilon}{b \cdot a}$. \\
\textbf{Proof: } If $\left\{c\right\} \times \Delta_l \times Z$ (resp. $\left\{c + \gamma \right\} \times \Delta_l \times Z$) are contained in the domain, where $g_{a,b,\varepsilon} =\tilde{g}_b$, the left (resp. the right) boundaries of $\pi_{\theta}\left(g_{a,b, \varepsilon,\delta}\left(\left[c,c+ \gamma \right] \times \Delta_l \times Z\right)\right)$ and $\pi_{\theta}\left(\tilde{g}_b\left(\left[c,c+ \gamma \right] \times \Delta_l \times Z\right)\right)$ coincide. Otherwise, i.e. $c \in \left(\frac{k}{a} + \varepsilon, \frac{k+1}{a} - \varepsilon\right)$ (resp. $c + \gamma \in \left(\frac{k}{a} + \varepsilon, \frac{k+1}{a} - \varepsilon\right)$) the sets $\pi_{\theta}\left(g_{a,b, \varepsilon,\delta}\left(\left\{c\right\} \times \Delta_l \times Z\right)\right)$ and $\pi_{\theta}\left(\tilde{g}_b\left(\left\{c\right\} \times \Delta_l \times Z\right)\right)$ (resp. $\pi_{\theta}\left(g_{a,b, \varepsilon,\delta}\left(\left\{c + \gamma \right\} \times \Delta_l \times Z\right)\right)$ and $\pi_{\theta}\left(\tilde{g}_b\left(\left\{c + \gamma \right\} \times \Delta_l \times Z\right)\right)$) differ by a length of at most $\frac{1-2\varepsilon}{a}$. Since $\pi_{\theta}\left(\tilde{g}_{b}\left(\left\{u\right\} \times \Delta_l \times Z\right)\right)$ for arbitrary $u \in \mathbb{S}^1$ has a length of $\frac{\varepsilon}{a}$ on the $\theta$-axis, this discrepancy will be equalised after at most $\frac{1-2\varepsilon}{a} : \frac{\varepsilon}{a} = \frac{1-2 \varepsilon}{\varepsilon}$ blocks $\Delta_l$ on the $r_1$-axis. Thus, the resulting interval on the $r_1$-axis of $K_{c,\gamma} \cap \tilde{g}^{-1}_b\left(L \times K \times Z\right)$ and the corresponding $r_1$-projection of $K_{c,\gamma} \cap g^{-1}_{a,b,\varepsilon}\left(L \times K \times Z\right)$ can differ by a length of at most $4 \cdot \frac{1-2\varepsilon}{\varepsilon} \cdot \frac{\varepsilon}{b \cdot a} = 4 \cdot \left(1-2\varepsilon\right)\frac{1}{b \cdot a}$. \qed

Therefore, we compute on the one side:
\begin{align*}
& \tilde{\mu}\left(Q\right) \leq  \left(b \cdot \lambda\left(K\right) + 2\right) \cdot \left(\frac{l_2 - \left(l_1 - \gamma\right)}{b} + 4 \cdot \frac{1-2 \varepsilon}{b \cdot a}\right) \cdot \mu^{(m-2)}\left(Z\right) \\
& = \left(\lambda\left(K\right) \cdot \tilde{\lambda} \left(L\right) + 2 \cdot \frac{\tilde{\lambda}\left(L\right)}{b} + \lambda\left(K\right) \cdot \gamma + \frac{2 \cdot \gamma}{b} + 4 \cdot \lambda(K) \cdot \frac{1-2 \varepsilon}{a} + 8 \cdot \frac{1-2 \varepsilon}{b \cdot a} \right) \cdot \mu^{(m-2)}\left(Z\right)
\end{align*}
and on the other side
\begin{align*}
& \tilde{\mu}\left(Q\right) \geq  \left(b \cdot \lambda\left(K\right) - 2\right) \cdot \left(\frac{l_2 - \left(l_1 - \gamma\right)}{b} - 4 \cdot \frac{1-2 \varepsilon}{b \cdot a} \right) \cdot \mu^{(m-2)}\left(Z\right) \\
& = \left(\lambda\left(K\right) \cdot \tilde{\lambda}\left(L\right) - 2 \cdot \frac{\tilde{\lambda}\left(L\right)}{b} + \lambda\left(K\right) \cdot \gamma - \frac{2 \cdot \gamma}{b} - 4 \cdot \lambda(K) \cdot \frac{1-2 \varepsilon}{a} + 8 \cdot \frac{1-2 \varepsilon}{b \cdot a} \right) \cdot \mu^{(m-2)}\left(Z\right).
\end{align*}
Both equations together yield:
\begin{align*}
& \left| \tilde{\mu}\left(Q\right) - \lambda\left(K\right) \cdot \tilde{\lambda} \left(L\right) \cdot \mu^{(m-2)}\left(Z\right) - \gamma \cdot \lambda\left(K\right) \cdot \mu^{(m-2)}\left(Z\right) - 8 \cdot \frac{1-2 \varepsilon}{b \cdot a}\cdot \mu^{(m-2)}\left(Z\right) \right| \\
& \leq \left(\frac{2}{b}\cdot \tilde{\lambda}\left(L\right) + \frac{2 \cdot \gamma}{b} + 4 \cdot \lambda(K) \cdot \frac{1-2 \varepsilon}{ a} \right) \cdot \mu^{(m-2)}\left(Z\right).
\end{align*}
The claim follows because
\begin{align*}
& \left| \tilde{\mu}\left(Q\right) - \lambda\left(K\right) \cdot \tilde{\lambda} \left(L\right) \cdot \mu^{(m-2)}\left(Z\right) \right| - \gamma \cdot \lambda\left(K\right) \cdot \mu^{(m-2)}\left(Z\right) - 8 \cdot \frac{1-2 \varepsilon}{b \cdot a} \cdot \mu^{(m-2)}\left(Z\right) \\
& \leq \left| \tilde{\mu}\left(Q\right) - \lambda\left(K\right) \cdot \tilde{\lambda} \left(L\right) \cdot \mu^{(m-2)}\left(Z\right)  - \gamma \cdot \lambda\left(K\right) \cdot \mu^{(m-2)}\left(Z\right) - 8 \cdot \frac{1-2 \varepsilon}{b \cdot a} \cdot \mu^{(m-2)}\left(Z\right) \right|.
\end{align*}
\end{pr}

\begin{lem} \label{lem:cube}
Let $n\geq5$, $g_n$ as in section \ref{subsection:g} and $\hat{I}_n \in \eta_n$, where $\eta_n$ is the partial partition constructed in section \ref{subsubsection:eta}. For the diffeomorphism $\phi_n$ constructed in section \ref{subsection:phi} and $m_n$ as in chapter \ref{section:distri} we consider $\Phi_n = \phi_n \circ R^{m_n}_{\alpha_{n+1}} \circ \phi^{-1}_{n}$ and denote $\pi_{\vec{r}} \left( \Phi_n \left(\hat{I}_n\right)\right)$ by J. \\
Then for every $m$-dimensional cube $S$ of side length $q^{-\sigma}_{n}$ lying in $\mathbb{S}^1 \times J$ we get
\begin{equation}
\left| \mu\left(\hat{I} \cap \Phi^{-1}_n \circ g^{-1}_{n}\left(S\right)\right) \cdot \tilde{\mu}\left(J\right) - \mu\left(\hat{I}\right) \cdot \mu\left(S\right) \right| \leq \frac{21}{n} \cdot \mu\left(\hat{I}\right) \cdot \mu\left(S\right).
\end{equation}
\end{lem}

In other words this Lemma tells us that a partition element is ``almost uniformly distributed'' under $g_n \circ \Phi_n$ on the whole manifold $M= \mathbb{S}^1 \times \left[0,1\right]^{m-1}$.

\begin{pr}
Let $S$ be a $m$-dimensional cube with side length $q^{-\sigma}_{n}$ lying in $\mathbb{S}^{1} \times J$. Furthermore, we denote:
\begin{equation*}
S_\theta = \pi_{\theta}\left(S\right) \ \ \ \ \ \ \ \ \ \  S_{r_1} = \pi_{r_1}\left(S\right) \ \ \ \ \ \ \ \ \ \ S_{\tilde{\vec{r}}} = \pi_{\left(r_2,...,r_{m-1}\right)} \left(S\right) \ \ \ \ \ \ \ \ \ \ S_r = S_{r_1} \times S_{\tilde{\vec{r}}} = \pi_{\vec{r}}\left(S\right)
\end{equation*}
Obviously: $\tilde{\lambda}\left(S_{\theta}\right) = \lambda\left(S_{r_1}\right) = q^{-\sigma}_{n}$ and $\tilde{\lambda}\left(S_{\theta}\right) \cdot \lambda\left(S_{r_1}\right) \cdot \mu^{(m-2)}\left(S_{\tilde{\vec{r}}}\right)= \mu\left(S\right) = q^{-m \sigma}_{n}$. \\
According to Lemma \ref{lem:distri} $\Phi_n$ $\left(\frac{1}{n \cdot q^{m}_{n}}, \frac{1}{n^4},\frac{1}{n}\right)$-distributes the partition element $\hat{I}_n \in \eta_n$, in particular $\Phi_n\left(\hat{I}_n\right) \subseteq \left[c,c+\gamma\right] \times J$ for some $c \in \mathbb{S}^1$ and some $\gamma \leq \frac{1}{n \cdot q^m_n}$. Furthermore, we saw in the proof of Lemma \ref{lem:g-phi} that $\left[c,c+\gamma\right] \times J$ is contained in the interior of the step-by-step domains of the map $g_n$ and that on its boundary $g_n = \tilde{g}_{\left[n q^{\sigma}_{n}\right]}$ holds. Particularly it follows $\gamma \geq \frac{1-2\varepsilon}{a}$ in case of $g_n = g_{a,b, \varepsilon,\delta}$. For $l \in \mathbb{Z}$, $0\leq l \leq \frac{b\cdot a}{\varepsilon}-1$ we introduce the set $\Delta_l = \left[\frac{l \varepsilon}{ba}, \frac{(l+1) \varepsilon}{ba}\right]$ and therewith we consider
\begin{equation*}
\tilde{S}_{r_1} \coloneqq \bigcup_{\Delta_l \subseteq S_{r_1}} \Delta_l; \ \ \ \ \tilde{S}_r \coloneqq \bigcup_{\Delta_l \subseteq S_{r_1}} \Delta_l \times S_{\tilde{\vec{r}}} \ \ \ \ \text{ as well as } \ \ \ \  \tilde{S} \coloneqq S_{\theta} \times \tilde{S}_r \subseteq S
\end{equation*} 
Using the triangle inequality we obtain
\begin{align*}
& \left| \mu\left(\hat{I} \cap \Phi^{-1}_n\left(g^{-1}_{n}(S)\right)\right) \cdot \tilde{\mu}\left(J\right) - \mu\left(\hat{I}\right) \cdot \mu\left(S\right) \right| \\ 
\leq  & \left| \mu\left(\hat{I} \cap \Phi^{-1}_n\left(g^{-1}_{n}(S)\right)\right) - \mu\left(\hat{I} \cap \Phi^{-1}_n\left(g^{-1}_{n}\left(\tilde{S}\right)\right)\right) \right| \cdot \tilde{\mu}\left(J\right) \\
& + \left| \mu\left(\hat{I} \cap \Phi^{-1}_n\left(g^{-1}_{n}\left(\tilde{S}\right)\right)\right) \cdot \tilde{\mu}\left(J\right) - \mu\left(\hat{I}\right) \mu\left(\tilde{S}\right)\right| + \mu\left(\hat{I}\right) \cdot \left|\mu\left(\tilde{S}\right) - \mu\left(S\right) \right|
\end{align*}
Here $\left|\mu\left(\tilde{S}\right) - \mu\left(S\right) \right| = \mu\left(S \setminus \tilde{S}\right) \leq 2 \cdot \frac{\varepsilon}{b \cdot a} \cdot \tilde{\lambda}\left(S_{\theta}\right) \cdot \mu^{(m-2)}\left(S_{\tilde{\vec{r}}}\right) \leq 2 \cdot \frac{\varepsilon}{a} \cdot \mu\left(S\right)$, where we used $b=\left[n \cdot q^{\sigma}_{n}\right] \geq q^{\sigma}_{n}$ in case of $n>4$. Since $\Phi_n$ and $g_n$ are measure-preserving, we additionally obtain: $\left| \mu\left(\hat{I} \cap \Phi^{-1}_n\left(g^{-1}_{n}(S)\right)\right) - \mu\left(\hat{I} \cap \Phi^{-1}_n\left(g^{-1}_{n}\left(\tilde{S}\right)\right)\right) \right| \leq \mu\left(S \setminus \tilde{S}\right) \leq 2 \cdot \frac{\varepsilon}{a} \cdot \mu\left(S\right)$. \\
In the proof of Lemma \ref{lem:g-phi} we observed $\mu\left(\Phi_n\left(\hat{I}\right)\right)=\frac{1}{a} \cdot \left(1-\frac{2}{26n^4}\right) \cdot \tilde{\mu}\left(J\right)$. Hence:
\begin{align*}
& \left| \mu\left(\hat{I} \cap \Phi^{-1}_n\left(g^{-1}_{n}(S)\right)\right) - \mu\left(\hat{I} \cap \Phi^{-1}_n\left(g^{-1}_{n}\left(\tilde{S}\right)\right)\right) \right| \cdot \tilde{\mu}\left(J\right) \leq 2 \cdot \frac{\varepsilon}{a} \cdot \mu\left(S\right) \cdot \tilde{\mu}\left(J\right) \\
&= 2 \cdot \frac{\varepsilon}{1-\frac{2}{26n^4}} \cdot \mu\left(S\right) \cdot \mu\left(\Phi_n\left(\hat{I}\right)\right) \leq 4 \cdot \varepsilon \cdot \mu\left(S\right) \cdot \mu\left(\Phi_n\left(\hat{I}\right)\right) = 4 \cdot \varepsilon \cdot \mu\left(S\right) \cdot \mu\left(\hat{I}\right)
\end{align*}
Thus, we obtain:
\begin{equation} \label{triangle}
\begin{split}
& \left| \mu\left(\hat{I} \cap \Phi^{-1}_n\left(g^{-1}_{n}(S)\right)\right) \cdot \tilde{\mu}\left(J\right) - \mu\left(\hat{I}\right) \cdot \mu\left(S\right) \right| \\
& \leq \left| \mu\left(\hat{I} \cap \Phi^{-1}_n\left(g^{-1}_{n}\left(\tilde{S}\right)\right)\right) \cdot \tilde{\mu}\left(J\right) - \mu\left(\hat{I}\right) \mu\left(\tilde{S}\right)\right| + 5 \cdot \varepsilon \cdot \mu\left(S\right) \cdot \mu\left(\hat{I}\right)
\end{split}
\end{equation}
Next, we want to estimate the first summand. By construction of the map $g_n = g_{a,b, \varepsilon,\delta}$ and the definition of $\tilde{S}$ it holds: $\Phi_n\left(\hat{I}\right) \cap g^{-1}_{n} \left(\tilde{S}\right) \subseteq \left[c,c+\gamma\right] \times \tilde{S}_r \eqqcolon K_{c, \gamma}$. Considering the proof of Lemma \ref{lem:g-phi} again, we obtain $g_n\left(K_{c, \gamma}\right) = \tilde{g}_{\left[n q^{\sigma}_{n}\right]}\left(K_{c, \gamma}\right)$ (since $c$ and $c+\gamma$ are in the domain where $g_n = \tilde{g}_{\left[n q^{\sigma}_{n}\right]}$ holds). \\
Because of Lemma \ref{lem:distri} $2\gamma \leq \frac{2}{n \cdot q^{m}_{n}} < q^{-\sigma}_{n}$ for $n>2$. So we can define a cuboid $S_1 \subseteq \tilde{S}$, where $S_1 \coloneqq \left[ s_1 + \gamma, s_2 - \gamma \right] \times \tilde{S}_r$ using the notation $S_{\theta} = \left[s_1, s_2\right]$. We examine the two sets
\begin{equation*}
Q := \pi_{\vec{r}}\left( K_{c, \gamma} \cap g^{-1}_{n}\left(S_{\theta} \times \tilde{S}_{r} \right)\right) \ \ \ \ \ \ \ \ \  Q_1 := \pi_{\vec{r}}\left( K_{c, \gamma} \cap g^{-1}_{n}\left(\left[s_1 + \gamma, s_2 - \gamma \right] \times \tilde{S}_{r} \right)\right)
\end{equation*}
As seen above $\Phi_n\left(\hat{I}\right) \cap g^{-1}_{n} \left(\tilde{S}\right) \subseteq K_{c, \gamma}$. Hence $\Phi_n\left(\hat{I}\right) \cap g^{-1}_{n} \left(\tilde{S}\right) \subseteq \Phi_n\left(\hat{I}\right) \cap g^{-1}_{n} \left(\tilde{S}\right) \cap K_{c, \gamma}$, which implies $\Phi_n \left(\hat{I}\right) \cap g^{-1}_{n}\left(\tilde{S}\right) \subseteq \Phi_n \left(\hat{I}\right) \cap \left(\mathbb{S}^{1} \times Q\right)$. \\
\textbf{Claim: } On the other hand: $\Phi_n \left(\hat{I}\right) \cap \left(\mathbb{S}^{1} \times Q_1\right) \subseteq \Phi_n \left(\hat{I}\right) \cap g^{-1}_{n}\left(\tilde{S}\right)$. \\
\textbf{Proof of the claim: } For $\left(\theta, \vec{r}\right) \in \Phi_n \left(\hat{I}\right) \cap \left(\mathbb{S}^{1} \times Q_1\right)$ arbitrary it holds $\left(\theta, \vec{r}\right) \in \Phi_n \left(\hat{I}\right)$, i.e. $\theta \in \left[c, c + \gamma\right]$, and $\vec{r} \in \pi_{\vec{r}}\left( K_{c, \gamma} \cap g^{-1}_{n}\left(\left[s_1 + \gamma, s_2 - \gamma \right] \times \tilde{S}_{r} \right)\right)$, i.e. in particular $\vec{r} \in \tilde{S}_r$. This implies the existence of $\bar{\theta} \in \left[c, c + \gamma\right]$ satisfying $\left(\bar{\theta}, \vec{r}\right) \in K_{c,\gamma} \cap g^{-1}_{n} \left(S_1\right)$. Hence, there are $\beta \in \left[s_1 + \gamma, s_2 - \gamma \right]$ and $\vec{r}_1 \in \tilde{S}_r$, such that $g_n \left( \bar{\theta}, \vec{r} \right) = \left(\beta,\vec{r}_1\right)$. Because of $\bar{\theta} \in \left[c,c+\gamma\right]$ and $\vec{r} \in \tilde{S}_r$ the point $\left( \bar{\theta}, \vec{r} \right)$ is contained in one cuboid of the form $\Delta_{a,b, \varepsilon}$. Since $\theta \in \left[c,c+\gamma\right]$, $\left( \theta, \vec{r} \right)$ is contained in the same $\Delta_{a,b, \varepsilon}$. Thus, $\pi_{\vec{r}} \left(g_n\left(\theta, \vec{r}\right)\right) \in \tilde{S}_r$. Furthermore, $g_n\left(\theta, \vec{r}\right)$ and $g_n\left(\bar{\theta}, \vec{r}\right)$ are in a distance of at most $\gamma$ on the $\theta$-axis, because $\theta, \bar{\theta} \in \left[c,c+\gamma\right]$, i.e. $\left|\theta - \bar{\theta} \right| \leq \gamma$, $g_n \left(K_{c, \gamma}\right)= \tilde{g}_{\left[n q^{\sigma}_{n}\right]}\left(K_{c, \gamma}\right)$ and the map $\tilde{g}_{\left[n q^{\sigma}_{n}\right]}$ preserves the distances on the $\theta$-axis. Thus, there are $\bar{\beta} \in \left[s_1,s_2\right]$ and $\vec{r}_2 \in \tilde{S}_r$ such that $g_n\left(\theta, \vec{r}\right)= \left(\bar{\beta}, \vec{r}_2\right)$. So $\left(\theta, \vec{r}\right) \in \Phi_n\left(\hat{I}\right)\cap g^{-1}_{n}\left(\tilde{S}\right)$. \qed \\
Altogether, the following inclusions are true:
\begin{equation*}
\Phi_n \left(\hat{I}\right) \cap \left(\mathbb{S}^{1} \times Q_1\right) \subseteq \Phi_n \left(\hat{I}\right) \cap g^{-1}_{n}\left(\tilde{S}\right) \subseteq \Phi_n \left(\hat{I}\right) \cap \left(\mathbb{S}^{1} \times Q\right).
\end{equation*}
Thus, we obtain:
\begin{equation} \label{big}
\begin{split}
& \left| \mu\left(\hat{I} \cap \Phi^{-1}_n\left(g^{-1}_{n}(\tilde{S})\right)\right) \cdot \tilde{\mu}\left(J\right) - \mu\left(\hat{I}\right) \cdot \mu\left(\tilde{S}\right) \right| \\
\leq \max  \Bigg{(} & \left|\mu\left(\hat{I} \cap \Phi^{-1}_n\left(\mathbb{S}^{1} \times Q\right)\right) \cdot \tilde{\mu}\left(J\right) - \mu\left(\hat{I}\right) \cdot \mu\left(\tilde{S}\right) \right| , \\
 & \left|\mu\left(\hat{I} \cap \Phi^{-1}_n\left(\mathbb{S}^{1} \times Q_1\right)\right) \cdot \tilde{\mu}\left(J\right) - \mu\left(\hat{I}\right) \cdot \mu\left(\tilde{S}\right) \right| \Bigg{)}
\end{split}
\end{equation}
We want to apply Lemma \ref{lem:help} for $K = \tilde{S}_{r_1}$, $L=S_{\theta}$, $Z = S_{\tilde{\vec{r}}}$ and $b=\left[n\cdot q^{\sigma}_{n}\right]$ (note that $4 \cdot \frac{1-2\varepsilon}{a}-\gamma \leq 3 \cdot \frac{1-2\varepsilon}{a} \leq \frac{3}{n \cdot q^{m}_{n}} < \frac{1}{q^{\sigma}_{n}} = \tilde{\lambda}\left(L\right)$ because of the mentioned relation $\gamma\geq\frac{1-2\varepsilon}{a}$ and for $n > 4$: $b \cdot \lambda(K) = \left[n q^{\sigma}_{n}\right] \cdot q^{-\sigma}_{n} \geq \frac{1}{2} n q^{\sigma}_{n} \cdot q^{-\sigma}_{n} >2 $):
\begin{align*}
& \left| \tilde{\mu}\left(Q\right) - \mu\left(\tilde{S}\right) \right| \\
& \leq \left(\frac{2}{\left[n \cdot q^{\sigma}_{n}\right] } \cdot \tilde{\lambda}\left(S_{\theta}\right) + \frac{2 \gamma}{\left[n\cdot q^{\sigma}_{n}\right]} + \gamma \cdot \lambda\left(\tilde{S}_{r_1}\right) + 4 \cdot \frac{1-2\varepsilon}{a} \lambda\left(\tilde{S}_{r_1}\right) + 8 \cdot \frac{1-2 \varepsilon}{\left[n q^{\sigma}_{n}\right] \cdot a} \right) \cdot \mu^{(m-2)}\left(S_{\tilde{\vec{r}}}\right) \\
& \leq \left(\frac{4}{n \cdot q^{\sigma}_{n} } \tilde{\lambda}\left(S_{\theta}\right) + \frac{4}{n \cdot q^{\sigma}_{n} \cdot q^{\sigma}_{n}} + \frac{1}{n \cdot q^{\sigma}_{n}} \lambda\left(S_{r_1}\right) + 4 \cdot \frac{1-2\varepsilon}{n \cdot q^{m}_{n}} \lambda\left(S_{r_1}\right) + \frac{16 \cdot \left(1-2\varepsilon\right)}{n \cdot q^{\sigma}_{n} \cdot n \cdot q^{m}_{n}} \right) \cdot \mu^{(m-2)}\left(S_{\tilde{\vec{r}}}\right) \\
& \leq \frac{14}{n} \cdot \mu\left(S\right). 
\end{align*}
In particular, we receive from this estimate: $\frac{14}{n} \cdot \mu\left(S\right) \geq \tilde{\mu}\left(Q\right) - \mu\left( \tilde{S} \right) \geq \tilde{\mu}\left(Q\right) - \mu\left(S\right)$, hence: $\tilde{\mu}\left(Q\right) \leq \left(1 + \frac{14}{n} \right) \cdot \mu\left(S\right) \leq 4 \cdot \mu\left(S\right)$. \\
Analogously, we obtain: $\tilde{\mu}\left(Q_1\right) \leq 4 \cdot \mu \left(S\right)$ as well as $\left|\tilde{\mu}\left(Q_1\right) - \mu\left(\tilde{S}_1\right)\right| \leq \frac{14}{n} \cdot \mu\left(S\right)$. \\
Since $Q$ as well as $Q_1$ are a finite union of disjoint $\left(m-1\right)$-dimensional intervals contained in $J$ and $\Phi_n$ $\left(\frac{1}{n \cdot q^{m}_{n}}, \frac{1}{n^4}, \frac{1}{n}\right)$-distributes the interval $\hat{I}$, we get:
\begin{equation*}
\left| \mu\left(\hat{I} \cap \Phi^{-1}_n\left(\mathbb{S}^{1} \times Q\right)\right) \cdot \tilde{\mu}\left(J\right) - \mu\left(\hat{I}\right) \cdot \tilde{\mu}\left(Q\right) \right| \leq \frac{1}{n} \cdot \mu\left(\hat{I}\right) \cdot \tilde{\mu}\left(Q\right) \leq \frac{4}{n} \cdot \mu\left(\hat{I}\right) \cdot \mu\left(S\right)
\end{equation*}
as well as
\begin{equation*}
\left| \mu\left(\hat{I} \cap \Phi^{-1}_n\left(\mathbb{S}^{1} \times Q_1\right)\right) \cdot \tilde{\mu}\left(J\right) - \mu\left(\hat{I}\right) \cdot \tilde{\mu}\left(Q_1\right) \right| \leq \frac{1}{n} \cdot \mu\left(\hat{I}\right) \cdot \tilde{\mu}\left(Q_1\right) \leq \frac{4}{n} \cdot \mu\left(\hat{I}\right) \cdot \mu\left(S\right).
\end{equation*}
Now we can proceed
\begin{align*}
& \left| \mu\left(\hat{I} \cap \Phi^{-1}_n\left(\mathbb{S}^{1} \times Q\right)\right) \cdot \tilde{\mu}\left(J\right) - \mu\left(\hat{I}\right) \cdot \mu\left(\tilde{S}\right) \right| \\
& \leq \left| \mu\left(\hat{I} \cap \Phi^{-1}_n\left(\mathbb{S}^{1} \times Q\right)\right) \cdot \tilde{\mu}\left(J\right) - \mu\left(\hat{I}\right) \cdot \tilde{\mu}\left(Q\right) \right| + \mu\left(\hat{I}\right) \cdot \left| \tilde{\mu}\left(Q\right) - \mu\left(\tilde{S}\right) \right| \\
& \leq \frac{4}{n} \cdot \mu\left(\hat{I}\right) \cdot \mu\left(S\right) + \mu\left(\hat{I}\right) \cdot \frac{14}{n} \cdot \mu\left(S\right) = \frac{18}{n} \cdot \mu\left(\hat{I}\right) \cdot \mu\left(S\right).
\end{align*}
Noting that $\mu\left(S_1\right)= \mu\left(\tilde{S}\right) - 2 \gamma \cdot \tilde{\mu}\left(\tilde{S}_r\right)$ and so $\mu\left(\tilde{S}\right) - \mu\left(S_1\right) \leq 2 \cdot \frac{1}{n \cdot q^{\sigma}_{n}} \cdot \tilde{\mu}\left(\tilde{S}_r\right) \leq \frac{2}{n} \cdot \mu\left(S\right)$ we obtain in the same way as above:
\begin{equation*}
\left| \mu\left(\hat{I} \cap \Phi^{-1}_n\left(\mathbb{S}^{1} \times Q_1\right)\right) \cdot \tilde{\mu}\left(J\right) - \mu\left(\hat{I}\right) \cdot \mu\left(\tilde{S}\right) \right| \leq \frac{20}{n} \cdot \mu\left(\hat{I}\right) \cdot \mu\left(S\right).
\end{equation*}
Using equation (\ref{big}) this yields:
\begin{equation*}
\left| \mu\left(\hat{I} \cap \Phi^{-1}_n\left(g^{-1}_{n}\left(\tilde{S}\right)\right)\right) \cdot \tilde{\mu}\left(J\right) - \mu\left(\hat{I}\right) \cdot \mu\left(\tilde{S}\right) \right| \leq \frac{20}{n} \cdot \mu\left(\hat{I}\right) \cdot \mu\left(S\right).
\end{equation*}
Finally, we conclude with the aid of equation (\ref{triangle}) because of $\varepsilon = \frac{1}{8n^4}$:
\begin{equation*}
\left| \mu\left(\hat{I} \cap \Phi^{-1}_n\left(g^{-1}_{n}(S)\right)\right) \cdot \tilde{\mu}\left(J\right) - \mu\left(\hat{I}\right) \cdot \mu\left(S\right) \right| \leq \frac{21}{n} \cdot \mu\left(\hat{I}\right) \cdot \mu\left(S\right).
\end{equation*}
\end{pr}

Now we are able to prove the desired criterion for weak mixing.

\begin{prop}[Criterion for weak mixing] \label{prop:crit}
Let $f_n = H_n \circ R_{\alpha_{n+1}} \circ H^{-1}_{n}$ and the sequence $\left(m_n\right)_{n \in \mathbb{N}}$ be constructed as in the previous sections. Suppose additionally that $d_0 \left(f^{m_n}, f^{m_n}_{n} \right) < \frac{1}{2^n}$ for every $n \in \mathbb{N}$, $\left\| DH_{n-1} \right\|_0 \leq \frac{\ln\left(q_n\right)}{n}$ and that the limit $f= \lim_{n \rightarrow \infty} f_n$ exists. \\
Then $f$ is weakly mixing.
\end{prop}

\begin{pr}
To apply Lemma \ref{lem:app2} we consider the partial partitions $\nu_n \coloneqq H_{n-1} \circ g_n \left(\eta_n\right)$. As proven in Lemma \ref{lem:points} these partial partitions satisfy $\nu_n \rightarrow \varepsilon$. We have to establish equation (\ref{2}). To do so, let $\varepsilon >0$ and a $m$-dimensional cube $A\subseteq \mathbb{S}^1 \times \left(0,1\right)^{m-1}$ be given. There exists $N \in \mathbb{N}$ such that $A \subseteq \mathbb{S}^1 \times \left[\frac{1}{n^4}, 1-\frac{1}{n^4}\right]^{m-1}$ for every $n \geq N$. Because of Lemma \ref{lem:distri} and the properties of a $\left(\frac{1}{n \cdot q^{m}_{n}}, \frac{1}{n^4}, \frac{1}{n}\right)$-distribution we obtain for every $\hat{I}_n \in \eta_n$ that $\pi_{\vec{r}}\left(\Phi_n\left( \hat{I}_n\right)\right) \supseteq \left[\frac{1}{n^4}, 1-\frac{1}{n^4}\right]^{m-1}$. Furthermore, we note that $f^{m_n}_{n} = H_n \circ R^{m_n}_{\alpha_{n+1}} \circ H^{-1}_{n} = H_{n-1} \circ g_n \circ \Phi_n \circ g^{-1}_{n} \circ H^{-1}_{n-1}$. \\
Let $S_n$ be a $m$-dimensional cube of side length $q^{-\sigma}_{n}$ contained in $\mathbb{S}^1 \times \left[\frac{1}{n^4}, 1-\frac{1}{n^4}\right]^{m-1}$. We look at $C_n \coloneqq H_{n-1}\left(S_n\right)$, $\Gamma_n \in \nu_n$, and compute (since $g_n$ and $H_{n-1}$ are measure-preserving):
\begin{align*}
& \left| \mu \left(\Gamma_n \cap f^{-m_n}_{n}\left(C_n\right)\right) - \mu \left(\Gamma_n\right) \cdot \mu\left(C_n\right) \right| = \left| \mu\left(\hat{I}_n \cap \Phi^{-1}_n \circ g^{-1}_{n}\left(S_n\right)\right) - \mu\left(\hat{I}_n\right) \cdot \mu\left(S_n\right) \right| \\
& \leq \frac{1}{\tilde{\mu}\left(J\right)} \cdot \left| \mu\left(\hat{I}_n \cap \Phi^{-1}_n \circ g^{-1}_{n}\left(S_n\right)\right) \cdot \tilde{\mu}\left(J\right) - \mu\left(\hat{I}_n\right) \cdot \mu\left(S_n\right) \right| + \frac{1-\tilde{\mu}\left(J\right)}{\tilde{\mu}\left(J\right)} \cdot \mu\left(\hat{I}_n\right) \cdot \mu\left(S_n\right).
\end{align*}
Bernoulli's inequality yields: $\tilde{\mu}(J) \geq \left( 1-\frac{1}{n}\right)^{m-1} \geq 1 + \left(m-1\right) \cdot \left(-\frac{1}{n} \right) = 1 - \frac{m-1}{n}$. Hence we obtain for $n>2\cdot(m-1)$: $\tilde{\mu}\left(J\right)\geq\frac{1}{2}$ and so: $\frac{1- \tilde{\mu}\left(J\right)}{\tilde{\mu}\left(J\right)} \leq 2 \cdot \left(1- \tilde{\mu}\left(J\right) \right) \leq \frac{2 \cdot \left(m-1\right)}{n}$. We continue by applying Lemma \ref{lem:cube}:
\begin{align*}
\left| \mu \left(\Gamma_n \cap f^{-m_n}_{n}\left(C_n\right)\right) - \mu \left(\Gamma_n\right) \cdot \mu\left(C_n\right) \right| & \leq 2 \cdot \frac{21}{n} \cdot \mu\left(\hat{I}_n\right) \cdot \mu\left(S_n\right) + \frac{2 \cdot \left(m-1\right)}{n} \cdot \mu\left(\hat{I}_n\right) \cdot \mu\left(S_n\right) \\
& = \frac{40+2\cdot m}{n} \cdot \mu\left(\hat{I}_n\right) \cdot \mu \left( S_n \right)
\end{align*}
Moreover, it holds that diam$\left(C_n\right) \leq \left\|DH_{n-1} \right\|_0 \cdot \text{diam}\left(S_n\right) \leq \sqrt{m} \cdot \frac{\ln\left(q_n\right)}{q^{\sigma}_{n}}$, i.e. diam$\left(C_n\right)\rightarrow0$ as $n\rightarrow \infty$. Thus, we can approximate $A$ by a countable disjoint union of sets $C_n=H_{n-1}\left(S_n\right)$ with $S_n \subseteq \mathbb{S}^1 \times \left[\frac{1}{n^4}, 1-\frac{1}{n^4}\right]^{m-1}$ a $m$-dimensional cube of side length $q^{-\sigma}_{n}$ with given precision, assuming that $n$ is chosen to be large enough. Consequently for sufficiently large  $n$ there are sets $A_1=\dot{\bigcup}_{i \in \Sigma^{1}_{n}} C^{i}_{n}$ and $A_2=\dot{\bigcup}_{i \in \Sigma^{2}_{n}} C^{i}_{n}$ with countable sets $\Sigma^{1}_{n}$ and $\Sigma^{2}_{n}$ of indices satisfying $A_1 \subseteq A \subseteq A_2$ as well as $ \left| \mu(A) - \mu(A_i) \right| \leq \frac{\epsilon}{3} \cdot \mu(A)$ for $i=1,2$. \\
Additionally we choose $n$ such that $\frac{40 + 2 \cdot m}{n} < \frac{\epsilon}{3}$ holds. It follows that
\begin{align*}
& \mu \left(\Gamma_n \cap f^{-m_n}_{n}\left(A\right)\right) - \mu \left(\Gamma_n\right) \cdot \mu\left(A\right) \\
& \leq \mu \left(\Gamma_n \cap f^{-m_n}_{n}\left(A_2\right)\right) - \mu \left(\Gamma_n\right) \cdot \mu\left(A_2\right) + \mu \left(\Gamma_n\right) \cdot \left( \mu\left(A_2\right) - \mu\left(A\right) \right) \\
& \leq \sum_{i \in \Sigma^{2}_{n}} \left(\mu \left(\Gamma_n \cap f^{-m_n}_{n}\left(C^{i}_{n}\right)\right) - \mu \left(\Gamma_n\right) \cdot \mu\left(C^{i}_{n}\right) \right) + \frac{\epsilon}{3} \cdot \mu \left(\Gamma_n\right) \cdot \mu\left(A\right) \\
& \leq \sum_{i \in \Sigma^{2}_{n}} \left(\frac{40 + 2 \cdot m}{n} \cdot \mu\left(\hat{I}_n\right) \cdot \mu\left(S^{i}_{n}\right) \right) + \frac{\epsilon}{3} \cdot \mu \left(\Gamma_n\right) \cdot \mu\left(A\right) \\
& = \frac{40 + 2 \cdot m}{n} \cdot \mu \left(\Gamma_n\right) \cdot \mu\left(\bigcup_{i \in \Sigma^{2}_{n}} C^{i}_{n}\right) + \frac{\epsilon}{3} \cdot \mu \left(\Gamma_n\right) \cdot \mu\left(A\right) \leq \frac{\epsilon}{3} \cdot \mu \left(\Gamma_n\right) \cdot \mu\left(A_2\right) + \frac{\epsilon}{3} \cdot \mu \left(\Gamma_n\right) \cdot \mu\left(A\right) \\
&  = \frac{\epsilon}{3} \cdot \mu \left(\Gamma_n\right) \cdot \mu\left(A\right) + \frac{\epsilon}{3} \cdot \mu \left(\Gamma_n\right) \cdot \left(\mu\left(A_2\right)-\mu\left(A\right)\right) + \frac{\epsilon}{3} \cdot \mu \left(\Gamma_n\right) \cdot \mu\left(A\right) \leq \epsilon \cdot \mu \left(\Gamma_n\right) \cdot \mu\left(A\right).
\end{align*}
Analogously, we estimate that $\mu \left(\Gamma_n \cap f^{-m_n}_{n}\left(A\right)\right) - \mu \left(\Gamma_n\right) \cdot \mu\left(A\right) \geq - \epsilon \cdot \mu \left(\Gamma_n\right) \cdot \mu\left(A\right)$. Both estimates enable us to conclude that $\left| \mu \left(\Gamma_n \cap f^{-m_n}_{n}\left(A\right)\right) - \mu \left(\Gamma_n\right) \cdot \mu\left(A\right) \right| \leq \epsilon \cdot \mu \left(\Gamma_n\right) \cdot \mu\left(A\right)$.
\end{pr}

\section{Convergence of $\left(f_n\right)_{n \in \mathbb{N}}$ in Diff$^{\infty}\left(M\right)$} \label{section:conv}
In the following we show that the sequence of constructed measure-preserving smooth diffeomorphisms $f_n = H_n \circ R_{\alpha_{n+1}} \circ H^{-1}_{n}$ converges. For this purpose, we need a couple of results concerning the conjugation maps.

\subsection{Properties of the conjugation maps $\phi_n$ and $H_n$}
In order to find estimates on the norms $\left|\left\| H_n \right\|\right|_k$ we will need the next technical result which is an application of the chain rule:

\begin{lem} \label{lem:derphi}
Let $\phi := \tilde{\phi}^{(m)}_{\lambda_m, \mu_m} \circ ... \circ \tilde{\phi}^{(2)}_{\lambda_2, \mu_2}$, $j \in \left\{1,...,m\right\}$ and $k \in \mathbb{N}$. For any multi-index $\vec{a}$ with $\left|\vec{a}\right|=k$ the partial derivative $D_{\vec{a}} \left[ \phi \right]_j$ consists of a sum of products of at most $(m-1) \cdot k$ terms of the form 
\begin{equation*}
D_{\vec{b}} \left(\left[\tilde{\phi}^{(i)}_{\lambda_i, \mu_i}\right]_l\right) \circ \tilde{\phi}^{(i-1)}_{\lambda_{i-1}, \mu_{i-1}} \circ ... \circ \tilde{\phi}^{(2)}_{\lambda_2, \mu_2},
\end{equation*}
where $l \in \left\{1,...,m\right\}$, $i \in \left\{2,...,m\right\}$ and $\vec{b}$ is a multi-index with $\left|\vec{b}\right|\leq k$.
\end{lem}

In the same way we obtain a similar statement holding for the inverses:
\begin{lem} \label{lem:derphiinv}
Let $\psi := \left(\tilde{\phi}^{(2)}_{\lambda_2, \mu_2}\right)^{-1} \circ ... \circ \left(\tilde{\phi}^{(m)}_{\lambda_m, \mu_m}\right)^{-1}$, $j \in \left\{1,...,m\right\}$ and $k \in \mathbb{N}$. For any multi-index $\vec{a}$ with $\left|\vec{a}\right|=k$ the partial derivative $D_{\vec{a}} \left[ \psi \right]_j$ consists of a sum of products of at most $(m-1) \cdot k$ terms of the following form 
\begin{equation*}
D_{\vec{b}} \left(\left[\left(\tilde{\phi}^{(i)}_{\lambda_i, \mu_i}\right)^{-1}\right]_l\right) \circ \left(\tilde{\phi}^{(i+1)}_{\lambda_{i+1}, \mu_{i+1}}\right)^{-1} \circ ... \circ \left(\tilde{\phi}^{(m)}_{\lambda_m, \mu_m}\right)^{-1},
\end{equation*}
where $l \in \left\{1,...,m\right\}$, $i \in \left\{2,...,m\right\}$ and $\vec{b}$ is a multi-index with $\left|\vec{b}\right|\leq k$.
\end{lem}

\begin{rem} \label{rem:faa}
In the proof of the following lemmas we will use the formula of Faà di Bruno in several variables. It can be found in the paper \textit{``A multivariate Faà di Bruno formula with applications''} (\cite{Fa}) for example. \\
For this we introduce an ordering on $\mathbb{N}^d_0$: For multiindices $\vec{\mu} = \left(\mu_1,...,\mu_d\right)$ and $\vec{\nu} = \left(\nu_1,...,\nu_d\right)$ in $\mathbb{N}^d_0$ we will write $\vec{\mu} \prec \vec{\nu}$, if one of the following properties is satisfied:
\begin{enumerate}
	\item $\left| \vec{\mu} \right| < \left| \vec{\nu} \right|$, where $\left| \vec{\mu} \right| = \sum^{d}_{i=1} \mu_i$.
	\item $\left| \vec{\mu} \right| = \left| \vec{\nu} \right|$ and $\mu_1 < \nu_1$.
	\item $\left| \vec{\mu} \right| = \left| \vec{\nu} \right|$, $\mu_i = \nu_i$ for $1\leq i \leq k$ and $\mu_{k+1} < \nu_{k+1}$ for some $1 \leq k < d$.
\end{enumerate}
In other words, we compare by order and then lexicographically.
Additionally we will use these notations:
\begin{itemize}
	\item For $\vec{\nu}=\left(\nu_1,...,\nu_d\right) \in \mathbb{N}^d_0$:
	\begin{equation*}
	\vec{\nu}! = \prod^{d}_{i=1} \nu_i !
	\end{equation*}
	\item For $\vec{\nu}=\left(\nu_1,...,\nu_d\right) \in \mathbb{N}^d_0$ and $\vec{z} = \left(z_1,...,z_d\right) \in \mathbb{R}^d$:
	\begin{equation*}
	\vec{z}^{\ \vec{\nu}} = \prod^{d}_{i=1} z^{\nu_i}_{i}
	\end{equation*}
\end{itemize}
Then we get for the composition $h\left(x_1,...,x_d\right) := f\left(g^{(1)}\left(x_1,...,x_d\right),..., g^{(m)}\left(x_1,...,x_d\right)\right)$ with sufficiently differentiable functions $f: \mathbb{R}^m \rightarrow \mathbb{R}$, $g^{(i)}: \mathbb{R}^d \rightarrow \mathbb{R}$ and a multi-index $\vec{\nu} \in \mathbb{N}^d_0$ with $\left| \vec{\nu} \right| = n$:
\begin{equation*}
D_{\vec{\nu}}h = \sum_{\vec{\lambda} \in \mathbb{N}^m_0 \text{ with } 1\leq \left|\vec{\lambda} \right|\leq n} D_{\vec{\lambda}}f \cdot \sum^{n}_{s=1}\ \  \sum_{p_s\left(\vec{\nu},\vec{\lambda}\right)} \vec{\nu}! \cdot \prod^{s}_{j=1}\frac{\left[D_{\vec{l}_j} \vec{g}\right]^{\vec{k}_j}}{\vec{k}_j ! \cdot \left(\vec{l}_j !\right)^{\left|\vec{k}_j\right|}}
\end{equation*}
Here $\left[D_{\vec{l}_j} \vec{g}\right]$ denotes $\left(D_{\vec{l}_j} g^{(1)},...,D_{\vec{l}_j}g^{(m)}\right)$ and
\begin{align*}
& p_s\left(\vec{\nu}, \vec{\lambda}\right) := \\
& \left\{ \left(\vec{k}_1,...,\vec{k}_s, \vec{l}_1,...,\vec{l}_s\right) : \vec{k}_i \in \mathbb{N}^{m}_{0}, \left| \vec{k}_i \right| >0, \vec{l}_i \in \mathbb{N}^d_0, 0 \prec \vec{l}_1 \prec ... \prec \vec{l}_s, \sum^{s}_{i=1} \vec{k}_i = \vec{\lambda} \text{ and } \sum^{s}_{i=1} \left| \vec{k}_i \right| \cdot \vec{l}_i = \vec{\nu} \right\}
\end{align*}
\end{rem}
With the aid of these technical results we can prove an estimate on the norms of the map $\phi_n$:

\begin{lem} \label{lem:normphi}
For every $k \in \mathbb{N}$ it holds that
\begin{equation*}
||| \phi_n |||_k \leq C \cdot q^{\left(m-1\right)^2 \cdot k \cdot n \cdot \left(n+1\right)}_{n},
\end{equation*}
where $C$ is a constant depending on $m$, $k$ and $n$, but is independent of $q_n$.
\end{lem}

\begin{pr}
First of all we consider the map $\tilde{\phi}_{\lambda, \mu} \coloneqq \tilde{\phi}_{\lambda, \varepsilon, i,j,\mu, \delta, \varepsilon_2} = C^{-1}_{\lambda} \circ \psi_{\mu,\delta,i,j,\varepsilon_2} \circ \varphi_{\varepsilon, i, j} \circ C_{\lambda}$ introduced in subsection \ref{subsection:phi}:
	\begin{align*}
	& \tilde{\phi}_{\lambda, \mu}\left(x_1,...,x_m\right) = \\
	& \left(\frac{1}{\lambda} \left[\psi_{\mu} \circ \varphi_{\varepsilon}\right]_1 \left(\lambda x_1,x_2,...,x_m\right), \left[\psi_{\mu} \circ \varphi_{\varepsilon}\right]_2\left(\lambda x_1,x_2,...,x_m\right),...,\left[\psi_{\mu} \circ \varphi_{\varepsilon}\right]_m\left(\lambda x_1,x_2,...,x_m\right) \right).
	\end{align*}
Let $k \in \mathbb{N}$. We compute for a multi-index $\vec{a}$ with $0\leq \left| \vec{a}\right| \leq k$: $\left\|D_{\vec{a}}\left[\tilde{\phi}_{\lambda, \mu}\right]_1\right\|_0 \leq \lambda^{k-1} \cdot ||| \psi_{\mu} \circ \varphi_{\varepsilon} |||_k$ and for $r \in \left\{2,...,m\right\}$: $\left\|D_{\vec{a}}\left[\tilde{\phi}_{\lambda, \mu}\right]_r\right\|_0 \leq \lambda^{k} \cdot ||| \psi_{\mu} \circ \varphi_{\varepsilon} |||_k$. \\
Therefore, we examine the map $\psi_{\mu}$. For any multi-index $\vec{a}$ with $0\leq \left| \vec{a}\right| \leq k$ and $r \in \left\{1,...,m\right\}$ we obtain: $\left\|D_{\vec{a}} \left[\psi_{\mu}\right]_r \right\|_0 \leq \mu^{k-1} \cdot ||| \varphi_{\varepsilon_2} |||_k = C_{k, \varepsilon_2} \cdot \mu^{k-1}$ and analogously $\left\|D_{\vec{a}} \left[\psi^{-1}_{\mu}\right]_r \right\|_0 \leq C_{k, \varepsilon_2} \cdot \mu^{k-1}$. Hence: $|||\psi_{\mu}|||_k \leq C \cdot \mu^{k-1}$. \\
In the next step we use the formula of Faà di Bruno mentioned in remark \ref{rem:faa}. With it we compute for any multi-index $\vec{\nu}$ with $\left|\vec{\nu}\right| = k$: 
\begin{align*}
\left\|D_{\vec{\nu}} \left[\left(\psi_{\mu} \circ \varphi_{\varepsilon}\right)^{-1}\right]_r \right\|_0 & = \left\|D_{\vec{\nu}} \left[ \varphi^{-1}_{\varepsilon} \circ \psi^{-1}_{\mu}\right]_r \right\|_0 \\
& = \left\| \sum_{\vec{\lambda} \in \mathbb{N}^m_0 , 1\leq \left|\vec{\lambda} \right|\leq k} D_{\vec{\lambda}}\left[\varphi^{-1}_{\varepsilon}\right]_r \sum^{k}_{s=1}\ \  \sum_{\left(\vec{k}_1,...,\vec{k}_s, \vec{l}_1,...,\vec{l}_s\right) \in p_s\left(\vec{\nu},\vec{\lambda}\right)} \vec{\nu}! \prod^{s}_{j=1}\frac{\left[D_{\vec{l}_j} \psi^{-1}_{\mu} \right]^{\vec{k}_j}}{\vec{k}_j ! \cdot \left(\vec{l}_j !\right)^{\left|\vec{k}_j\right|}} \right\|_0 \\ 
& = \left\| \sum_{\vec{\lambda} \in \mathbb{N}^m_0, 1\leq \left|\vec{\lambda} \right|\leq k} D_{\vec{\lambda}}\left[\varphi^{-1}_{\varepsilon}\right]_r \cdot \sum^{k}_{s=1}\ \  \sum_{p_s\left(\vec{\nu},\vec{\lambda}\right)} \vec{\nu}! \cdot \prod^{s}_{j=1}\frac{\prod^{m}_{t=1}\left(D_{\vec{l}_j} \left[ \psi^{-1}_{\mu}\right]_t\right)^{\vec{k}_{j_t}}}{\vec{k}_j ! \cdot \left(\vec{l}_j !\right)^{\left|\vec{k}_j\right|}} \right\|_0 \\ \displaybreak[0]
& \leq \sum_{\vec{\lambda} \in \mathbb{N}^m_0, 1\leq \left|\vec{\lambda} \right|\leq k} \left\|D_{\vec{\lambda}} \left[\varphi^{-1}_{\varepsilon}\right]_r \right\|_0 \cdot \sum^{k}_{s=1} \ \ \sum_{p_s\left(\vec{\nu},\vec{\lambda}\right)}\vec{\nu}! \cdot \prod^{s}_{j=1}\frac{\prod^{m}_{t=1}\left\| D_{\vec{l}_j} \left[\psi^{-1}_{\mu}\right]_t \right\|^{\vec{k}_{j_t}}_0}{\vec{k}_j ! \cdot \left(\vec{l}_j !\right)^{\left|\vec{k}_j\right|}} \\ \displaybreak[0]
& \leq \sum_{\vec{\lambda} \in \mathbb{N}^m_0 \text{ with } 1\leq \left|\vec{\lambda} \right|\leq k} \left\|D_{\vec{\lambda}} \left[\varphi^{-1}_{\varepsilon}\right]_r \right\|_0 \cdot \sum^{k}_{s=1} \ \ \sum_{p_s\left(\vec{\nu},\vec{\lambda}\right)}\vec{\nu}! \cdot \prod^{s}_{j=1}\frac{|||\psi^{-1}_{\mu}|||^{\sum^{m}_{t=1} \vec{k}_{j_t}}_{\left|\vec{l}_j\right|}}{\vec{k}_j ! \cdot \left(\vec{l}_j !\right)^{\left|\vec{k}_j\right|}} \\ 
& = \sum_{\vec{\lambda} \in \mathbb{N}^m_0 \text{ with } 1\leq \left|\vec{\lambda} \right|\leq k} \left\|D_{\vec{\lambda}} \left[\varphi^{-1}_{\varepsilon}\right]_r \right\|_0 \cdot \sum^{k}_{s=1} \ \ \sum_{p_s\left(\vec{\nu},\vec{\lambda}\right)}\vec{\nu}! \cdot \prod^{s}_{j=1}\frac{|||\psi^{-1}_{\mu}|||^{\left|\vec{k}_j\right|}_{\left|\vec{l}_j\right|}}{\vec{k}_j ! \cdot \left(\vec{l}_j !\right)^{\left|\vec{k}_j\right|}}
\end{align*}
As seen above: $|||\psi^{-1}_{\mu}|||^{\left|\vec{k}_j\right|}_{\left|\vec{l}_j\right|} \leq C \cdot \mu^{\left|\vec{k}_j\right| \cdot \left|\vec{l}_j\right|}$. Hereby: $\prod^{s}_{j=1} |||\psi^{-1}_{\mu} |||^{\left| \vec{k}_j \right|}_{\left|\vec{l}_j\right|} \leq \hat{C} \cdot \mu^{\sum^{s}_{j=1} \left|\vec{l}_j\right| \cdot \left|\vec{k}_j\right|}$,
where $\hat{C}$ is independent of $\mu$. By definition of the set $p_s\left(\vec{\nu}, \vec{\lambda}\right)$ we have $\sum^{s}_{i=1} \left| \vec{k}_i \right| \cdot \vec{l}_i = \vec{\nu}$. Hence:
\begin{equation*}
k = \left| \vec{\nu} \right| = \left| \sum^{s}_{i=1} \left| \vec{k}_i\right| \cdot \vec{l}_i \right| = \sum^{m}_{t=1} \left(\sum^{s}_{i=1} \left| \vec{k}_i\right| \cdot \vec{l}_i \right)_t = \sum^{m}_{t=1} \sum^{s}_{i=1} \left| \vec{k}_i \right| \cdot \vec{l}_{i_t} = \sum^{s}_{i=1} \left| \vec{k}_i \right| \cdot \left( \sum^{m}_{t=1} \vec{l}_{i_t} \right) = \sum^{s}_{i=1} \left| \vec{k}_i \right| \cdot \left| \vec{l}_i \right|
\end{equation*}
This shows $\prod^{s}_{j=1} |||\psi^{-1}_{\mu} |||^{\left| \vec{k}_j \right|}_{\left|\vec{l}_j\right|} \leq \hat{C} \cdot \mu^k$ and finally $\left\|D_{\vec{\nu}} \left[\left(\psi_{\mu} \circ \varphi_{\varepsilon}\right)^{-1}\right]_r \right\|_0 \leq C \cdot \mu^k$. Analogously we compute $\left\|D_{\vec{\nu}} \left[\psi_{\mu} \circ \varphi_{\varepsilon}\right]_r \right\|_0 \leq C \cdot |||\psi_{\mu} |||_k \leq C \cdot \mu^{k-1}$. Altogether, we obtain $|||\psi_{\mu} \circ \varphi_{\varepsilon}|||_k \leq C \cdot \mu^k$. Hereby, we estimate $\left\| D_{\vec{a}} \left[\tilde{\phi}_{\lambda,  \mu} \right]_r \right\|_0 \leq C \cdot \lambda^k \cdot \mu^k$ and analogously $\left\| D_{\vec{a}} \left[\tilde{\phi}^{-1}_{\lambda,  \mu} \right]_r \right\|_0 \leq C \cdot \lambda^k \cdot \mu^k$. In conclusion this yields $||| \tilde{\phi}_{\lambda,  \mu} |||_k \leq C \cdot \mu^k \cdot \lambda^k$. \\
In the next step we consider $\phi \coloneqq \tilde{\phi}^{(m)}_{\lambda_m, \mu_m} \circ ... \circ \tilde{\phi}^{(2)}_{\lambda_2, \mu_2}$. Let $\lambda_{\text{max}} \coloneqq \max \left\{ \lambda_2,..., \lambda_m \right\}$ as well as $\mu_{\text{max}} \coloneqq \max \left\{ \mu_2,..., \mu_m \right\}$. Inductively we will show $||| \phi |||_k \leq \tilde{C} \cdot \lambda^{(m-1) \cdot k}_{\text{max}} \cdot \mu^{(m-1) \cdot k}_{\text{max}}$ for every $k \in \mathbb{N}$, where $\tilde{C}$ is a constant independent of $\lambda_i$ and $\mu_i$. \\
\textit{Start: $k=1$ } \\
Let $l \in \left\{1,...,m\right\}$ be arbitrary. By Lemma \ref{lem:derphi} a partial derivative of $\left[ \phi \right]_l$ of first order consists of a sum of products of at most $m-1$ first order partial derivatives  of functions $\tilde{\phi}^{(j)}_{\lambda_j, \mu_j}$. Therewith, we obtain using $|||\tilde{\phi}^{(j)}_{\lambda_j, \mu_j}|||_1 \leq C \cdot \lambda_{\text{max}} \cdot \mu_{\text{max}}$ the estimate $\left\| D_i \left[ \phi \right]_l \right\|_0 \leq C_1 \cdot \lambda^{m-1}_{\text{max}} \cdot \mu^{m-1}_{\text{max}}$ for every $i \in \left\{1,...,m\right\}$, where $C_1$ is a constant independent of $\lambda$ and $\mu$.\\
With the aid of Lemma \ref{lem:derphiinv}	we obtain the same statement for $\phi^{-1} = \left(\tilde{\phi}^{(2)}_{\lambda_2, \mu_2}\right)^{-1} \circ ... \circ \left(\tilde{\phi}^{(m)}_{\lambda_m, \mu_m} \right)^{-1}$. Hence, we conclude: $||| \phi |||_1 \leq \tilde{C}_1 \cdot \lambda^{m-1}_{\text{max}} \cdot \mu^{m-1}_{\text{max}}$. \\
\textit{Assumption: } The claim is true for $k \in \mathbb{N}$. \\
\textit{Induction step $k \rightarrow k+1$:}\\
In the proof of Lemma \ref{lem:derphi} one observes that at the transition $k\rightarrow k+1$ in the product of at most $(m-1) \cdot k$ terms of the form $D_{\vec{b}} \left(\left[\tilde{\phi}^{(i)}_{\lambda_i, \mu_i}\right]_l\right) \circ \tilde{\phi}^{(i-1)}_{\lambda_{i-1}, \mu_{i-1}} \circ ... \circ \tilde{\phi}^{(2)}_{\lambda_2, \mu_2}$ one is replaced by a product of a term $\left(D_j D_{\vec{b}} \left[\tilde{\phi}^{(i)}_{\lambda_i, \mu_i}\right]_l \right) \circ \tilde{\phi}^{(i-1)}_{\lambda_{i-1}, \mu_{i-1}} \circ ... \circ \tilde{\phi}^{(2)}_{\lambda_2, \mu_2}$ with $j \in \left\{1,...,m\right\}$ and at most $m-2$ partial derivatives of first order. Because of $||| \tilde{\phi}^{(i)}_{\lambda_i, \mu_i}|||_{k+1} \leq C \cdot \lambda^{k+1}_{\text{max}} \cdot \mu^{k+1}_{\text{max}}$ and $||| \tilde{\phi}^{(j)}_{\lambda_j, \mu_j}|||_1 \leq C \cdot \lambda_{\text{max}} \cdot \mu_{\text{max}}$ the $\lambda_{\text{max}}$-exponent as well as the $\mu_{\text{max}}$-exponent increase by at most $1+(m-2) \cdot 1= m-1$. \\
In the same spirit one uses the proof of Lemma \ref{lem:derphiinv} to show that also in case of $\phi^{-1}$ the $\lambda_{\text{max}}$-exponent as well as the $\mu_{\text{max}}$-exponent increase by at most $m-1$. \\
Using the assumption we conclude 
\begin{equation*}
||| \phi |||_{k+1} \leq \hat{C} \cdot \lambda^{k \cdot (m-1) + m-1}_{\text{max}} \cdot \mu^{k \cdot (m-1) + m-1}_{\text{max}} = \hat{C}\cdot \lambda^{(k+1) \cdot (m-1)}_{\text{max}} \cdot \mu^{(k+1) \cdot (m-1)}_{\text{max}}.
\end{equation*}
So the proof by induction is completed. \\
In the setting of our explicit construction of the map $\phi_n$ in section \ref{subsection:phi} we have $\varepsilon_1 = \frac{1}{60 \cdot n^4}$, $\varepsilon_2 = \frac{1}{22 \cdot n^4}$, $\lambda_{\text{max}} =n \cdot q^{1+(m-1) \cdot \frac{n \cdot \left(n-1\right)}{2} + (m-2) \cdot n}_{n}$ and $\mu_{\text{max}} = q^{n}_n$. Thus: 
\begin{align*}
||| \phi_n |||_k & \leq \tilde{C}\left(m,k,n\right) \cdot \left(n \cdot q^{1+(m-1) \cdot \frac{n \cdot \left(n-1\right)}{2} + (m-2) \cdot n}_{n}\right)^{(m-1) \cdot k} \cdot \left(q^{n}_{n}\right)^{(m-1) \cdot k} \\
& \leq C\left(m,k,n\right) \cdot q^{\left(m-1\right)^2 \cdot k \cdot n \cdot \left(n+1\right)}_{n},
\end{align*}	
where $C\left(m,k,n\right)$ is a constant independent of $q_n$.
\end{pr}

In the next step we consider the map $h_n = g_n \circ \phi_n$, where $g_n$ is constructed in section \ref{subsection:g}:

\begin{lem} \label{lem:normh}
For every  $k \in \mathbb{N}$ it holds:
\begin{equation*}
||| h_n |||_k \leq \bar{C} \cdot q^{3 \cdot \left(m-1\right)^2 \cdot k \cdot n \cdot \left(n+1\right)}_{n},
\end{equation*}
where $\bar{C}$ is a constant depending on $m$, $k$ and $n$, but is independent of $q_n$.
\end{lem}

\begin{pr}
Outside of $\mathbb{S}^1 \times \left[\delta,1-\delta\right]^{m-1}$, i.e. $g_n=\tilde{g}_{\left[nq^{\sigma}_n\right]}$, we have:
\begin{align*}
& h_n \left(x_1,...,x_m\right) = g_n \circ \phi_n \left(x_1,...,x_m\right) \\
& = \left( \left[\phi_n \left(x_1,...,x_m\right)\right]_1 + \left[ n \cdot q^{\sigma}_n \right] \cdot \left[\phi_n \left(x_1,...,x_m\right)\right]_2, \left[\phi_n \left(x_1,...,x_m\right)\right]_2,...,\left[\phi_n \left(x_1,...,x_m\right)\right]_m \right)
\end{align*}
and
\begin{align*}
& h^{-1}_n \left(x_1,...,x_m\right) = \phi^{-1}_n \circ g^{-1}_n \left(x_1,...,x_m\right) \\
& = \left( \left[\phi^{-1}_n \left(x_1-\left[ n \cdot q^{\sigma}_n \right] \cdot x_2,x_2,...,x_m\right)\right]_1,...,\left[\phi_n \left(x_1-\left[ n \cdot q^{\sigma}_n \right] \cdot x_2,x_2,...,x_m\right)\right]_m \right).
\end{align*}
Since $\sigma < 1$ we can estimate:
\begin{equation*}
||| h_n |||_k \leq 2 \cdot \left[n \cdot q^{\sigma}_{n}\right]^k \cdot ||| \phi_n |||_k  \leq \bar{C}\left(m,k,n\right) \cdot q^{\sigma \cdot k}_{n} \cdot q^{\left(m-1\right)^2 \cdot k \cdot n \cdot \left(n+1\right)}_{n} \leq \bar{C}\left(m,k,n\right) \cdot q^{2 \cdot \left(m-1\right)^2 \cdot k \cdot n \cdot \left(n+1\right)}_{n}
\end{equation*}
with a constant $\bar{C}\left(m,k,n\right)$ independent of $q_n$. \\
In the other case we have
\begin{equation*}
g_n \circ \phi_n \left(x_1,...,x_m\right) = \left( \left[g_{a,b,\varepsilon} \left(\left[\phi_n\right]_1,\left[\phi_n\right]_2\right)\right]_1, \left[g_{a,b,\varepsilon} \left(\left[\phi_n\right]_1,\left[\phi_n\right]_2\right)\right]_2,\left[\phi_n \right]_3,...,\left[\phi_n \right]_m \right).
\end{equation*}
We will use the formula of Faà di Bruno as above for any multi-index $\vec{\nu}$ with $\left|\vec{\nu}\right|=k$ and $r\in\left\{1,...,m\right\}$:
\begin{align*}
\left\|D_{\vec{\nu}} \left[h_n\right]_r \right\|_0 & = \left\|D_{\vec{\nu}} \left[g_{a,b,\varepsilon} \circ \phi_n \right]_r \right\|_0 \\
& \leq \sum_{\vec{\lambda} \in \mathbb{N}^m_0 \text{ with } 1\leq \left|\vec{\lambda} \right|\leq k} \left\|D_{\vec{\lambda}} \left[g_{a,b,\varepsilon}\right]_r \right\|_0 \cdot \sum^{k}_{s=1} \ \ \sum_{p_s\left(\vec{\nu},\vec{\lambda}\right)}\vec{\nu}! \cdot \prod^{s}_{j=1}\frac{||| \phi_n |||^{\left|\vec{k}_j\right|}_{\left|\vec{l}_j\right|}}{\vec{k}_j ! \cdot \left(\vec{l}_j !\right)^{\left|\vec{k}_j\right|}}
\end{align*}
By Lemma \ref{lem:normphi} we have $||| \phi_n |||_k \leq C \cdot q^{\left(m-1\right)^2 \cdot k \cdot n \cdot \left(n+1\right)}_{n}$, where $C$ is a constant independent of $q_n$. As above we show $\prod^{s}_{j=1} ||| \phi_n |||^{\left| \vec{k}_j \right|}_{\left|\vec{l}_j\right|} \leq \hat{C} \cdot q^{\left(\sum^{s}_{j=1} \left|\vec{l}_j\right| \cdot \left|\vec{k}_j\right| \right) \cdot \left(m-1\right)^2 \cdot n \cdot \left(n+1\right)}_{n} = \hat{C} \cdot q^{\left(m-1\right)^2 \cdot k \cdot n \cdot \left(n+1\right)}_{n}$, where $\hat{C}$ is a constant independent of $q_n$.\\
Furthermore, we examine the map $g_{a,b, \varepsilon,\delta}=D^{-1}_{a,b, \varepsilon} \circ g_{\varepsilon} \circ D_{a,b,\varepsilon}$ for $a,b \in \mathbb{Z}$ and obtain 
\begin{equation*}
||| g_{a,b, \varepsilon,\delta} |||_k \leq \left(\frac{b \cdot a}{\varepsilon}\right)^{k} \cdot ||| g_{\varepsilon} |||_k = C_{\varepsilon,k} \cdot b^k \cdot a^{k}.
\end{equation*}
By our constructions in section \ref{subsection:g} we have $b = \left[ n \cdot q^{\sigma}_{n} \right] \leq n \cdot q^{\sigma}_{n}$, $a\leq n \cdot q^{1+(m-1) \cdot \frac{n \cdot \left(n+1\right)}{2}}_{n}$ and $\varepsilon = \frac{1}{8n^4}$. Hence: $||| g_n |||_k \leq C_{n,k}\cdot q^{\sigma \cdot k}_{n} \cdot q^{k+ k \cdot (m-1) \cdot \frac{n \cdot \left(n+1\right)}{2}}_{n} \leq C_{n,k} \cdot q^{2 \cdot k \cdot (m-1) \cdot n \cdot \left(n+1\right)}_{n}$. Finally, we conclude: $\left\|D_{\vec{\nu}} \left[h_n\right]_r \right\|_0 \leq C \cdot q^{2 \cdot k \cdot (m-1) \cdot n \cdot \left(n+1\right)}_{n} \cdot q^{k \cdot \left(m-1\right)^2 \cdot  n \cdot \left(n+1\right)}_{n} \leq C \cdot q^{3 \cdot k \cdot \left(m-1\right)^2 \cdot n \cdot \left(n+1\right)}_{n}$. \\
In the next step we consider $h^{-1}_{n}=\phi^{-1}_{n} \circ g^{-1}_{a,b,\varepsilon}$. For $r \in \left\{1,...,m\right\}$ and any multi-index $\vec{\nu}$ with $\left| \vec{\nu} \right| = k$ we obtain using the formula of Faà di Bruno again:
\begin{align*}
\left\|D_{\vec{\nu}} \left[h^{-1}_n\right]_r \right\|_0 & = \left\|D_{\vec{\nu}} \left[\phi^{-1}_n \circ g^{-1}_n \right]_r \right\|_0 \\
& \leq \sum_{\vec{\lambda} \in \mathbb{N}^m_0 \text{ with } 1\leq \left|\vec{\lambda} \right|\leq k} \left\|D_{\vec{\lambda}} \left[\phi^{-1}_n\right]_r \right\|_0 \cdot \sum^{k}_{s=1} \ \ \sum_{p_s\left(\vec{\nu},\vec{\lambda}\right)}\vec{\nu}! \cdot \prod^{s}_{j=1}\frac{||| g_n |||^{\left|\vec{k}_j\right|}_{\left|\vec{l}_j\right|}}{\vec{k}_j ! \cdot \left(\vec{l}_j !\right)^{\left|\vec{k}_j\right|}}
\end{align*}
As above we show $\prod^{s}_{j=1} ||| g_n |||^{\left| \vec{k}_j \right|}_{\left|\vec{l}_j\right|} \leq \hat{C} \cdot q^{2 \cdot k \cdot (m-1) \cdot n \cdot \left(n+1\right)}_{n}$, where $\hat{C}$ is a constant independent of $q_n$. Since $||| \phi_n |||_k \leq C \cdot q^{k \cdot \left(m-1\right)^2 \cdot n \cdot \left(n+1\right)}_{n}$ we get
\begin{equation*}
\left\|D_{\vec{\nu}} \left[h^{-1}_n\right]_r \right\|_0 \leq \check{C} \cdot q^{2 \cdot k \cdot (m-1) \cdot n \cdot \left(n+1\right)}_{n} \cdot q^{k \cdot \left(m-1\right)^2 \cdot n \cdot \left(n+1\right)}_{n} \leq \check{C}  \cdot q^{3 \cdot k  \cdot \left(m-1\right)^2 \cdot n \cdot \left(n+1\right)}_{n},
\end{equation*}
where $\check{C}$ is a constant independent of $q_n$. \\
Thus, we finally obtain $||| h_n |||_k \leq C(n,k,m) \cdot q^{3 \cdot \left(m-1\right)^2 \cdot k \cdot n \cdot \left(n+1\right)}_{n}$.
\end{pr}

Finally, we are able to prove an estimate on the norms of the map $H_n$:

\begin{lem} \label{lem:normH}
For every $k \in \mathbb{N}$ we get:
\begin{equation*}
||| H_n |||_k \leq \breve{C} \cdot q^{3 \cdot \left(m-1\right)^2 \cdot k \cdot n \cdot \left(n+1\right)}_{n},
\end{equation*}
where $\breve{C}$ is a constant depending solely on $m$, $k$, $n$ and $H_{n-1}$. Since $H_{n-1}$ is independent of $q_n$ in particular, the same is true for $\breve{C}$.
\end{lem}

\begin{pr}
Let $k \in \mathbb{N}$, $r \in \left\{1,...,m\right\}$ and $\vec{\nu} \in \mathbb{N}^{m}_{0}$ be a multi-index with $\left| \vec{\nu} \right| =k$. As above we estimate:
\begin{align*}
\left\|D_{\vec{\nu}} \left[H_n\right]_r \right\|_0 & = \left\|D_{\vec{\nu}} \left[H_{n-1} \circ h_n \right]_r \right\|_0 \\
& \leq \sum_{\vec{\lambda} \in \mathbb{N}^m_0 \text{ with } 1\leq \left|\vec{\lambda} \right|\leq k} \left\|D_{\vec{\lambda}} \left[H_{n-1}\right]_r \right\|_0 \cdot \sum^{k}_{s=1} \ \ \sum_{p_s\left(\vec{\nu},\vec{\lambda}\right)}\vec{\nu}! \cdot \prod^{s}_{j=1}\frac{|||h_n|||^{\left|\vec{k}_j\right|}_{\left|\vec{l}_j\right|}}{\vec{k}_j ! \cdot \left(\vec{l}_j !\right)^{\left|\vec{k}_j\right|}}
\end{align*}
and compute using Lemma \ref{lem:normh}: $\prod^{s}_{j=1} ||| h_n |||^{\left| \vec{k}_j \right|}_{\left|\vec{l}_j\right|} \leq \hat{C} \cdot q^{3 \cdot \left(m-1\right)^2 \cdot k \cdot n \cdot \left(n+1\right)}_{n}$, where $\hat{C}$ is a constant independent of $q_n$. Since $H_{n-1}$ is independent of $q_n$ we conclude: 
\begin{equation*}
\left\|D_{\vec{\nu}} \left[H_n\right]_r \right\|_0 \leq \check{C} \cdot q^{3 \cdot \left(m-1\right)^2 \cdot k \cdot n \cdot \left(n+1\right)}_{n}, 
\end{equation*}
where $\check{C}$ is a constant independent of $q_n$. \\
In the same way we prove an analogous estimate of $\left\|D_{\vec{\nu}} \left[H^{-1}_n\right]_r \right\|_0$ and verify the claim.
\end{pr}

In particular, we see that this norm can be estimated by a power of $q_n$. 

\subsection{Proof of convergence}
For the proof of the convergence of the sequence $\left(f_n\right)_{n\in \mathbb{N}}$ in the Diff$^{\infty}\left(M\right)$-topology the next result, that can be found in \cite[Lemma 4]{FSW}, is very useful.

\begin{lem} \label{lem:konj}
Let $k \in \mathbb{N}_0$ and $h$ be a C$^{\infty}$-diffeomorphism on $M$. Then we get for every $\alpha,\beta \in \mathbb{R}$:
\begin{equation*}
d_k\left(h \circ R_{\alpha} \circ h^{-1}, h \circ R_{\beta} \circ h^{-1}\right) \leq C_k \cdot ||| h |||^{k+1}_{k+1} \cdot \left| \alpha - \beta \right|,
\end{equation*}
where the constant $C_k$ depends solely on $k$ and $m$. In particular $C_0 = 1$.
\end{lem}

In the following Lemma we show that under some assumptions on the sequence $\left(\alpha_n\right)_{n\in \mathbb{N}}$ the sequence $\left(f_n\right)_{n \in \mathbb{N}}$ converges to $f \in \mathcal{A}_{\alpha}\left(M\right)$ in the Diff$^{\infty}\left(M\right)$-topology. Afterwards, we will show that we can fulfil these conditions (see Lemma \ref{lem:conv}).

\begin{lem} \label{lem:convgen}
Let $\varepsilon > 0$ be arbitrary and $\left(k_n\right)_{n \in \mathbb{N}}$ be a strictly increasing sequence of natural numbers satisfying $\sum^{\infty}_{n=1}\frac{1}{k_n} < \varepsilon$. Furthermore, we assume that in our constructions the following conditions are fulfilled:
$$
\left| \alpha - \alpha_1 \right| < \varepsilon \quad\text{ and   }\quad\left| \alpha - \alpha_n \right| \leq \frac{1}{2 \cdot k_n \cdot C_{k_n} \cdot ||| H_n |||^{k_n+1}_{k_n+1}}\text{ for every }n \in \mathbb{N},
$$
where $C_{k_n}$ are the constants from Lemma \ref{lem:konj}.
\begin{enumerate}
	\item Then the sequence of diffeomorphisms $f_n = H_n \circ R_{\alpha_{n+1}} \circ H^{-1}_{n}$ converges in the Diff$^{\infty}(M)$-topology to a measure-preserving smooth diffeomorphism $f$, for which $d_{\infty}\left(f,R_{\alpha}\right)< 3 \cdot \varepsilon$ holds.
	\item Also the sequence of diffeomorphisms $\hat{f}_n = H_n \circ R_{\alpha} \circ H^{-1}_{n} \in \mathcal{A}_{\alpha}\left(M\right)$ converges to $f$ in the Diff$^{\infty}(M)$-topology. Hence $f \in \mathcal{A}_{\alpha}\left(M\right)$. 
\end{enumerate}
\end{lem}

\begin{pr}
\begin{enumerate}
	\item According to our construction it holds $h_n \circ R_{\alpha_n} = R_{\alpha_n} \circ h_n$ and hence
	\begin{align*}
	f_{n-1} & = H_{n-1} \circ R_{\alpha_n} \circ H^{-1}_{n-1} = H_{n-1} \circ R_{\alpha_n} \circ h_n \circ h^{-1}_{n} \circ H^{-1}_{n-1} \\
	& = H_{n-1} \circ h_n \circ R_{\alpha_n} \circ h^{-1}_{n} \circ H^{-1}_{n-1} = H_n \circ R_{\alpha_n} \circ H^{-1}_{n}.
	\end{align*}
	Applying Lemma \ref{lem:konj} we obtain for every $k,n \in \mathbb{N}$:
	\begin{equation}  \label{est1}
	d_k\left(f_n, f_{n-1}\right) = d_k\left(H_n \circ R_{\alpha_{n+1}} \circ H^{-1}_{n}, H_n \circ R_{\alpha_n} \circ H^{-1}_{n}\right) \leq C_k \cdot ||| H_n |||^{k+1}_{k+1} \cdot \left| \alpha_{n+1} - \alpha_n \right|
	\end{equation}
	In section \ref{subsection:first} we assumed $\left|\alpha - \alpha_n \right| \stackrel{n \rightarrow \infty}{\longrightarrow} 0$ monotonically. Using the triangle inequality we obtain $\left| \alpha_{n+1} - \alpha_n \right| \leq \left| \alpha_{n+1} - \alpha \right| + \left| \alpha - \alpha_n \right| \leq 2 \cdot \left| \alpha - \alpha_n \right|$ and therefore equation (\ref{est1}) becomes:
	\begin{equation*}
	d_k\left(f_n, f_{n-1}\right) \leq C_k \cdot ||| H_n |||^{k+1}_{k+1} \cdot 2 \cdot \left| \alpha_{n} - \alpha \right|.
	\end{equation*}
	By the assumptions of this Lemma it follows for every $k \leq k_n$:
	\begin{equation} \label{est4}
	d_k\left(f_n,f_{n-1}\right) \leq d_{k_n}\left(f_n,f_{n-1}\right) \leq C_{k_n} \cdot ||| H_n |||^{k_n+1}_{k_n +1} \cdot 2 \cdot \frac{1}{2 \cdot k_n \cdot C_{k_n} \cdot ||| H_n |||^{k_n+1}_{k_n +1} } \leq \frac{1}{k_n}
	\end{equation}
	In the next step we show that for arbitrary $k \in \mathbb{N}$ $\left(f_n\right)_{n \in \mathbb{N}}$ is a Cauchy sequence in Diff$^k\left(M\right)$, i.e. $\lim_{n,m\rightarrow \infty} d_k\left(f_n,f_m\right) = 0$. For this purpose, we calculate:
	\begin{equation} \label{est2}
	\lim_{n \rightarrow \infty} d_k\left(f_n,f_m\right) \leq \lim_{n \rightarrow \infty} \sum^{n}_{i=m+1} d_k\left(f_i, f_{i-1}\right) = \sum^{\infty}_{i=m+1}d_k\left(f_i, f_{i-1}\right).
	\end{equation}
	We consider the limit process $ m \rightarrow \infty$, i.e. we can assume $k\leq k_m$ and obtain from equations (\ref{est4}) and (\ref{est2}):
	\begin{equation*}
	\lim_{n,m \rightarrow \infty} d_k\left(f_n,f_m\right) \leq \lim_{m \rightarrow \infty} \sum^{\infty}_{i=m+1} \frac{1}{k_i} = 0.
	\end{equation*}
	Since Diff$^{k}\left(M\right)$ is complete, the sequence $\left(f_n\right)_{n \in \mathbb{N}}$ converges consequently in Diff$^k\left(M\right)$ for every $k \in \mathbb{N}$. Thus, the sequence converges in Diff$^{\infty}\left(M\right)$ by definition. \\
	\\
	Furthermore, we estimate:
	\begin{equation} \label{est3}
	d_{\infty}\left(R_{\alpha},f \right)= d_{\infty}\left(R_{\alpha}, \lim_{n\rightarrow \infty} f_n\right)
 \leq d_{\infty}\left(R_{\alpha}, R_{\alpha_1}\right) + \sum^{\infty}_{n=1} d_{\infty}\left(f_n, f_{n-1}\right),  
	\end{equation} 
	where we used the notation $f_0=R_{\alpha_1}$. \\
	By explicit calculations we obtain $d_k\left(R_{\alpha}, R_{\alpha_1}\right) = d_0\left(R_{\alpha}, R_{\alpha_1}\right) = \left| \alpha - \alpha_1 \right|$ for every $k \in \mathbb{N}$, hence 
	\begin{equation*}
	d_{\infty}\left(R_{\alpha}, R_{\alpha_1}\right) = \sum^{\infty}_{k=1} \frac{\left| \alpha - \alpha_1 \right|}{2^k \cdot \left(1+d_k\left(R_{\alpha}, R_{\alpha_1}\right)\right)} \leq \left| \alpha - \alpha_1 \right| \cdot \sum^{\infty}_{k=1} \frac{1}{2^k} = \left| \alpha - \alpha_1 \right|.
	\end{equation*}
	Additionally it holds:
	\begin{align*}
	\sum^{\infty}_{n=1} d_{\infty}\left(f_n, f_{n-1}\right) & = \sum^{\infty}_{n=1} \sum^{\infty}_{k=1} \frac{d_k\left(f_n, f_{n-1}\right)}{2^k \cdot \left(1+d_k\left(f_n,f_{n-1}\right)\right)} \\
	&= \sum^{\infty}_{n=1}\left( \sum^{k_n}_{k=1}\frac{d_k\left(f_n, f_{n-1}\right)}{2^k \cdot \left(1+d_k\left(f_n,f_{n-1}\right)\right)} + \sum^{\infty}_{k=k_n +1}\frac{d_k\left(f_n, f_{n-1}\right)}{2^k \cdot \left(1+d_k\left(f_n,f_{n-1}\right)\right)}\right)
	\end{align*}
	As seen above $d_k\left(f_n,f_{n-1}\right)\leq \frac{1}{k_n}$ for every $k\leq k_n$. Hereby, it follows further:
	\begin{align*}
	\sum^{\infty}_{n=1} d_{\infty}\left(f_n, f_{n-1}\right) & \leq \sum^{\infty}_{n=1}\left( \frac{1}{k_n} \cdot \sum^{k_n}_{k=1}\frac{1}{2^k} + \sum^{\infty}_{k=k_n +1}\frac{d_k\left(f_n, f_{n-1}\right)}{2^k \cdot \left(1+d_k\left(f_n,f_{n-1}\right)\right)}\right) \\
	& \leq \sum^{\infty}_{n=1} \frac{1}{k_n} + \sum^{\infty}_{n=1}\sum^{\infty}_{k=k_n +1} \frac{1}{2^k}.
	\end{align*}
	Because of $\sum^{\infty}_{k=k_n +1} \frac{1}{2^k} = 2 - \sum^{k_n}_{k=0} \frac{1}{2^k} = \left(\frac{1}{2}\right)^{k_n} \leq \frac{1}{k_n}$ we conclude:
	\begin{equation*}
	\sum^{\infty}_{n=1} d_{\infty}\left(f_n, f_{n-1}\right)\leq \sum^{\infty}_{n=1} \frac{1}{k_n} + \sum^{\infty}_{n=1}\frac{1}{k_n} < 2 \cdot \varepsilon.
	\end{equation*}
	Hence, using equation (\ref{est3}) we obtain the desired estimate $d_{\infty}\left(f, R_{\alpha}\right) < 3 \cdot \varepsilon$.
	\item We have to show: $\hat{f}_n \rightarrow f$ in Diff$^{\infty}\left(M\right)$. \\
	For it we compute with the aid of Lemma \ref{lem:konj} for every $n \in \mathbb{N}$ and $k \leq k_n$:
	\begin{align*}
	d_k\left(f_n, \hat{f}_n\right) & \leq d_{k_n}\left(H_n \circ R_{\alpha_{n+1}} \circ H^{-1}_n, H_n \circ R_{\alpha} \circ H^{-1}_n\right) \\
	& \leq C_{k_n} \cdot ||| H_n |||^{k_n +1}_{k_n+1} \cdot \left| \alpha_{n+1} - \alpha \right| \leq C_{k_n} \cdot ||| H_n |||^{k_n +1}_{k_n+1} \cdot \left| \alpha_{n} - \alpha \right| \\
	& \leq C_{k_n} \cdot ||| H_n |||^{k_n +1}_{k_n+1} \cdot \frac{1}{2 \cdot k_n \cdot C_{k_n} \cdot ||| H_n |||^{k_n +1}_{k_n+1}} = \frac{1}{2 \cdot k_n} \leq \frac{1}{k_n}.
	\end{align*}
	Fix some $k \in \mathbb{N}$. \\
	\textbf{Claim: }$\forall \delta >0 \ \ \exists N \ \  \forall n \geq N: \ \ d_k\left(f, \hat{f}_{n}\right) < \delta$, i.e. $\hat{f}_n \rightarrow f$ in Diff$^k\left(M\right)$.\\
	\textbf{Proof: } Let $\delta >0$ be given. Since $f_n \rightarrow f$ in Diff$^{\infty}\left(M\right)$ we have  $f_n \rightarrow f$ in Diff$^{k}\left(M\right)$ in particular. Hence, there is $n_1 \in \mathbb{N}$, such that $d_k\left(f,f_n\right)< \frac{\delta}{2}$ for every $n \geq n_1$. Because of $k_n \rightarrow \infty$ we conclude the existence of $n_2 \in \mathbb{N}$, such that $\frac{1}{k_n}<\frac{\delta}{2}$ for every $n\geq n_2$, as well as the existence of $n_3 \in \mathbb{N}$, such that $k_n\geq k$ for every $n\geq n_3$. Then we obtain for every $n\geq \max\left\{n_1,n_2,n_3\right\}$:
	\begin{equation*}
	d_k\left(f, \hat{f}_n\right) \leq d_k\left(f,f_n\right) + d_k\left(f_n, \hat{f}_n\right) < \frac{\delta}{2} + d_{k_n}\left(f_n, \hat{f}_n\right) \leq \frac{\delta}{2} + \frac{1}{k_n} < \frac{\delta}{2} + \frac{\delta}{2} = \delta.
	\end{equation*}
	Hence, the claim is proven.\\
	\\
	In the next step we show: $\lim_{n\rightarrow \infty}d_{\infty}\left(\hat{f}_n,f\right) = 0$. For this purpose, we examine:
	\begin{align*}
	d_{\infty}\left(f_n, \hat{f}_n\right) & = \sum^{k_n}_{k=1} \frac{d_k\left(f_n, \hat{f}_n\right)}{2^k \cdot \left(1+d_k\left(f_n, \hat{f}_n\right)\right)} + \sum^{\infty}_{k=k_n +1} \frac{d_k\left(f_n, \hat{f}_n\right)}{2^k \cdot \left(1+d_k\left(f_n, \hat{f}_n\right)\right)} \\
	& \leq \frac{1}{k_n} \cdot \sum^{k_n}_{k=1} \frac{1}{2^k} + \sum^{\infty}_{k=k_n +1} \frac{1}{2^k} \leq \frac{1}{k_n} + \left(\frac{1}{2}\right)^{k_n}.
	\end{align*}
	Consequently $\lim_{n \rightarrow \infty} d_{\infty}\left(f_n, \hat{f}_n\right) = 0$. With it we compute:	
	\begin{align*}
	\lim_{n \rightarrow \infty} d_{\infty}\left(f, \hat{f}_n\right) & = \lim_{n \rightarrow \infty} d_{\infty}\left(\lim_{m \rightarrow \infty}f_m, \hat{f}_n\right) = \lim_{n \rightarrow \infty} \lim_{m \rightarrow \infty} d_{\infty}\left(f_m, \hat{f}_n\right) \\
	&\leq \lim_{n \rightarrow \infty} \lim_{m \rightarrow \infty} \left( \sum^{m}_{i=n+1} d_{\infty}\left(f_i,f_{i-1}\right) + d_{\infty}\left(f_n, \hat{f}_n\right)\right) \\
	&= \lim_{n \rightarrow \infty}  \sum^{\infty}_{i=n+1} d_{\infty}\left(f_i,f_{i-1}\right)  + \lim_{n \rightarrow \infty} d_{\infty}\left(f_n, \hat{f}_n\right) = 0.
	\end{align*}
	As asserted we obtain: $\lim_{n\rightarrow \infty}d_{\infty}\left(\hat{f}_n,f\right) = 0$.
\end{enumerate}
\end{pr}

As announced we show that we can satisfy the conditions from Lemma \ref{lem:convgen} in our constructions:

\begin{lem} \label{lem:conv}
Let $\left(k_n\right)_{n \in \mathbb{N}}$ be a strictly increasing sequence of natural numbers with $\sum^{\infty}_{n=1} \frac{1}{k_n} < \infty$ and $C_{k_n}$ be the constants from Lemma \ref{lem:konj}. For any Liouvillean number $\alpha$ there exists a sequence $\alpha_n = \frac{p_n}{q_n}$ of rational numbers with the property that $260n^4$ divides $q_n$, such that our conjugation maps $H_n$ constructed in section \ref{subsection:g} and \ref{subsection:phi} fulfil the following conditions:
\begin{enumerate}
	\item For every $n \in \mathbb{N}$: 
	\begin{equation*}
	\left| \alpha - \alpha_n \right| < \frac{1}{2 \cdot k_n \cdot C_{k_n} \cdot ||| H_n |||^{k_n +1}_{k_n +1}}.
	\end{equation*}
	\item For every $n \in \mathbb{N}$: 
	\begin{equation*}
	\left| \alpha - \alpha_n \right| < \frac{1}{2^{n+1} \cdot q_n \cdot ||| H_n |||_1}.
	\end{equation*}
	\item For every $n \in \mathbb{N}$:
	\begin{equation*}
	\left\|DH_{n-1} \right\|_0 < \frac{\ln\left(q_n\right)}{n}.
	\end{equation*}
\end{enumerate}
\end{lem}

\begin{pr}
In Lemma \ref{lem:normH} we saw $||| H_n |||_{k_n+1} \leq \breve{C}_n \cdot q^{3 \cdot \left(m-1\right)^2 \cdot \left(k_n+1\right) \cdot n \cdot \left(n+1\right)}_{n}$, where the constant $\breve{C}_n$ was independent of $q_n$. Thus, we can choose $q_n \geq \breve{C}_n$ for every $n \in \mathbb{N}$. Hence, we obtain: $||| H_n |||_{k_n+1} \leq q^{4 \cdot \left(m-1\right)^2 \cdot \left(k_n+1\right) \cdot n \cdot \left(n+1\right)}_{n}$. \\
Besides $q_n \geq \breve{C}_n$ we keep the mentioned condition $q_n \geq 64 \cdot 260 \cdot n^4 \cdot (n-1)^{11} \cdot q^{(m-1) \cdot (n-1)^2+3}_{n-1}$ in mind. Furthermore, we can demand $\left\|DH_{n-1}\right\|_0 < \frac{\ln\left(q_n\right)}{n}$ from $q_n$ because $H_{n-1}$ is independent of $q_n$. Since $\alpha$ is a Liouvillean number, we find a sequence of rational numbers $\tilde{\alpha}_n = \frac{\tilde{p}_n}{\tilde{q}_n}$, $\tilde{p}_n, \tilde{q}_n$ relatively prime, under the above restrictions (formulated for $\tilde{q}_n$) satisfying:
\begin{equation*}
  \left| \alpha - \tilde{\alpha}_n \right| = \left| \alpha - \frac{\tilde{p}_n}{\tilde{q}_n} \right|< \frac{\left|\alpha-\alpha_{n-1}\right|}{2^{n+1} \cdot k_n \cdot C_{k_n} \cdot \left(260n^4\right)^{1+4 \cdot \left(m-1\right)^2 \cdot \left(k_n+1\right)^2 \cdot n \cdot \left(n+1\right)} \cdot \tilde{q}^{1+4 \cdot \left(m-1\right)^2 \cdot \left(k_n+1\right)^2 \cdot n \cdot \left(n+1\right)}_{n}}
\end{equation*}
Put $q_n \coloneqq 260 n^4 \cdot \tilde{q}_n$ and $p_n \coloneqq 260n^4 \cdot \tilde{p}_n$. Then we obtain: 
\begin{equation*}
\left| \alpha - \alpha_n \right| < \frac{\left|\alpha-\alpha_{n-1}\right|}{2^{n+1} \cdot k_n \cdot C_{k_n} \cdot q^{1+4 \cdot \left(m-1\right)^2 \cdot \left(k_n+1\right)^2 \cdot n \cdot \left(n+1\right)}_{n}}.
\end{equation*}
So we have $\left|\alpha-\alpha_{n}\right| \stackrel{n\rightarrow \infty}{\longrightarrow} 0$ monotonically. Because of $||| H_n |||^{k_n+1}_{k_n+1} \leq q^{4 \cdot \left(m-1\right)^2 \cdot \left(k_n+1\right)^2 \cdot n \cdot \left(n+1\right)}$ this yields: $\left| \alpha - \alpha_n \right| < \frac{1}{2^{n+1} \cdot q_n \cdot k_n \cdot C_{k_n} \cdot ||| H_n |||^{k_n+1}_{k_n +1}}$. Thus, the first property of this Lemma is fulfilled. \\
Furthermore, we note $k_n\geq1$ and $C_{k_n} \geq 1$ by Lemma \ref{lem:konj}. Thus, $q_n \cdot k_n \cdot C_{k_n} \geq q_n$. Moreover, $||| H_n |||_1 \geq \left\|H_n\right\|_0 =1$, because $H_n: \mathbb{S}^1 \times \left[0,1\right]^{m-1} \rightarrow \mathbb{S}^1 \times \left[0,1\right]^{m-1}$ is a diffeomorphism.  Hence, $||| H_n |||^{k_n+1}_{k_n +1}\geq ||| H_n |||_1$. Altogether, we conclude $2^{n+1} \cdot q_n \cdot k_n \cdot C_{k_n} \cdot ||| H_n |||^{k_n+1}_{k_n +1} \geq 2^{n+1} \cdot q_n \cdot ||| H_n |||_1$ and so: 
\begin{equation} \label{alpha}
\left| \alpha - \alpha_n \right| < \frac{1}{2^{n+1} \cdot q_n \cdot k_n \cdot C_{k_n} \cdot ||| H_n |||^{k_n+1}_{k_n +1}} \leq \frac{1}{2^{n+1} \cdot q_n \cdot ||| H_n |||_1},
\end{equation}
i.e. we verified the second property.
\end{pr}
\begin{rem}
Lemma \ref{lem:conv} shows that the conditions of Lemma \ref{lem:convgen} are satisfied. Therefore, our sequence of constructed diffeomorphisms $f_n$ converges in the Diff$^{\infty}(M)$-topology to a diffeomorphism $f \in \mathcal{A}_{\alpha}(M)$.
\end{rem}

To apply Proposition \ref{prop:crit} we need another result:

\begin{lem} \label{lem:m}
Let $\left(\alpha_n \right)_{n \in \mathbb{N}}$ be constructed as in Lemma \ref{lem:conv}. Then it holds for every $n \in \mathbb{N}$ and for every $\tilde{m} \leq q_{n+1}$:
\begin{equation*}
d_0\left(f^{\tilde{m}},f^{\tilde{m}}_{n}\right)\leq\frac{1}{2^n}.
\end{equation*}
\end{lem}

\begin{pr}
In the proof of Lemma \ref{lem:convgen} we observed $f_{i-1}= H_{i} \circ R_{\alpha_{i}} \circ H^{-1}_{i}$ for every $i \in \mathbb{N}$. Hereby and with the help of Lemma \ref{lem:konj} we compute:
\begin{equation*}
d_0\left(f^{\tilde{m}}_{i},f^{\tilde{m}}_{i-1}\right) = d_0\left(H_{i} \circ R_{\tilde{m} \cdot \alpha_{i+1}} \circ H^{-1}_{i},  H_{i} \circ R_{\tilde{m} \cdot \alpha_{i}} \circ H^{-1}_{i}\right) \leq ||| H_{i} |||_1 \cdot \tilde{m} \cdot 2 \cdot \left|\alpha-\alpha_{i}\right|.
\end{equation*}
Since $\tilde{m} \leq q_{n+1} \leq q_i$ we conclude for every $i > n$ using equation (\ref{alpha}) :
\begin{equation*}
d_0 \left( f^{\tilde{m}}_{i}, f^{\tilde{m}}_{i-1}\right) \leq ||| H_i |||_1 \cdot \tilde{m} \cdot 2 \cdot \left| \alpha - \alpha_i \right| \leq ||| H_i |||_1 \cdot \tilde{m} \cdot 2 \cdot \frac{1}{2^{i+1} \cdot q_i \cdot ||| H_i |||_1} \leq \frac{\tilde{m}}{q_i} \cdot \frac{1}{2^i} \leq \frac{1}{2^i}.
\end{equation*}
Thus, for every $\tilde{m}\leq q_{n+1}$ we get the claimed result:
\begin{equation*}
d_0\left(f^{\tilde{m}}, f^{\tilde{m}}_{n}\right) = \lim_{k \rightarrow \infty} d_0\left(f^{\tilde{m}}_{k}, f^{\tilde{m}}_{n}\right) \leq \lim_{k \rightarrow \infty} \sum^{k}_{i=n+1} d_0\left(f^{\tilde{m}}_{i}, f^{\tilde{m}}_{i-1}\right)\leq  \sum^{\infty}_{i=n+1} \frac{1}{2^{i}} = \left(\frac{1}{2}\right)^n.
\end{equation*}
\end{pr}

\begin{rem}
Note that the sequence $\left(m_n\right)_{n \in \mathbb{N}}$ defined in section \ref{section:distri} meets the mentioned condition $m_{n} \leq q_{n+1}$ and hence Lemma \ref{lem:m} can be applied to it.
\end{rem}

Concluding we have checked that all the assumptions of Proposition \ref{prop:crit} are satisfied. Thus, this criterion guarantees that the constructed diffeomorphism $f \in \mathcal{A}_{\alpha}(M)$ is weakly mixing. In addition, for every $\varepsilon >0$ we can choose the parameters by Lemma \ref{lem:convgen} in such a way, that $d_{\infty}\left(f, R_{\alpha}\right) < \varepsilon$ holds.  

\section{Construction of the measurable $f$-invariant Riemannian metric} \label{section:metric}
Let $\omega_0$ denote the standard Riemannian metric on $M = \mathbb{S}^1 \times \left[0,1\right]^{m-1}$. The following Lemma shows that the conjugation map $h_n = g_n \circ \phi_n$ constructed in section \ref{section:constr} is an isometry with respect to $\omega_0$ on the elements of the partial partition $\zeta_n$.

\begin{lem} \label{lem:isom}
Let $\check{I}_n \in \zeta_n$. Then $h_n |_{\check{I}_n}$ is an isometry with respect to $\omega_0$.
\end{lem}

\begin{pr}
Let $\check{I}_{n,k} \in \zeta_n$ be a partition element on $\left[\frac{k-1}{n \cdot q_n}, \frac{k}{n \cdot q_n}\right] \times \left[0,1\right]^{m-1}$. This element $\check{I}_{n,k}$ is positioned in such a way that all the occurring maps $\varphi_{\varepsilon, 1,j}$ and $\varphi^{-1}_{\varepsilon_2,1,j}$ act as rotations on it. Thus, $\phi_n |_{\check{I}_{n,k}}$ is an isometry and $\phi_n\left(\check{I}_{n,k}\right)$ is equal to
\begin{equation*}
\begin{split}
& \Bigg{[}\frac{k-1}{n \cdot q_n} + \frac{j^{(1)}_{1}}{n \cdot q^2_n}+...+ \frac{j^{\left((m-1) \cdot \frac{k \cdot \left(k-1\right)}{2}\right)}_{1}+1}{n \cdot q^{(m-1) \cdot \frac{k \cdot \left(k-1\right)}{2}+1}_n}-\frac{j^{(1)}_2}{n \cdot q^{(m-1) \cdot \frac{k \cdot \left(k-1\right)}{2}+2}_n}-...-\frac{j^{(k)}_2}{n \cdot q^{(m-1) \cdot \frac{k \cdot \left(k-1\right)}{2}+k+1}_n}\\
& \quad - \frac{j^{(1)}_3}{n \cdot q^{(m-1) \cdot \frac{k \cdot \left(k-1\right)}{2}+k+2}_n}-...-\frac{j^{(k)}_m+1}{n \cdot q^{(m-1) \cdot \frac{k \cdot \left(k+1\right)}{2}+1}_n}+\frac{1}{n^5 \cdot q^{(m-1) \cdot \frac{k \cdot \left(k+1\right)}{2}+1}_n},\\
& \quad \frac{k-1}{n \cdot q_n} + \frac{j^{(1)}_{1}}{n \cdot q^2_n}+...-\frac{j^{(k)}_m}{n \cdot q^{(m-1) \cdot \frac{k \cdot \left(k+1\right)}{2}+1}_n}-\frac{1}{n^5 \cdot q^{(m-1) \cdot \frac{k \cdot \left(k+1\right)}{2}+1}_n} \Bigg{]} \\
\times & \Bigg{[} \frac{j^{\left((m-1) \cdot \frac{k \cdot \left(k-1\right)}{2}+1\right)}_1}{q_n}+ ...\frac{j^{\left((m-1) \cdot \frac{k \cdot \left(k-1\right)}{2}+k\right)}_1}{q^k_n}+\frac{j^{(k+1)}_2}{q^{k+1}_n} +...+ \frac{j^{\left((m-1) \cdot \frac{k \cdot \left(k+1\right)}{2}+1\right)}_2}{q^{1+(m-1) \cdot \frac{k \cdot \left(k+1\right)}{2}}_n}+ \\
& \quad \frac{j^{\left((m-1) \cdot \frac{k \cdot \left(k+1\right)}{2}+2\right)}_2 }{8 n^5 \cdot q^{1+(m-1) \cdot \frac{k \cdot \left(k+1\right)}{2}}_n \cdot \left[n q^{\sigma}_{n}\right]} +  \frac{1}{8  n^9 \cdot q^{1+(m-1) \cdot \frac{k \cdot \left(k+1\right)}{2}}_n \cdot \left[ n q^{\sigma}_{n}\right]}, \\
& \quad \frac{j^{\left((m-1) \cdot \frac{k \cdot \left(k-1\right)}{2}+1\right)}_1}{q_n}+ ...+ \frac{j^{\left((m-1) \cdot \frac{k \cdot \left(k+1\right)}{2}+2\right)}_2 + 1}{8 n^5 \cdot q^{1+(m-1) \cdot \frac{k \cdot \left(k+1\right)}{2}}_n \cdot \left[n q^{\sigma}_{n}\right]} -  \frac{1}{8  n^9 \cdot q^{1+(m-1) \cdot \frac{k \cdot \left(k+1\right)}{2}}_n \cdot \left[ n q^{\sigma}_{n}\right]}\Bigg{]} \\
\times & \prod^{m}_{i=3} \Bigg{[} \frac{j^{\left((m-1) \cdot \frac{k \cdot \left(k-1\right)}{2}+(i-2) \cdot k+1\right)}_1}{q_n} + ... + \frac{j^{\left((m-1) \cdot \frac{k \cdot \left(k-1\right)}{2}+(i-1) \cdot k\right)}_1}{q^k_n} + \frac{1}{n^4 \cdot q^{k}_n}, \\
& \quad \quad \frac{j^{\left((m-1) \cdot \frac{k \cdot \left(k-1\right)}{2}+(i-2) \cdot k+1\right)}_1}{q_n} + ... + \frac{j^{\left((m-1) \cdot \frac{k \cdot \left(k-1\right)}{2}+(i-1) \cdot k\right)}_1 + 1}{q^k_n} - \frac{1}{n^4 \cdot q^{k}_n}\Bigg{]}. 
\end{split}
\end{equation*}
Then we have to examine the application of $g_n = g_{n \cdot q^{1+(m-1) \cdot \frac{\left(k+1\right) \cdot k}{2}}_n, \left[n \cdot q^{\sigma}_{n}\right], \frac{1}{8n^4}, \frac{1}{32n^4}}$. In particular, we have $\frac{\varepsilon}{b \cdot a} = \frac{1}{8n^4 \cdot \left[n \cdot q^{\sigma}_{n}\right] \cdot n \cdot q^{1+(m-1) \cdot \frac{\left(k+1\right) \cdot k}{2}}_n}$. Since $4 \cdot \varepsilon = \frac{1}{2n^4}< \frac{1}{n^4}$, $g_n$ acts as translation on $\phi_n\left(\check{I}_{n,k}\right)$.
\end{pr}

This Lemma implies that $h^{-1}_{n} |_{h_n\left(\check{I}_n\right)}$ is an isometry as well.\\
In the following we construct the $f$-invariant measurable Riemannian metric. This construction parallels the approach in \cite[section 4.8]{GK}. For it we put $\omega_n \coloneqq \left(H^{-1}_n\right)^{\ast} \omega_0$. Each $\omega_n$ is a smooth Riemannian metric because it is the pullback of a smooth metric via a $C^{\infty}\left(M\right)$-diffeomorphism. Since $R^{\ast}_{\alpha_{n+1}}\omega_0 = \omega_0$ the metric $\omega_n$ is $f_n$-invariant:
\begin{align*}
f^{\ast}_{n} \omega_n & = \left(H_n \circ R_{\alpha_{n+1}} \circ H^{-1}_{n}\right)^{\ast} \left(H^{-1}_{n}\right)^{\ast}\omega_0 = \left(H^{-1}_{n}\right)^{\ast} R^{\ast}_{\alpha_{n+1}}  H^{\ast}_{n} \left(H^{-1}_{n}\right)^{\ast} \omega_0 = \left(H^{-1}_{n}\right)^{\ast} R^{\ast}_{\alpha_{n+1}} \omega_0 \\
& = \left(H^{-1}_{n}\right)^{\ast} \omega_0 = \omega_n.
\end{align*}
With the succeeding Lemmas we show that the limit $\omega_{\infty} \coloneqq \lim_{n \rightarrow \infty} \omega_n$ exists $\mu$-almost everywhere and is the desired $f$-invariant Riemannian metric.

\begin{lem} \label{lem:ae}
The sequence $\left(\omega_n\right)_{n \in \mathbb{N}}$ converges $\mu$-a.e. to a limit $\omega_{\infty}$
\end{lem}

\begin{pr}
For every $N \in \mathbb{N}$ we have for every $k>0$:
\begin{equation*}
\omega_{N+k} = \left(H^{-1}_{N+k}\right)^{\ast} \omega_0 = \left(h^{-1}_{N+k} \circ ... \circ h^{-1}_{N+1} \circ H^{-1}_{N}\right)^{\ast} \omega_0 = \left(H^{-1}_{N}\right)^{\ast} \left(h^{-1}_{N+k} \circ ... \circ h^{-1}_{N+1}\right)^{\ast} \omega_0.
\end{equation*}
Since the elements of the partition $\zeta_n$ cover $M$ except a set of measure at most $\frac{4m}{n^2}$ by Remark \ref{rem:überd}, Lemma \ref{lem:isom} shows that $\omega_{N+k}$ coincides with $\omega_N = \left( H^{-1}_{N}\right)^{\ast} \omega_0$ on a set of measure at least $1-\sum^{\infty}_{n= N+1} \frac{4m}{n^2}$. As this measure approaches $1$ for $N\rightarrow \infty$, the sequence $\left(\omega_n\right)_{n \in \mathbb{N}}$ converges on a set of full measure.
\end{pr}

\begin{lem}
The limit $\omega_{\infty}$ is a measurable Riemannian metric.
\end{lem}

\begin{pr}
The limit $\omega_{\infty}$ is a measurable map because it is the pointwise limit of the smooth metrics $\omega_n$, which in particular are measurable. By the same reasoning $\omega_{\infty} |_p$ is symmetric for $\mu$-almost every $p \in M$. Furthermore, $\omega_{\infty}$ is positive definite because $\omega_n$ is positive definite for every $n \in \mathbb{N}$ and $\omega_{\infty}$ coincides with $\omega_N$ on $T_1M \otimes T_1M$ minus a set of measure at most $\sum^{\infty}_{n= N+1} \frac{4m}{n^2}$. Since this is true for every $N \in \mathbb{N}$, $\omega_{\infty}$ is positive definite on a set of full measure. 
\end{pr}

\begin{rem} \label{rem:Egoroff}
In the proof of the subsequent Lemma we will need Egoroff's theorem (for example \cite[§21, Theorem A]{Ha2}): Let $\left(N,d\right)$ denote a separable metric space. Given a sequence $\left(\varphi_n\right)_{n \in \mathbb{N}}$ of $N$-valued measurable functions on a measure space $\left(X, \Sigma, \mu\right)$ and a measurable subset $A \subseteq X$, $\mu\left(A\right) < \infty$, such that $\left(\varphi_n\right)_{n \in \mathbb{N}}$ converges $\mu$-a.e. on $A$ to a limit function $\varphi$. Then for every $\varepsilon > 0$ there exists a measurable subset $B \subset A$ such that $\mu\left(B\right) < \varepsilon$ and $\left(\varphi_n\right)_{n \in \mathbb{N}}$ converges to $\varphi$ uniformly on $A \setminus B$.
\end{rem}

\begin{lem}
$\omega_{\infty}$ is $f$-invariant, i.e. $f^{\ast}\omega_{\infty} = \omega_{\infty}$ $\mu$-a.e..
\end{lem}

\begin{pr}
By Lemma \ref{lem:ae} the sequence $\left(\omega_n\right)_{n \in \mathbb{N}}$ converges in the C$^{\infty}$-topology pointwise almost everywhere. Hence, we obtain using Egoroff's theorem: For every $\delta > 0$ there is a set $C_{\delta} \subseteq M$ such that $\mu\left(M \setminus C_{\delta}\right) < \delta$ and the convergence $\omega_n \rightarrow \omega_{\infty}$ is uniform on $C_{\delta}$. \\
The function $f$ was constructed as the limit of the sequence $\left(f_n\right)_{n \in \mathbb{N}}$ in the C$^{\infty}$-topology. Thus, $\tilde{f}_n \coloneqq f^{-1}_{n} \circ f \rightarrow \text{id}$ in the C$^{\infty}$-topology. Since $M$ is compact, this convergence is uniform too. \\
Furthermore, the smoothness of $f$ implies $f^{\ast} \omega_{\infty} = f^{\ast}\lim_{n \rightarrow \infty} \omega_n = \lim_{n \rightarrow \infty} f^{\ast} \omega_n$. Therewith, we compute on $C_{\delta}$:
$f^{\ast} \omega_{\infty} = \lim_{n \rightarrow \infty} \left( \left(f_n \tilde{f}_n \right)^{\ast} \omega_n \right) = \lim_{n \rightarrow \infty} \left( \tilde{f}^{\ast}_{n} f^{\ast}_{n} \omega_n \right) = \lim_{n \rightarrow \infty} \tilde{f}^{\ast}_{n} \omega_n = \omega_{\infty}$, where we used the uniform convergence on $C_{\delta}$ in the last step. As this holds on every set $C_{\delta}$ with $\delta > 0$, it also holds on the set $\bigcup_{\delta > 0} C_{\delta}$. This is a set of full measure and therefore the claim follows.
\end{pr}

Hence, the desired $f$-invariant measurable Riemannian metric $\omega_{\infty}$ is constructed and thus Proposition \ref{prop:satz2} is proven.

\end{document}